
\documentstyle{amsppt}
\magnification=1200
\loadeufm

\font\dotless=cmr10 

\def\umi{{\"{\dotless\char'020}}}

\def\hati{{\^{\dotless\char'020}}}
\def\Hati{\hbox{$\hat{\imath}$}}
\def\Hatj{\hbox{$\hat{\jmath}$}}
\font\boldtitlefont=cmssdc10 scaled\magstep2
\font\authorfont=cmss10 scaled\magstep1
\font\sectionfont=cmss10 scaled\magstep1
\font\sans=cmss10

\NoRunningHeads
\pagewidth{6.5 true in}
\pageheight{8.9 true in}
\loadeusm

\catcode`\@=11
\def\logo@{}
\catcode`\@=13

\newif\iffirstul

\def\ullspecial#1{\ifx\empty#1 \let\next=\relax \else
    \iffirstul\firstulfalse\else%
    \discretionary{}{}{\hbox{\ }}\fi%
$\underline{\vphantom{y}%
  \hbox{#1}}$\let\next=\ullspecial\fi 
 \next}

\def\nspace{\lineskip=1pt\baselineskip=12pt%
     \lineskiplimit=0pt}
\def\dspace{\lineskip=2pt\baselineskip=18pt%
     \lineskiplimit=0pt}

\def\smbul{{\sssize\bullet}}
\def\oplusop{\operatornamewithlimits{\oplus}\limits}
\def\otimesop{\operatornamewithlimits{\otimes}
     \limits}
\def\colim{\operatornamewithlimits{\text{\rm colim}}}
\def\w{{\mathchoice{\,{\scriptstyle\wedge}\,}
{{\scriptstyle\wedge}}{{\scriptscriptstyle\wedge}}%
     {{\scriptscriptstyle\wedge}}}}
\def\Le{{\mathchoice{\,{\scriptstyle\le}\,}
{\,{\scriptstyle\le}\,}{\,{\scriptscriptstyle\le}\,}%
{\,{\scriptscriptstyle\le}\,}}}
\def\Ge{{\mathchoice{\,{\scriptstyle\ge}\,}
{\,{\scriptstyle\ge}\,}{\,{\scriptscriptstyle\ge}\,}%
{\,{\scriptscriptstyle\ge}\,}}}
\def\tworightarrows#1#2{\lower2.0pt\hbox{\,\,\vbox{\baselineskip=0pt
    \hbox{\lower6.0pt\hbox to #1pt{\rightarrowfill}}\kern-6pt
    \hbox to #2pt{\rightarrowfill}}\,\,}}
\def\mapright#1{\smash{\mathop{\,\longrightarrow\,}%
  \limits^{#1}}}
\def\longmapright#1#2{\,\smash{\mathop
     {\hbox to #1pt{\rightarrowfill}}\limits^{#2}}\,}
\def\longmaprightsub#1#2{\,\smash{\mathop
     {\hbox to #1pt{\rightarrowfill}}\limits_{#2}}\,}
\def\simover{\smash{\mathop{\,\longrightarrow\,}
  \limits^{\lower1.5pt\hbox{$\scriptstyle\sim$}}}}
\def\longhook{\lhook\!\longrightarrow}
\def\Righthook#1{\scriptscriptstyle%
  \subset\kern-4.00pt\lower2.85pt
     \hbox to #1pt{\rightarrowfill}}
\def\Heading#1#2{\bigskip\par\noindent{\bf#1}
     \medskip\parindent=#2pt\nspace}
\def\Ref[#1]{\par\smallskip\hang\indent\llap{\hbox to\parindent
     {[#1]\hfil\enspace}}\ignorespaces}
\def\ref#1{\par\smallskip\hang\indent\llap{\hbox to\parindent
     {#1\hfil\enspace}}\ignorespaces}

\def\Item#1{\par\smallskip\hang\indent
  \llap{\hbox to\parindent{#1\hfill\enspace}}
     \ignorespaces}
\def\Amalg{\,\shortparallel\kern-10.15 true pt
\hbox{\lower.5 true pt\hbox to 4.0 true pt{\hrulefill}}\,\,}
\def\upvee{^{\sssize\vee}}
\def\uphat{^{\sssize\wedge}}
\def\CUP{{\,\ssize \cup\,}}
\def\CAP{{\,\ssize \cap\,}}
\def\Coprod{{\ssize\coprod}}

\redefine\Supset{\raise1.75pt
     \hbox{$\,\,\ssize\supset\,\,$}}
\redefine\Subset{\raise1.75pt
     \hbox{$\,\,\ssize\subset\,\,$}}
\redefine\mod{\,\text{\rm mod\,}}
\redefine\Cup{\raise.80pt
     \hbox{$\,\ssize\cup\,$}}

\def\DR{\text{\rm DR}}
\def\Gal{\text{\rm Gal}}
\def\Ext{\text{\rm Ext}}
\def\End{\text{\rm End}}
\def\ab{\text{\rm ab}}
\def\ad{\text{\rm ad}}
\def\comp{\text{\rm comp}}
\def\Dec{\text{\rm Dec}}
\def\distr{\text{\rm distr}}
\def\est{\text{\rm est}}
\def\Spec{\text{\rm Spec}}
\def\Gr{\text{\rm Gr}} 
\def\GR{\text{\rm gr}}
\def\Sym{\text{\rm Sym}} 
\def\GL{\text{\rm GL}} 
\def\Hom{\text{\rm Hom}} 
\def\HOM{{\scrH}\text{\rm om}}
\def\image{\text{\rm image}}
\def\Ind{\text{\rm Ind}}

\def\Lie{\text{\rm Lie}}

\def\Lib{\text{\rm Lib}}

\def\ieme{\text{\rm i\`eme}}
\def\iemes{\text{\rm i\`emes}}
\def\im{\text{\rm Im}}
\def\ev{\text{\rm ev}}
\def\et{\text{\rm et}}
\def\Fil{\text{\rm Fil}}
\def\real{\text{\rm real}}
\def\mot{\text{\rm mot}}

\def\MAT{\text{\rm MAT}}

\def\MT{\text{\rm MT}}
\def\DM{\text{\rm DM}}
\def\DMT{\text{\rm DMT}}
\def\SmCor{\text{\rm SmCor}}
\def\Sing{\text{\rm Sing}}
\def\torsion{\text{\rm torsion}}
\def\pr{\text{\rm pr}}
\def\Ker{\text{\rm Ker}}
\def\kar{\text{\rm kar}}
\def\Aut{\text{\rm Aut}} 
\def\Rep{\text{\rm Rep}}

\def\dch{\text{\rm dch}}

\def\coeff{\text{\rm coeff}}

\def\eps{{\varepsilon}}

\def\abar{\bar{a}} 
\def\dbar{\bar{d}} 

\def\kbar{\bar{k}} 
\def\Kbar{\bar{K}}
\def\Xbar{\bar{X}}
\def\xbar{\bar{x}}
\def\Ybar{\bar{Y}}
\def\ybar{\bar{y}} 
\def\zbar{\bar{z}}
\def\Zbar{\bar{Z}}
\def\Gammabar{\bar{\Gamma}}

\def\zetabar{{\bar{\zeta}}}
\def\dbQbar{\bar{\dbQ}}
\def\underbark{\hbox{\b{$k$}}}
\def\underbart{\hbox{\b{$t$}}}

\def\khat{\Hat{k}}
\def\xhat{\Hat{x}}

\def\bfs{\bold{s}}
\def\bfv{\bold{v}}

\def\Etil{\widetilde{E}}
\def\Mtil{\widetilde{M}}
\def\uptil{\kern-1pt\tilde{\phantom{a}}}
\def\upcheck{^{\sssize\vee}}
\def\Lamuptil{\Lambda\raise1pt\hbox{$\!\uptil$}}

\def\dbC{{\Bbb C}} 
\def\dbG{{\Bbb G}} 
\def\dbH{{\Bbb H}} 
\def\dbR{{\Bbb R}} 
\def\dbP{{\Bbb P}}
\def\dbQ{{\Bbb Q}}
\def\dbZ{{\Bbb Z}}

\def\scr#1{{\fam\eusmfam\relax#1}}

\def\scrA{{\scr A}} 
\def\scrC{{\scr C}} 
\def\scrD{{\scr D}}
\def\scrE{{\scr E}}
\def\scrF{{\scr F}}
\def\scrH{{\scr H}}
\def\scrK{{\scr K}}
\def\scrI{{\scr I}}
\def\scrL{{\scr L}}

\def\scrO{{\scr O}}
\def\scrR{{\scr R}}
\def\scrS{{\scr S}}
\def\scrT{{\scr T}}
\def\scrU{{\scr U}}

\def\gr#1{{\fam\eufmfam\relax#1}}


\def\gru{{\gr u}} 
\def\grg{{\gr g}} 
\def\grl{{\gr l}}

\def\grZ{{\gr Z}}
\def\grgl{{\grg\grl}}

\parindent=20pt
\frenchspacing
\footline={\hfill}

\topmatter
\title\nofrills
{\boldtitlefont Groupes fondamentaux motiviques
de Tate mixte}
\endtitle
\bigskip
\author
{\authorfont P. Deligne et A.~B.~Goncharov$^*$}
\endauthor
\endtopmatter
\footnote""{$^*$Research partially supported by NSF
Grant DMS-0099390}

\NoBlackBoxes
\document

{\narrower{\narrower{\narrower{
\bigskip\bigskip\bigskip\noindent
{\sans R\'esum\'e.}\enspace
Nous d\'efinissons la cat\'egorie des motifs de
Tate mixte sur l'anneau des $S$-entiers d'un 
corps de nombres, et le groupe fondamental
motivique (rendu unipotent) d'une vari\'et\'e
unirationnelle sur un corps de nombres.
Nous consid\'erons plus en d\'etail le groupe
fondamental motivique de la droite projective
moins $0$, $\infty$ et les racines $N$-i\`emes
de l'unit\'e.\bigskip}}}}

\bigskip\bigskip
{\narrower{\narrower{\narrower{\bigskip\noindent
{\sans Abstract.}\enspace
We define the category of mixed Tate motives
over the ring of $S$-integers of a number
field.
We define the motivic fundamental group (made
unipotent) of a unirational variety over a
number field.
We apply this to the study of the motivic
fundamental group of the projective line
punctured at zero, infinity and all $N$-th
roots of unity.\bigskip}}}}

\newpage
\footline={\hss\folio\hss}
\pageno=-1

\def\leaderfill{\leaders\hbox to 1em{\hss.\hss}\hfill}

\topmatter
\toc\nofrills{\sectionfont Table of Contents}
\bigskip\bigskip
\widestnumber\head{2a.}
\head 0. Introduction\page{1}\endhead
\bigskip
\head 1. Motifs de Tate mixte sur l'anneau
     des $S$-entiers de $k$ \page{3}\endhead
\bigskip
\head 2. Le point de vue tannakien \page{13}\endhead
\bigskip
\head 3.  Cohomologie et groupe fondamental
  \page{25}\endhead
\bigskip
\head 4. Le groupe fondamental unipotent
motivique d'une\hfill\break
vari\'et\'e rationnelle\page{40}\endhead
\bigskip
\head 5. Le groupe fondamental motivique du
compl\'ement,\hfill\break
      dans $\dbP^1$, de $0$, $\infty$ et
     $\mu_N$\page{50}\endhead
\bigskip
\head 6. Profondeur\page{67}\endhead
\bigskip
\specialhead{} Appendice. Rappels sur les
groupes alg\'ebriques 
   unipotents \page{76}\endspecialhead
\bigskip
\specialhead{} Bibliographie\page{82}\endspecialhead
\endtoc

\endtopmatter

\newpage

\footline={\hss$\phantom{\displaystyle\sum}\atop\folio$\hss}
\pageno=1

\subhead
0. Introduction
\endsubhead

\dspace
Si $k$ est un corps de nombres, on dispose
d'une cat\'egorie tannakienne $\MT(k)$
des motifs de Tate mixte sur $k$, pour laquelle
les groupes d'extensions de $\dbQ(0)$
par $\dbQ(n)$ ont la relation d\'esir\'ee avec
les groupes de $K$-th\'eorie de $k$.
C'est le coeur d'une $t$-structure sur une
sous-cat\'egorie de la cat\'egorie
triangul\'ee motivique de Levine (1998) ou
Voevodsky (2000). 
La conjecture d'annulation de
Beilinson-Soul\'e est vraie sur $k$, et c'est
ce qui permet de d\'efinir la $t$-structure
requise.
Par Levine (1998) et Huber (2000), on
dispose sur cette cat\'egorie tannakienne de
foncteurs fibres ``r\'ealisation''
correspondant aux th\'eories de cohomologie
usuelles (sauf que la cohomologie cristalline
n'a jusqu'\`a pr\'esent pas \'et\'e
consid\'er\'ee).

Variante: soient $S$ un ensemble de places
finies de $k$ et $\scrO_S$ l'anneau des
$S$-entiers de $k$.
De fa\c{c}on peut-\^{e}tre artificielle, mais
commode, on peut d\'efinir la cat\'egorie
tannakienne $\MT(\scrO_S)$ de motifs de Tate
mixte sur $\scrO_S$ comme une
sous-cat\'egorie de $\MT(k)$.
Voir 1.6.

Autre variante: par descente d'une extension
de $k$ \`a $k$, on d\'efinit la cat\'egorie
tannakienne $\MAT(k)$ des motifs  d'Artin-Tate 
mixte, qui deviennent
de Tate mixte sur une extension finie de $k$.
Elle contient celle des motifs d'Artin
(repr\'esentations rationnelles de
$\Gal(\kbar/k)$).

Dans Deligne (1989), l'un de nous a
d\'efini, pour certaines vari\'et\'es
alg\'ebriques $X$ sur $k$, un ``syst\`eme de
r\'ealisations'' du groupe fondamental
$\pi_1(X,0)$ rendu unipotent: une alg\`ebre de
Hopf commutative dans une cat\'egorie de
$\Ind$-syst\`emes de r\'ealisations -- ou
plut\^{o}t son spectre.
Si $X$ est unirationelle, nous construisons un
groupe fondamental rendu unipotent motivique
$\pi_1(X,0)_{\mot}$: une alg\`ebre de Hopf
commutative dans la cat\'egorie 
$\Ind$-$\MAT(k)$, des $\Ind$-objets de
$\MAT(k)$ dont la
pr\'ec\'edente se d\'eduise par application du
foncteur r\'ealisation.
Ceci implique des bornes sup\'erieures sur la
taille de l'action de Galois sur le
compl\'et\'e $\ell$-adique du $\pi_1$ 
(cf. Hain, Matsumoto (2001)), et sur le degr\'e de 
transcendance de corps engendr\'es par des
p\'eriodes.

Nous consid\'erons aussi le cas de points base
``\`a l'infini''.
Nous le faisons \`a peu de frais, en utilisant
que le foncteur ``r\'ealisations'' est
pleinement fid\`ele et r\'efl\`ete les
sous-quotients: pour prouver qu'un syst\`eme
de r\'ealisations est motivique, i.e. provient
d'un objet de $\MAT(k)$, il suffit de
l'exhiber comme sous-quotient d'un autre
syst\`eme de r\'ealisations dont on sache
d\'ej\`a qu'il est motivique. 

Dans la fin de l'article, nous utilisons ces
constructions dans le cas o\`u $X$ est le
compl\'ement, dans $\dbP^1$, de $0$, $\infty$ et
du groupe $\mu_N$ des racines $N^{\iemes}$ de
l'unit\'e.
Nous obtenons des r\'esultats sur la
structure de l'action 
du groupe de  Galois motivique sur le $\pi_1$.
Pour $N=1$, ils impliquent les r\'esultats de
Terasoma (2002) donnant la partie ``borne
sup\'erieure'' de la conjecture de Zagier sur
le nombre de valeurs multiz\^{e}ta de poids
$w$ lin\'eairement ind\'ependantes sur
$\dbQ$.

\newpage

\dspace
\subhead
1. Motifs de Tate mixte sur l'anneau des
$S$-entiers d'un corps de nombres
\endsubhead

\subhead
1.1
\endsubhead
Soit $k$ un corps. 
Nombre de r\'esultats cit\'es d\'ependant de
la r\'esolution des singularit\'es, nous
supposons $k$ de caract\'eristique $0$.
Hanamura (1995), Levine (1998) et Voevodsky
(2000) ont chacun d\'efini une cat\'egorie
triangul\'ee de motifs sur $k$, et Levine
(1998) VI 2.5.5 construit une \'equivalence
entre sa cat\'egorie triangul\'ee et celle de
Voevodsky.
Cette derni\`ere sera pour nous la plus
commode.
Nous la noterons $\DM(k)$ et noterons
$\DM(k)_{\dbQ}$ celle qui s'en d\'eduit par
tensorisation avec $\dbQ$.
On dispose dans $\DM(k)$ d'objets de Tate
$\dbZ(n)$ ($n\in\dbZ$), dont les images dans
$\DM(k)_{\dbQ}$ seront not\'ees $\dbQ(n)$.
Seule nous importera la sous-cat\'egorie
triangul\'ee $\DMT(k)_{\dbQ}$ de
$\DM(k)_{\dbQ}$ engendr\'ee par les $\dbQ(n)$:
celle des ``extensions it\'er\'ees'' de
$\dbQ(n)[m]$.
Rappelons que $E$ est {\it extension} de $B$
par $A$ s'il existe un triangle distingu\'e
$A\to E\to B\to A[1]$.

On dispose dans $\DM(k)$ d'un produit
tensoriel associatif et commutatif \`a unit\'e
$\otimes$, compatible \`a la structure
triangul\'ee.
L'automorphisme de sym\'etrie de
$\dbZ(1)\otimes\dbZ(1)$ est l'identit\'e
(Voevodsky (2000) 2.1.5), $M\mapsto
M(1):=M\otimes \dbZ(1)$ est une \'equivalence
(loc. cit. 4.1.3), et $\dbZ(n)$ est
$\dbZ(1)^{\otimes n}$ ($n\in\dbZ$).
Le produit tensoriel $\otimes$ est rigide
(loc. cit. 4.3.7): existence pour tout objet
$M$ d'un dual $M^*$ muni de $\ev\colon\,
M^*\otimes M\to \dbZ(0)$ et $\delta\colon\,
\dbZ(0)\to M\otimes M^*$ tels que les compos\'es
$M\to M\otimes M^*\otimes M\to M$ et $M^*\to
M^*\otimes M\otimes M^*\to M^*$ soient l'identit\'e.

Rappelons, d'apr\`es Levine (1993), comment,
lorsque la conjecture d'annulation de
Beilinson Soul\'e est v\'erifi\'ee
pour $k$, on peut extraire de $\DMT(k)_{\dbQ}$
une cat\'egorie tannakienne $\MT(k)$ qui
m\'erite le nom de cat\'egorie des motifs de
Tate mixte sur $k$.

Ecrivons $\Hom^j(M,N)$ pour $\Hom(M,N[j])$.
Les $\Hom^*(\dbZ(a),\dbZ(b))$ ne d\'ependant
que de $i=b-a$.
Ils sont donn\'es par les groupes de Chow
sup\'erieurs convenablement num\'erot\'es de
$\Spec(k)$ (loc. cit. 4.2.9 pour
$X=\Spec(k)$).
Pour $i< 0$ ils sont nuls. 
Pour $i=0$,
$$
\align
&\matrix
\Hom^j(\dbZ(0),\dbZ(0))
   &=\dbZ &\qquad \text{ si $j=0$},\\
   &\hfill 0 &\qquad \text{ sinon}.\hfill
\endmatrix\tag1.1.1\\
\intertext{Pour $i=1$,}
&\matrix
\Hom^j(\dbZ(0),\dbZ(1))&=k^* 
   &\qquad\text{ si $j=1$,}\\
   &\,\,\,\,0 &\qquad \text{ sinon}.\hfill
\endmatrix\tag1.1.2
\endalign
$$

Apr\`es tensorisation avec $\dbQ$, ces groupes
sont encore donn\'ees par les sous-espaces
propres des op\'erations d'Adams dans les
groupes de $K$-th\'eorie de $k$ (Levine
(1998) II 3.6.6; r\'ef\'erence originale:
Bloch (1986) compl\'et\'e par Bloch (1994);
voir aussi Levine (1994)).
La conjecture d'annulation de
Beilinson-Soul\'e est que
$$
\Hom^j(\dbQ(0),\dbQ(i))=0\quad
\text{pour}\quad i>0\quad\text{et}\quad
j\Le 0.
\tag1.1.3
$$
Lorsque tel est le cas, Beilinson Bernstein
Deligne (1982) 1.3.14, appliqu\'e \`a la
sous-cat\'egorie pleine de $\DMT(k)_{\dbQ}$
d'objets les $\dbQ(n)$, montre que les
extensions it\'er\'ees de $\dbQ(n)$
($n\in\dbZ$) forment une cat\'egorie
ab\'elienne, coeur d'une $t$-structure sur
$\DMT(k)_{\dbQ}$.
On l'appelle la cat\'egorie
$\MT(k)$ des motifs de Tate mixte sur $k$.
Elle est stable par produit tensoriel et
passage au dual.
Ces foncteurs \'etant exacts, cela r\'esulte
de ce que $\dbQ(n)\otimes\dbQ(m)=\dbQ(n+m)$ et
que $\dbQ(n)^*=\dbQ(-n)$.
On obtient ainsi sur $\MT(k)$ un produit
tensoriel exact en chaque variable,
associatif, commutatif, \`a unit\'e et rigide.

De ce que $\MT(k)$ est le coeur d'une
$t$-structure sur $\DMT(k)_{\dbQ}$ r\'esulte
que pour $A$ et $B$ dans $\MT(k)$, on a
$$
\align
\Ext^1(A,B) &=
\Hom^1(A,B)\qquad\text{et}\tag1.1.4\\
\Ext^2(A,B) &\hookrightarrow \Hom^2(A,B).\tag1.1.5
\endalign
$$
En particulier,
$$
\Ext^1(\dbQ(n),\dbQ(m))=0\quad
\text{pour}\quad m\Le n.\tag1.1.6
$$
Pour $A$ dans $\MT(k)$, il r\'esulte de
(1.1.6), par induction sur une description de
$A$ comme extension it\'er\'ee, que $A$ admet
une unique filtration ``par le poids''
$W$, finie, croissante et index\'ee par les
entiers pairs, telle que
$$
\Gr_{-2n}^W(M):= W_{-2n}(M)/W_{-2(n+1)}(M)
$$
soit une somme de copies de $\dbQ(n)$.
Posons
$$
\omega_n(M):=\Hom(\dbQ(n),\Gr_{-2n}^W(M)).
$$
On a donc 
$$
\Gr_{-2n}^W(M)=\dbQ(n)\otimes_\dbQ
\omega_n(M).
$$
La filtration $W$ est fonctorielle, exacte et
compatible au produit tensoriel.
Parce que l'automorphisme de symm\'etrie de
$\dbQ(1)\otimes\dbQ(1)$ est l'identit\'e, on
dispose d'isomorphismes 
$\dbQ(n)\otimes\dbQ(m)\simover\dbQ(n+m)$
compatibles \`a l'associativit\'e et \`a la
commutativit\'e de $\otimes$.
Le foncteur exact $M\mapsto
\omega(M):=\oplus\omega_n(M)$ est donc
un $\otimes$-foncteur: c'est un foncteur fibre
et la cat\'egorie $MT(k)$ est tannakienne.

\subhead
1.2
\endsubhead
Pour $M$ un motif de Tate mixte, le
sous-quotient $W_{-2n}(M)/W_{-2(n+2)}(M)$ de
$M$ est une extension $M_{n,n+1}$ de
$\Gr_{-2n}^W(M)=\dbQ(n)\otimes\omega_n(M)$ par
$\Gr_{-2(n+1)}^W(M)=\dbQ(n+1)\otimes\omega_{n+1}(M)$.
La classe de cette extension est un \'el\'ement
$$
e_n\in \Ext^1(\dbQ(n),\dbQ(n+1))\otimes\Hom
(\omega_n(M),\omega_{n+1}(M)).
$$
La tensorisation avec $\dbQ(n)$ est un
isomorphisme
$$
\Ext^1(\dbQ(0),\dbQ(1))\simover
\Ext^1(\dbQ(n),\dbQ(n+1))\tag1.2.1
$$
et les $e_n$ d\'efinissent
$$
e_M\colon\,\omega_*(M)\to\Ext^1(\dbQ(0),
\dbQ(1))\otimes\omega_*(M).
\tag1.2.2
$$

Pour $x\colon\,\dbQ(n)\to\Gr_{-2n}^W(M)$ dans
$\omega_n(M)$, l'image inverse par $x$ de
l'extension $M_{n,n+1}$ est une extension de
$\dbQ(n)$ par $\Gr_{-2(n+1)}^W(M)=\dbQ(n+1)
\otimes\omega_{n+1}(M)$, et $e_M(x)$ est la
classe de cette extension dans
$$
\Ext^1(\dbQ(0),\dbQ(1))\otimes
\omega_{n+1}(M)=\Ext^1(\dbQ(n),
\dbQ(n+1))\otimes\omega_{n+1}(M).
$$

On regardera $e_M$ comme une coaction, au sens
des coalg\`ebres de Lie.
Plus pr\'ecis\'ement, comme une coaction, sur
$\omega(M):=\oplus \omega_*(M)$, de la
coalg\`ebre de Lie librement coengendr\'ee par
$\Ext^1(\dbQ(0),\dbQ(1))$.

\proclaim{Proposition 1.3}
La coaction $e_M$ est fonctorielle en $M$ et
compatible au produit tensoriel.
\endproclaim

\demo{Preuve}
La v\'erification de la fonctorialit\'e en $M$
est laiss\'ee au lecteur.
Pour $x\in\omega_m(M)$ et $y\in\omega_n(N)$, on
doit montrer que
$$
e_{M\otimes N}(x\otimes y)=e_M(x)\otimes y
+x\otimes e_N(y).
\tag1.3.1
$$
Par fonctorialit\'e, on peut pour le
v\'erifier remplacer $M$ par $W_{-2m}(M)$, une
extension de $\dbQ(m)\otimes\omega_m(M)$ par
$W_{-2(m+1)}(M)$, puis par l'image inverse de
cette extension par
$x\colon\,\dbQ(m)\to\dbQ(m)\otimes\omega_m(M)$.
De m\^{e}me pour $N$, $n$ et $y$.
On se ram\`ene ainsi \`a supposer que
$M=W_{-2m}(M)$, $\omega_m(M)=\dbQ$, $x=1$, et
de m\^{e}me pour $N$, $n$ et $y$.
On a alors
$$
\align
\omega_{m+n+1}(M\otimes N)\simover
  &\omega_{m+n+1}(M\otimes\text{(quotient
$\dbQ(n)$ de $N$)})\\
\oplus\,\, &\omega_{m+n+1}\text{((quotient
$\dbQ(m)$ de $M$)}\otimes N).
\endalign
$$
et, par fonctorialit\'e en $M$ et $N$ de
$e_{M\otimes N}$, il suffit de v\'erifier
(1.3.1) pour $M$ remplac\'e par son quotient
$\dbQ(m)$, ainsi que pour $N$ remplac\'e par
son quotient $\dbQ(n)$.
Ces cas se r\'eduisent \`a la d\'efinition de
l'isomorphisme (1.2.1).
\enddemo

Si $E$ est une extension de $\dbQ(0)$ par
$\dbQ(1)$, il r\'esulte de 1.3 que la classe
de l'extension $E^*(1)$ de $\dbQ(0)$ par
$\dbQ(1)$ est l'oppos\'ee de la classe de $E$.

\definition{Definition 1.4}
{\it Pour $\Gamma$ un sous-espace vectoriel de
$\Ext^1(\dbQ(0),\dbQ(1))$, 
$MT(k)_\Gamma$ est la cat\'egorie des motifs de
Tate mixte tels que la coaction (1.2.2) se
factorise par une coaction de $\Gamma$}.
\enddefinition

Pour que $M$ soit dans $\MT(k)_\Gamma$, il
faut et il suffit que pour tout sous-quotient
$E$ de $M$ qui est une extension d'un
$\dbQ(n)$ par $\dbQ(n+1)$, la classe de $E$
dans 
$$
\Ext^1(\dbQ(0),\dbQ(1)\simover \Ext^1
(\dbQ(n),\dbQ(n+1))
$$ 
soit dans $\Gamma$. 
D'apr\`es 1.3, cette sous-cat\'egorie de
$MT(k)$, qui contient l'objet unit\'e $\dbQ(0)$,
est stable par sous-quotients, produits
tensoriels et duaux.
C'est une sous-cat\'egorie tannakienne de
$MT(k)$.

\subhead
1.5
\endsubhead
Quel que soit $k$ (de caract\'erisque $0$), on
dispose d'un foncteur ``r\'ealisation'' de
$\DM(k)$ dans la cat\'egorie triangul\'ee
``syst\`emes de r\'ealisations'' de Huber (1995).
L'\'etat de la litt\'erature \`a ce sujet
n'est pas satisfaisant, sp\'ecialement en ce qui
concerne la r\'ealisation de Hodge.
Dans Huber (2000), il est suppos\'e \`a tort
dans la preuve de 2.1.4, p.~775 que les
hyperrecouvrements propres de $X$ forment une
cat\'egorie filtrante.
Cette erreur est corrig\'ee dans un erratum
dat\'e d'octobre 2002.
Dans Levine (1998) V. 2.3.7 le foncteur Dec
de d\'ecalage d'une filtration par le poids
est appliqu\'e \`a un complexe de faisceaux,
alors qu'il devrait \^{e}tre appliqu\'e
apr\`es passage \`a $R\Gamma$.

Nous esquissons ci-dessous une fa\c{c}on
possible de proc\'eder, adapt\'ee \`a la
d\'efinition de Voevodsky
de $\DM(k)$, que nous commen\c{c}ons par
rappeler.
Le foncteur obtenu sera contravariant
(cohomologique).
Dans la suite de l'article, c'est sa variante
covariante (homologique), d\'eduite par
passage au dual, qui nous sera utile.

Soit $\SmCor(k)$ la cat\'egorie d'objets les
sch\'emas s\'epar\'es lisses sur $k$, un
morphisme de $X$ dans $Y$ \'etant une
correspondance finie de $X$ vers $Y$:
$\Hom(X,Y)$ est le groupe ab\'elien librement
engendr\'e par les sous-sch\'emas r\'eduits
irr\'eductibles ferm\'es $Z$ de $X\times Y$,
finis sur $X$ et dominant une composante
connexe de $X$.
La cat\'egorie $\SmCor(k)$ est additive
(avec $\oplus=\Coprod$).
La cat\'egorie triangul\'ee $\DM(k)$ se d\'eduit
de la cat\'egorie triangul\'ee de complexes
born\'es $K^b(\SmCor(k))$ par

\medskip\noindent
{\sans (a)} localisation:
on divise par la sous-cat\'egorie \'epaisse
engendr\'ee par les $[X\times A^1]\to[X]$
(invariance 
d'homotopie) et les $[U\CAP V]\to[U]\oplus[V]
\to[X]$ pour $X$ la r\'eunion des ouverts $U$
et $V$ (Mayer-Vietoris); 

\smallskip\noindent
{\sans (b)}\enspace 
adjonction de facteurs directs images
d'endomorphismes idempotents;

\smallskip\noindent
{\sans (c)}\enspace
inversion formelle de $\dbQ(1)$.

\medskip
Il r\'esulte de Voevodsky (2000) 
3.2.6 et de (a), Mayer
Vietoris, qu'on obtiendrait une cat\'egorie
\'equivalente si dans $\SmCor(k)$ on se
limitait aux sch\'emas lisses et
quasi-projectifs sur $k$.
Pour d\'efinir le foncteur r\'ealisation,
l'essentiel est de d\'efinir la r\'ealisation
d'un complexe born\'e $X^*$ d'objets de
$\SmCor(k)$, et par la remarque ci-dessus, on
peut se ramener au cas o\`u les $X^i$ sont
quasi-projectifs.

\proclaim{Lemme 1.5.1}
Soit $\Gamma\colon\, X\to Y$ dans $\SmCor(k)$.
Si $\Ybar$ est une compactification projective
et lisse de $Y$, il existe une
compactification projective et lisse $\Xbar$
de $X$ telle que $\Gamma$ se prolonge en
$\Gammabar\colon\,\Xbar\to\Ybar$.
quand elle existe, l'extension $\Gammabar$ 
est unique.
\endproclaim

\demo{Preuve}
Il suffit de traiter du cas o\`u $X$ est
connexe et o\`u $\Gamma$ est un sous-sch\'ema
ferm\'e $Z$ de $X\times Y$, r\'eduit,
irr\'eductible, fini et dominant sur $X$.
Si $d$ est le degr\'e de $Z$ sur $X$, $Z$ d\'efinit
un morphisme de sch\'emas $z$ de $X$ dans
$\Sym^d(Y)$.
D'apr\`es Hironaka, il existe une
compactification projective et lisse $\Xbar$
de $X$ telle que $z$ se prolonge en
$\zbar\colon\, \Xbar\to\Sym^d(\Ybar)$.
Pour un tel $\Xbar$, l'adh\'erence $\Zbar$ de
$Z$ dans $\Xbar\times\Ybar$ est encore finie
sur $\Xbar$.
Elle fournit le prolongement cherch\'e.
L'unicit\'e est claire.
Noter que, par Hironaka encore, on peut
supposer que $X$ est le compl\'ement dans
$\Xbar$ d'un diviseur \`a croisements normaux.

Pour $X^*$ un complexe born\'e d'objets de
$\SmCor(k)$, on d\'eduit de 1.5.1 l'existence
d'un syst\`eme de compactifications $\Xbar^n$
des $X^n$ tel que \ (a) $\Xbar^n$ est
projectif et lisse, et $X^n$ est le
compl\'ement dans $\Xbar^n$ d'un diviseur \`a
croisements normaux $D^n$; \ (b) la
diff\'erentielle $d$ se prolonge de $X^*$ \`a
$\Xbar^*$.
Par unicit\'e du prolongement, le prolongement
v\'erifie encore $d^2=0$.
Les syst\`emes de compactifications $\Xbar$
forment une cat\'egorie filtrante.
Pour $f\colon\, X^*\to Y^*$ un morphisme de
complexes, et pour $\Ybar^*$ une
compactification de $Y^*$ v\'erifiant (a),
(b), il existe une compactification $\Xbar^*$
de $X^*$ v\'erifiant (a), (b) et telle que $f$
se prolonge en un morphisme de complexes de
$\Xbar^*$ dans $\Ybar^*$.
\enddemo

D\'efinissons la r\'ealisation de de Rham de
$X^*$ -- ingr\'edient essentiel de la
r\'ealisation de Hodge relative \`a un
plongement complexe de $k$ dans $\dbC$.
L'id\'ee est de choisir une compactification
$\Xbar^*$ de $X^*$ comme ci-dessus et de
prendre sur chaque $\Xbar^n$ le complexe de de
Rham \`a p\^{o}les logarithmiques
$\Omega_{\Xbar^n}^*(\log\,D^n)$.
Ce complexe est bifiltr\'e: par le filtration
de Hodge $F^p=\sigma_{\Ge p}$, et par la
filtration par le poids $W_*$, qui compte le
nombre de facteurs $\frac{dz_i}{z_i}$ pour
$z_i$ une \'equation locale d'une composante
lisse de $D^n$.
Il s'agit ensuite d'obtenir des complexes
bifiltr\'es $(K^{(n)},W,F)$ repr\'esentant les
$R\Gamma(\Xbar^n,\Omega^*(\log\,D^n))$ et
d'utiliser la diff\'erentielle de $\Xbar^*$
pour en faire un complexe double bifiltr\'e,
avec $K^{(n)}$ comme colonne de premier indice
$-n$.
La difficult\'e sera de d\'eduire de $d^2=0$
pour $\Xbar^*$ que $d'{^2}=0$ identiquement.
Ceci fait, il reste \`a

\medskip\noindent
{\sans (i)}\enspace
remplacer la filtration $W$ de $K^{(n)}$ par
son translat\'e $W'$ d\'efini par
$W'{^m}=W^{m-n}$;

\smallskip\noindent
{\sans (ii)}\enspace
passer au complexe simple associ\'e; on notera
que son $\Gr^W$ est la somme des\break
$\Gr^W(K^{(n)})[n]$;

\smallskip\noindent
{\sans (iii)}\enspace
remplacer la filtration $W$ par sa d\'ecal\'ee
$\Dec\,W$.

\medskip
Nous travaillerons avec la topologie \'etale.
Celle de Nisnevitch ferait aussi bien
l'affaire, mais non celle de Zariski.
Pour un faisceau quasi-coh\'erent, la
cohomologie est la m\^{e}me.
L'hypercohomologie d'un complexe de de Rham
aussi.
Gain: pour $p\colon\, Z\to X$ un morphisme
fini, le foncteur $p_*$ est exact.
Ceci est faux pour Zariski.
Soit $\Gamma\colon\, X\to Y$ dans $\SmCor(k)$:
$\Gamma=\sum n_i Z_i$ et $Z_i\Subset X\times
Y$ est fini sur $X$.
Soient des compactifications $X$ et $Y$ en $\Xbar$ et
$\Ybar$, avec \`a l'infini un diviseur \`a
croisement normaux, et supposons que $\Gamma$
se prolonge en $\Gammabar$, i.e. que chaque
adh\'erence $\Zbar_i$ de $Z_i$ dans
$\Xbar\times\Ybar$ soit finie sur $\Xbar$.
Soit enfin $\vert\Gammabar\vert$ le support
$\cup \Zbar_i$ de $\Gammabar$.
La correspondance $\Gammabar$ d\'efinit alors
un morphisme de complexes de faisceaux sur
$\vert\Gamma\vert$
$$
[\Gammabar]\colon\, \pr_2^*\Omega^*
(\log)\to\pr_1^*\Omega^*(\log).
$$
Ici aussi que la topologie soit \'etale ou de
Nisnevitch est essentiel.
Prendre garde que $\pr_2^*$ est une image
inverse au sens des faisceaux, et non pas au
sens des faisceaux coh\'erents (pour laquelle
l'image inverse de la diff\'erentielle de de
Rham $d$ ne serait pas d\'efinie).
Ce morphisme s'\'etend aux r\'esolutions
flasques canoniques $\scrC^*$ et il suffit de
repr\'esenter
$R\Gamma(\Xbar^n,\Omega^*(\log))$ par
$\Gamma(\Xbar^n,\scrC^*\Omega^*(\log))$.

\subhead
1.6
\endsubhead
En caract\'eristique z\'ero, les seuls
corps pour lesquels la conjecture d'annulation
de Beilinson-Soul\'e soit d\'emontr\'ee sont
les corps de nombres ($=$ extensions finies de
$\dbQ$), les corps de fonctions rationelles
d'une courbe de genre $0$ sur un corps de
nombres, et les limites inductives de tels
corps.
Le cas qui nous importe est celui des corps de
nombres.
Pour ceux-ci, $K_n(k)\otimes\dbQ$ est nul pour
$n$ pair $>0$ et on a
$$
\align
\Ext^1(\dbQ(0),\dbQ(n))
&=K_{2n-1}(k)\otimes\dbQ,\quad\text{et}\\
\Ext^2(\dbQ(0),\dbQ(n)) &=0,
\endalign
$$
car $\Hom^2(\dbQ(0),\dbQ(n))=0$ dans
$\DM(k)_{\dbQ}$.

Soit $k$ un corps de nombres.
Le foncteur ``r\'ealisation'' sur $\DM(k)$
induit sur $\MT(k)$ un $\otimes$-foncteur \`a
valeurs dans les syst\`emes de r\'ealisations
au sens de Deligne (1989) (aspect cristallin
exclus) ou Jannsen (1990).

Rappelons que
$\Ext^1(\dbQ(0),\dbQ(1))=k^*\otimes\dbQ$.
Si $S$ est un ensemble de places finies de
$k$, et que $\scrO_S$ est l'ensemble des
$S$-entiers de $k$ (int\'egralit\'e en dehors
de $S$), nous d\'efinissons artificiellement la
cat\'egorie des motifs de Tate mixte sur
$\scrO_S$ comme \'etant $\MT(k)_\Gamma$
pour $\Gamma=\scrO_S^*\otimes\dbQ\Subset
k^*\otimes\dbQ^*$.
Cas particulier: pour $S$ le compl\'ement
d'une place finie $v$, $\scrO_S$
est le localis\'e $\scrO_{(v)}$ en $v$ de
l'anneau des entiers $\scrO$ de $k$, et on
parle de motifs de Tate mixte {\it non ramifi\'es}
en $v$.
Un motif de Tate mixte sur $\scrO_S$ est donc
un motif de Tate mixte sur $k$ non ramifi\'e
en chaque place finie $v\notin S$.

\proclaim{Proposition 1.7}
Soit $\ell$ un nombre premier distinct de la
caract\'eristique r\'esiduelle de $v$.
Pour qu'un motif de Tate mixte $M$ sur $k$ soit
non ramifi\'e en $v$, il faut et il suffit que
sa r\'ealisation $\ell$-adique $M_\ell$
soit non ramifi\'ee en $v$.
\endproclaim

\demo{Preuve}
La r\'ealisation $\ell$-adique est un
$\dbQ_\ell$-espace vectoriel d\'ependant
fonctoriellement du choix d'une cl\^{o}ture
alg\'ebrique $\kbar$ de $k$.
La fonctorialit\'e en $\kbar$ en fait une
repr\'esentation continue de $\Gal(\kbar/k)$.
Nous utiliserons que 

\medskip\noindent
{\sans (a)}
la r\'ealisation l-adique de $\dbQ(1)$
est $\dbQ_\ell(1)
:=\dbZ_\ell(1)\otimes_{\dbZ_\ell}\dbQ_\ell$,
avec $\dbZ_\ell(1)=\lim\,\mu_{\ell^n}(\kbar)$;

\medskip\noindent
{\sans (b)}
Pour $n\not=m$, puisque $\Hom(\dbQ(n),\dbQ(m))=0$,
une extension de $\dbQ(n)$ par $\dbQ(m)$ est
d\'etermin\'ee \`a isomorphisme unique pr\`es
par sa classe dans $\Ext^1(\dbQ(n),\dbQ(m))$.
Pour $x\in k^*\Subset\Ext^1(\dbQ(0),\dbQ(1))$,
notons $K(x)$ l'extension de $\dbQ(0)$
par $\dbQ(1)$, dite de Kummer, de classe $x$.
Sa r\'ealisation l-adique est d\'eduite,
par tensorisation avec $\dbQ_\ell$, de
l'extension  de Kummer $K_\ell(x)$ de
$\dbZ_\ell$ par $\dbZ_\ell(1)$.
La donn\'ee d'une extension
$$
0\to\dbZ_\ell(1)\to E\mapright{f}
\dbZ_\ell\to 0
$$
de  $\dbZ_\ell$ par $\dbZ_\ell(1)$ \'equivaut
\`a celle du $\dbZ_\ell(1)$-torseur
$f^{-1}(1)$, et $K_\ell(x)$ correspond au
$\dbZ_\ell(1)$ torseur limite projective des
$\mu_{\ell^n}(\kbar)$-torseurs des racines
$(\ell^n)^{\iemes}$ de $x$ dans $\kbar$.

\medskip
La repr\'esentation $\ell$-adique $K_\ell(x)$
est non ramifi\'ee en $v$ si et seulement si
$x$ est une unit\'e en $v$, i.e. est dans
$\scrO_{(v)}^*$.
Pour $M$ de Tate mixte, il r\'esulte donc de
(b) que $M$ est dans $MT(\scrO_v)$ si et
seulement si les repr\'esentations
$\ell$-adiques
$W_{-2n}(M_\ell)/W_{-2(n+2)}(M_\ell)$ sont non
ramifi\'ees.
Il reste \`a v\'erifier le lemme suivant.

\proclaim{Lemme 1.8}
Soient $K_v$ le corps des fractions de
l'hens\'elis\'e $\scrO_{(v)}^h$ de $\scrO_v$
et $H$ une repr\'e-\break
sentation $\ell$-adique de
$\Gal(\Kbar_v/K_v)$, munie d'une filtration
finie croissante index\'ee par les entiers
pairs $W$, telle que $\Gr_{-2n}^W(H)$ soit
somme de copies de $\dbQ_\ell(n)$.
Pour que la repr\'esentation $H$ soit non
ramifi\'ee, il suffit que les
repr\'esentations $W_{-2n}/W_{-2(n+2)}$ le
soient.
\endproclaim

\demo{Preuve}
Soit $I_v\Subset \Gal(\Kbar_v/K_v)$ le groupe
d'inertie.
On sait que son plus grand pro-$\ell$-quotient
est canoniquement isomorphe \`a
$\dbZ_\ell(1)$.
Soit $t_\ell\colon\, I_v\to\dbZ_\ell(1)$ le
morphisme de passage au quotient.
Il est $\Gal(\Kbar_v/K_v)$-\'equivariant.
\enddemo

Supposons que $I_v$ agisse trivialement sur
les $W_{-2n}/W_{-2(n+2)}$ et prouvons par
r\'e-\break
currence sur $r\Ge 2$ qu'il agit
trivialement sur les $W_{-2n}/W_{-2(n+r)}$.
Pour $r=2$, c'est l'hypo-\break
th\`ese.
Pour $r>2$, l'hypoth\`ese  de r\'ecurrence
assure que l'action de $I_v$ est triviale sur
$W_{-2n}/W_{-2(n+r-1)}$ et sur
$W_{-2(n+1)}/W_{-2(n+r)}$.
L'action de $\sigma\in I_v$ est donc de la
forme $1+\nu(\sigma)$, avec pour $\nu(\sigma)$
un compos\'e
$$
W_{-2n}/W_{-2(n+r)}\to
\Gr_{-2n}^W\longmapright{35}{n(\sigma)}
\Gr_{-2(n+r-1)}^W\to W_{-2n}/W_{-2(n+r)},
$$
et $n(\sigma_1\sigma_2)=n(\sigma_1)+n(\sigma_2)$.
L'application $\sigma\mapsto n(\sigma)$
se factorise donc par $\dbZ_\ell(1)$, et est
de la forme
$$
n(\sigma)=n \cdot t_\ell(\sigma),\,\,
\text{\rm pour }n\colon\,\Gr_{-2n}^W(1)\to
\Gr_{-2(n+r-1)}^W.
$$
Le morphisme $n$ est
$\Gal(\Kbar_v/K_v)$-equivariant.
Puisque $r>2$,
$\Hom_{\Gal(\Kbar_v/K_v)}(\dbQ_\ell(n+1~),$\break
$\dbQ_\ell(n+r-1))=0$ et $n$ est nul, ce qui
conclut la r\'ecurrence.
\enddemo

\proclaim{Proposition 1.9}
Soit $\Gamma$ un sous-espace vectoriel de
$\Ext^1(\dbQ(0),\dbQ(1))$.

\medskip\noindent
{\sans (i)}
Pour $r\Ge 2$, $\Ext^1(\dbQ(n),\dbQ(n+r))$ est
le m\^{e}me dans $MT(k)$ et dans
$MT(k)_\Gamma$.

\noindent
{\sans (ii)}
Les $\Ext^2$ de Yoneda sont nuls dans
$MT(k)_\Gamma$.
\endproclaim

\demo{Preuve}\quad (i) Si $r\Ge 2$, une
extension de $\dbQ(n)$ par $\dbQ(n+r)$ est en
effet automatiquement dans $MT(k)_\Gamma$.

\noindent
(ii)
Il suffit de v\'erifier la nullit\'e des
$\Ext^2(\dbQ(n),\dbQ(m))$.
Pour toute classe $c$ dans cet $\Ext^2$, il
existe dans $MT(k)_\Gamma$ une extension
$$
0\to F\to E_1\to\dbQ(n)\to 0,
\tag1.8.1
$$
de classe not\'ee $c_1$, et une extension
$E_2$ de $F$ par $\dbQ(m)$, de classe $c_2$,
telle que $c$ soit le produit de Yoneda
$c_2c_1$.
Ce produit est l'image de $c_2$ par le cobord
$$
\partial\colon\, \Ext^1(F,\dbQ(m))\to
\Ext^2(\dbQ(n),\dbQ(m))
$$
d\'efini par (1.8.1).
C'est aussi l'obstruction \`a l'existence dans
$MT(k)_\Gamma$ de\break
$X=X_0\Supset X_1\Supset
X_2\Supset X_3=0$  de quotients successifs
$X_i/X_{i+1}$ ($i=0,1,2$) munis\break
d'isomorphismes avec $\dbQ(n)$, $F$ et
$\dbQ(m)$, tel que la classe de l'extension
$X/X_2$ de $\dbQ(n)$ par $F$ soit $c_1$, et
que la classe de l'extension $X_1$ de $F$ par
$\dbQ(m)$ soit $c_2$.

Relevons $\dbQ(n)$ dans  $\Gr_{-2n}^W(E_1)$ et
soit $E'_1$ l'image inverse dans
$W_{-2n}(E_1)$ de ce rel\`evement.
C'est une extension de $\dbQ(n)$ par $F':=
F\CAP E'_1=W_{-2(n+1)}(F)$.
Si $E'_2$ est l'image inverse de $F'$ dans
$E_2$, $c$ est encore le produit de Yoneda
$c'_2c'_1$ de la classe $c'_1$ de $E'_1$ par
la classe $c'_2$ de $E'_2$.
On a gagn\'e que $F'=W_{-2(n+1)}(F')$.
Un argument dual permet ensuite de se ramener
au cas o\`u $F=W_{-2(n+1)}(F)$ et
$W_{-2m}(F)=0$.
Le $\Ext^2$ ne peut donc \^{e}tre non nul que
si $m\Ge n+2$.
Puisque le $\Ext^2$ est nul dans $MT(k)$, il
existe dans $MT(k)$ une extension it\'er\'ee
$X=X_0\Supset X_1\Supset X_2\Supset X_3=0$ du
type consid\'er\'e plus haut, le sous-quotient
$W_{-2a}/W_{-2(a+2)}$ de $X$ coincide avec le
m\^{e}me sous-quotient de $X_0/X_3$ si $a+2\le
m$, et avec le m\^{e}me sous-quotient de $X_1$
si $a>n$.
Si, comme on peut le supposer, $m\Ge n+2$, on
est pour tout a dans un de ces deux cas, $X$
est donc dans $MT(k)_\Gamma$, et $c=0$.
\enddemo

\newpage


\dspace
\subhead
2. Le point de vue tannakien 
\endsubhead

\subhead
2.1
\endsubhead
Soient $k$ un corps de nombres et $S$ un
ensemble fini de places finies de $k$.
Le foncteur $\omega$, d\'eduit du foncteur
gradu\'e $\omega_*$ de 1.1 par oubli de la graduation,
est un foncteur fibre de la cat\'egorie \`a
produit tensoriel $\MT(\scrO_S)$ d\'efinie en
1.6.
Cette cat\'egorie est donc tannakienne, et 
si $G_\omega$ est le sch\'ema en
groupe des automorphismes du
$\otimes$-foncteur $\omega$, le foncteur
$\omega$ induit une \'equivalence de
$\MT(\scrO_S)$ avec la cat\'egorie des
repr\'esentations (lin\'eaires, de dimension
finie) de $G_\omega$.

L'action de $G_\omega$ sur
$\omega(\dbQ(1))=\dbQ$ d\'efinit un morphisme
$$
G_\omega\to \dbG_m.\tag2.1.1
$$
Soit $U_\omega$ son noyau.
Le groupe multiplicatif $\dbG_m$ agit sur le
$\otimes$-foncteur $\omega$, $\lambda\in\dbG_m$
agissant sur $\omega_n$ par multiplication par
$\lambda^n$.
Cette action est une section $\dbG_m\to
G_\omega$ de (2.1.1).
Elle fait de $G_\omega$ un produit semi-direct
$$
G_\omega=\dbG_m\ltimes U_\omega.
\tag2.1.2
$$

Si cela est n\'ecessaire pour \'eviter une
ambigu{\umi}t\'e, on \'ecrira
$G_\omega\left<\MT(\scrO_S)\right>$ et\break
$U_\omega\left<\MT(\scrO_S)\right>$ au lieu de
$G_\omega$ et $U_\omega$.

L'action de $G_\omega$ est fonctorielle en $M$.
Pour tout $M$ dans $\MT(\scrO_S)$, 
elle respecte donc la 
filtration par le poids, index\'ee
par les entiers pairs, de $\omega(M)$ par les
$$
\omega(W_{-2n}M)=\oplusop_{m\Ge n}\omega_m(M).
$$

Le sous-groupe $U_\omega$ agit trivialement
sur $\omega(Q(1))$,  donc sur les
$\omega(\dbQ(n))$, ainsi que sur $\Gr^W(\omega(M))
=\omega(\Gr^W(M))$.
Il en r\'esulte que le sous-groupe
alg\'ebrique $U_\omega\left<M\right>$ de
$\GL(\omega(M))$ image de $U_\omega$ est
unipotent.
L'action de $\dbG_m\Subset G_\omega$ sur
$\omega(M)$ normalise cette image.
L'alg\`ebre de Lie $\gru_\omega\left<M\right>$
de $U_\omega\left<M\right>$ est donc une
sous-alg\`ebre de Lie gradu\'ee de
$\grgl(\omega(M))$.
Puisque l'action de $U_\omega\left<M\right>$
sur $\omega(M)$ respecte la filtration par le
poids et est l'identit\'e sur le gradu\'e
associ\'e, l'alg\`ebre de Lie
$\gru_\omega\left<M\right>$ est \`a
degr\'es $>0$.

Pour $M$ un motif de Tate mixte, notons
$\left<M\right>$ la sous-cat\'egorie
tannakienne de $\MT(\scrO_S)$
{\it engendr\'ee} par $M$, d'objets les
sous-quotients de sommes de $M^{\otimes
p}\otimes (M\upvee)^{\otimes q}$. 
Si $M'$ est dans $\left<M\right>$,
$U_\omega\left<M'\right>$ est un quotient de
$U_\omega\left<M\right>$, et le sch\'ema en
groupe $U_\omega$ est la limite projective
des $U_\omega\left<M\right>$ pour
$\left<M\right>$ de plus en plus grand.
Le sch\'ema en groupe $U_\omega$ est donc
pro-unipotent.

Les alg\`ebres de Lie
$\gru_\omega\left<M\right>$ sont \`a
degr\'e $>0$ et, par 1.6 et 1.7, les groupes
$$
\alignat2
\Ext^1(\dbQ(0),\dbQ(n)) = &\,0 
   &\qquad&\text{si $n\Le 0$},\tag2.1.3\\
&\scrO_S^*\otimes\dbQ &\qquad&\text{si $n=1$},\\
&K_{2n-1}(k)\otimes\dbQ &\qquad&\text{si $n\Ge 2$}.
\endalignat
$$
sont de dimension finie.
On est donc dans la situation consid\'er\'ee
en A.14.
En chaque degr\'e $n$, le syst\`eme projectif
des $(\gru_\omega\left<M\right>)_n$ est
stationnaire et on d\'efinit l'alg\`ebre de
Lie gradu\'ee $\gru_\omega^{\GR}$ comme
\'etant la somme sur $n$ des 
$\lim(\gru_\omega\left<M\right>)_n$.
L'alg\`ebre de $\Lie\,\gru_\omega$ de
$U_\omega$ est la limite projective des
$\gru_\omega\left<M\right>$; c'est le produit
des $(\gru_\omega^{\GR})_n$.
Appliquant A.15, on obtient:

\proclaim{Proposition 2.2}\quad
(i) L'alg\`ebre de Lie gradu\'ee
$\gru_\omega^{\GR}$ est \`a degr\'es $>0$ et est en
chaque degr\'e de dimension finie.
Elle d\'etermine $U_\omega$ par
$$
U_\omega
=\lim\,\exp(\gru_\omega^{\GR}\text{/partie de
degr\'e $\Ge N$}).
$$

\noindent
(ii)
Le foncteur $\omega$ induit une \'equivalence
de $\MT(\scrO_S)$ avec la cat\'egorie des
espaces vectoriels de dimension finie 
gradu\'es, munis d'une
action de $\gru_\omega^{\GR}$ compatible aux
graduations.
\endproclaim

Que les $\Ext^2(\dbQ(0),\dbQ(n))$ soient nuls
admet la traduction suivante.

\proclaim{Proposition 2.3}
L'alg\`ebre de Lie $\gru_\omega^{\GR}$ est libre.
\endproclaim

On en obtient un syst\`eme g\'en\'erateur
homog\`ene libre en relevant une base de 
chaque
$$
\left(\gru_\omega^{\GR}\right)_n^{\ab}=
\text{ dual de (2.1.3)}
\tag2.3.1
$$
dans $(\gru_\omega^{\GR})_n$.

\subhead
2.4 Mise en garde
\endsubhead
On ne dispose pas d'un rel\`evement
canonique de $(\gru_\omega^{\GR})_n^{\ab}$ dans
$(\gru_\omega^{\GR})_n$.
On ne dispose donc pas d'un isomorphisme
canonique de $\gru_\omega^{\GR}$ avec l'alg\`ebre
de Lie librement engendr\'ee par l'espace
vectoriel gradu\'e
$$
(\gru_\omega^{\GR})^{\ab}=\oplus(\Ext^1(\dbQ(0),
\dbQ(n))\upvee\text{ en degr\'e $n$}).
\tag2.4.2
$$

\subhead
2.5
\endsubhead
La cat\'egorie tannakienne $\MT(k)$ est la
r\'eunion filtrante des sous-cat\'egories
pleines $\MT(\scrO_S)$, pour $S$ de plus en
plus grand.
On a donc pour $\MT(k)$
$$
\align
G_\omega\left<\MT(k)\right>:=
\Aut^\otimes(\omega)
&=\lim\,G_\omega\left<\MT(\scrO_S)\right>=\\
&=\dbG_m\ltimes\lim\,U_\omega
\left<\MT(\scrO_S)\right>.
\endalign
$$

\subhead
2.6
\endsubhead
Rappelons quelques d\'efinitions fondant la
g\'eom\'etrie alg\'ebrique dans une
cat\'egorie tannakienne $\scrT$.

La cat\'egorie $\Ind\,\scrT$ des $\Ind$-objets
de $\scrT$ est munie d'un produit tensoriel
associatif commutatif \`a unit\'e h\'erit\'e
de $\scrT$.
Ceci permet de d\'efinir la notion
d'alg\`ebre commutative \`a unit\'e de $\Ind\,\scrT$:
c'est un objet $A$ muni d'un produit $A\otimes
A\to A$ et d'une unit\'e $1\to A$ v\'erifiant
les axiomes usuels.
On d\'efinit la cat\'egorie des
$\scrT$-{\it sch\'emas affines}, aussi
appel\'es sch\'emas affines en $\scrT$, comme
\'etant la duale de celle des alg\`ebres
commutatives \`a unit\'e de $\Ind\,\scrT$.
On note $\Spec(A)$ le $\scrT$-sch\'ema affine
correspondant \`a l'alg\`ebre $A$, et on
appelle $A$ son {\it alg\`ebre affine}.

Un $\scrT$-sch\'ema en groupe affine est un
objet en groupe de la cat\'egorie des
$\scrT$-sch\'emas affines, i.e. le spectre
$G=\Spec(A)$ d'une alg\`ebre de Hopf
commutative $A$ de $\Ind\,\scrT$.
Une action de $G$ sur un objet $M$ de $\scrT$
est une structure de comodule $M\to A\otimes M$.

Pour $M$ un object de $\scrT$, on notera encore
$M$ le $\scrT$-{\it sch\'ema vectoriel}
$\Spec(\Sym^*M\upvee)$, o\`u
$\Sym^*N:=\colim_A\oplusop\limits_{n\Le A}\Sym^n(N)$. 
Un pro-objet $\lim\,M_i$ de $\scrT$ d\'efinit
de m\^{e}me le $\scrT$-{\it sch\'ema provectoriel}
$\Spec(\colim_{i,A}\oplusop_{n\Le A}
\Sym^n(M_i\upvee))$.

Pour un expos\'e plus d\'etaill\'e, nous
renvoyons \`a Deligne (1989) \S5.
Par loc. cit. 5.8, de nombreux \'enonc\'es
connus en g\'eom\'etrie alg\'ebrique restent
valables en $\scrT$-g\'eom\'etrie alg\'ebrique.
Nous les utiliserons librement.

\subhead
2.7
\endsubhead
Rappelons Deligne (1990) qu'on dispose,
dans toute cat\'egorie tannakienne $\scrT$,
d'un $\scrT$-sch\'ema en groupe affine
canonique, le {\it groupe fondamental}
$\pi(\scrT)$ de $\scrT$.
Il agit sur tout objet $M$ de $\scrT$, cette
action est fonctorielle en $M$, et pour tout
foncteur fibre $F$ elle d\'efinit un
isomorphisme 
$$
F(\pi(\scrT))\simover \Aut^\otimes(F).\tag2.7.1
$$

\subhead
2.8
\endsubhead
Certaines des constructions 2.1--2.5 sont
motiviques, i.e. l'image par $\omega$ de
constructions dans $\MT(\scrO_S)$.
D'autres d\'ependent de ce que le foncteur
fibre $\omega$ est gradu\'e.

Notons $G$ le groupe fondamental $\pi(\MT(\scrO_S))$.
Les constructions suivantes sont motiviques:

\medskip\noindent
{\sans (a)}\enspace
le morphisme (2.1.1) de $G_\omega$ dans
$\dbG_m$.
C'est l'image par $\omega$ du morphisme
$$
G\to\dbG_m
\tag2.8.1
$$
donnant l'action du groupe fondamental $G$ de
$\MT(\scrO_S)$ sur $\dbQ(1)$.
Pour donner un sens \`a (2.8.1), observer que
la cat\'egorie des espace vectoriels de
dimension finie est une sous-cat\'egorie de
$\MT(\scrO_S)$, par $V\mapsto V\otimes\dbQ(0)$,
et qu'un sch\'ema en groupe affine peut donc
\^{e}tre consid\'er\'e comme un sch\'ema en
groupe de $\MT(\scrO_S)$.

\medskip\noindent
{\sans (b)}\enspace
le groupe $U_\omega$ est motivique: c'est
l'image par $\omega$ du noyau $U$ de (2.8.1).
De m\^{e}me pour sa pro-alg\`ebre de Lie
$\gru$

\medskip\noindent
{\sans (c)}\enspace
la formule (2.4.2) se rel\`eve en
un isomorphisme de pro-objets de
$\MT(\scrO_S)$
$$
\Lie(U^{\ab})=\prod\limits_{n}
\Ext^1(\dbQ(0),\dbQ(n))\upvee\otimes\dbQ(n).
\tag2.8.2
$$
Dans (2.8.2), la projecton sur le $n^{\ieme}$
facteur donne l'action de $\Lie\,U$ sur les
extensions de $\dbQ(0)$ par $\dbQ(n)$.

Par contre, la structure de produit
semi-direct (2.1.2) et la graduation de
$\gru_\omega$ (i.e., la d\'ecomposition de ce
provectoriel en produit des
$\left(\gru_\omega^{\GR}\right)_n$)
ne sont pas motiviques.  
La graduation de $\gru_\omega^{\ab}$ l'est,
image de (2.8.2) par $\omega$, car
l'action int\'erieure de $G$ sur $\gru^{\ab}$
se factorise par le quotient $\dbG_m$ de $G$.

Si cela est n\'ecessaire pour \'eviter une
ambiguit\'e, on \'ecrira
$G\left<\MT(\scrO_S)\right>$ et
$U\left<\MT(\scrO_S)\right>$ au lieu de $G$
et $U$.

\subhead
2.9
\endsubhead
La r\'ealisation de de Rham $M_{\DR}$ de $M$
dans $\MT(k)$ est un $k$-espace vectoriel, muni
d'une filtration par le poids $W$ image de
celle de $M$, et d'une filtration
d\'ecroissante $F$, la filtration de Hodge.
Ces filtrations sont fonctorielles en $M$,
exactes et compatibles au produit tensoriel.

On a
$$
\dbQ(1)_{\DR}=k\text{ en poids $-2$ et
filtration de Hodge $-1$},
$$
et $\dbQ(n)_{\DR}$ est donc $k$ en poids $-2n$
et filtration de Hodge $-n$.
Pour $M$ quelconque, les filtrations $W$ et
$F$ sont donc oppos\'ees, et d\'efinissent une
graduation de $M_{\DR}$ telle que
$$
\align
W_{-2n}(M_{\DR}) &=\oplusop_{m\Ge n}
  (M_{\DR})_m\\
F^{-n}(M_{\DR}) &=\oplusop_{m\Le n}
  (M_{\DR})_m.
\endalign
$$
On a
$$
\align
(M_{\DR})_n &=\Gr_{-2n}^W(M_{\DR})=
  (\Gr_{-2n}^W(M))_{\DR}\\
&=(\omega_n(M)\otimes\dbQ(n))_{\DR}=
  \omega_n(M)\otimes k.
\endalign
$$
On en d\'eduit la

\proclaim{Proposition 2.10}
Le foncteur fibre $\DR$ sur $\MT(k)$ est
canoniquement isomorphe au foncteur d\'eduit
du foncteur fibre $\omega$ par extension des
scalaires de $\dbQ$ \`a $k$.
\endproclaim

\subhead
2.11
\endsubhead
Pour $C$ une cl\^{o}ture alg\'ebrique de $\dbR$
et $\sigma$ un plongement de $k$ dans $C$, \`a
la th\'eorie de cohomologie
$$
X\text{ sur } k\longmapsto H^*\text{ (espace
topologique $X(C),\dbQ$)}
$$
correspond un foncteur fibre $M\mapsto
M_\sigma$, fonctoriel en $C$.
Parler de ``plongement dans $C$'' plut\^{o}t
que de plongement complexe (i.e. dans $\dbC$)
a l'avantage que la fonctorialit\'e en $C$
fournit un isomorphisme
$M_\sigma\to M_{\bar{\sigma}}$.
Pour un plongement r\'eel, i.e. pour 
$\sigma=\bar{\sigma}$, cet isomorphisme est
une involution de $M_\sigma$, le {\it Frobenius 
r\'eel}.

Un isomorphisme de comparaison canonique d'un
foncteur fibre $\alpha$ vers un foncteur
fibre $\beta$ sera not\'e
$\comp_{\beta,\alpha}$, et son inverse
$\comp_{\alpha,\beta}$.
La relation entre la cohomologie de de Rham
et la cohomologie singuli\`ere fournit
$$
\comp_{\sigma,\DR}\colon\,\,M_{\DR}\otimes_{k,\sigma}
C\simover M_\sigma\otimes_\dbQ C,
\tag2.11.1
$$
fonctoriel en $C$.
Il induit
$$
\comp_{\sigma,\omega}\colon\,\,
M_\omega\otimes C\simover
M_\sigma\otimes C\tag2.11.2
$$
et permet de voir $M_\sigma$ comme une
$\dbQ$-structure sur
$M_{\DR}\otimes_{k,\sigma}C=M_\omega\otimes
C$.
On a
$$
\dbQ(1)_\sigma=2\pi i\dbQ\Subset
C=\dbQ(1)_{\DR}\otimes_{k,\sigma}C
$$
et donc $\dbQ(n)_\sigma=(2\pi i)^n\dbQ\Subset
C$.
Dans ces formules, le choix de la racine
carr\'ee $i$ de $-1$ n'importe pas.

On reconstitue la $\dbQ$-structure $M_\omega$
de l'espace vectoriel gradu\'e $M_{\DR}$ \`a
partir de $M_\sigma$ et $\comp_{\sigma,\DR}$
par
$$
(M_\omega)_n\Subset(M_{\DR})_n\Subset
(M_{\DR}\otimes_{k,\sigma} C)_n\quad 
\est~\frac{1}{(2\pi i)^n}\,\comp_{\DR,\sigma}
\Gr_{-2n}^W(M_\sigma).
\tag2.11.3
$$

Ecrivons $G_\omega=\dbG_m\ltimes U_\omega$
pour $G_\omega\left<\MT(k)\right>$ et sa
d\'ecomposition 2.5, et $\tau$ pour
l'inclusion de $\dbG_m$ dans $G_\omega$.

\proclaim{Proposition 2.12}
Il existe $a_\sigma$ dans $G_\omega(C)$ tel que
$$
M_\sigma=\comp_{\sigma,\omega}(a_\sigma
M_\omega).\tag2.12.1
$$
Si $i$ est une racine carr\'ee de $-1$ dans
$C$, on peut choisir les $a_\sigma$ de la forme
$$
a_\sigma=a_\sigma^0\tau(2\pi i)\tag2.12.2
$$
avec $a_\sigma^0\in U_\omega(C)$ et
$a_{\bar{\sigma}}^0=(a_\sigma^0)^-$, d'o\`u
pour $\sigma$ r\'eel, 
$$
a_\sigma^0\in U_\omega(\dbR)\qquad
\text{(pour $\sigma=\bar{\sigma})$}
\tag2.12.3
$$
\endproclaim

La preuve est la m\^{e}me que celle de
Deligne (1989) 8.10.
Noter que (2.12.1) d\'etermine $a_\sigma$ \`a
multiplication \`a droite pr\`es par un \'el\'ement
de $G_\omega(\dbQ)$.
Choisissons $i$ dans $C$.
La fonctorialit\'e en $C$ montre que si
$a_\sigma$ v'erifie (2.12.1), alors
$\abar_\sigma$ v\'erifie (2.12.1) pour
$\bar{\sigma}$.
Si on impose \`a $a_\sigma$ d'\^{e}tre de la
forme (2.12.2), avec $a_\sigma^0$ dans
$U_\omega(C)$, alors
$\abar_\sigma=\abar_\sigma^0\tau(-2\pi
i)=\abar_\sigma^0\tau(2\pi i)\tau(-1)$ et
$\abar_\sigma^0$ est bien un choix possible
pour $a_{\bar{\sigma}}$.
Pour la v\'erification qu'on peut m\^{e}me
prendre $a_\sigma^0\in U_\omega(\dbR)$ si
$\sigma=\bar{\sigma}$, on renvoie \`a loc.  cit.
Si on impose \`a
$a_\sigma$ d'\^{e}tre de la forme (2.12.2),
son ind\'etermination est r\'eduite \`a
$U_\omega(\dbQ)$.
Si $\sigma$ est r\'eel et qu'on impose
en outre (2.12.3), elle est r\'eduite \`a
$$
\align
U^+(\dbQ):&=U(\dbQ)\CAP\tau(2\pi i)^{-1}
  U(\dbR)\tau(2\pi i)\\
&=\exp\left(\prod \gru_{2n}\right).
\endalign
$$

\subhead
2.13
\endsubhead
Appelons ``syst\`eme de r\'ealisations de Hodge
de Tate mixte'' (sur $k$) la donn\'ee\break
suivante.

\medskip\noindent
{\sans (a)}\enspace
un $k$-espace vectoriel gradu\'e $V_{\DR}$.
La filtration (croissante et index\'ee par les
entiers pairs) $W_{-2n}:=\oplusop_{m\Ge
n}(V_{\DR})_m$ est la {\it filtration par le
poids} et la filtration d\'ecroissante
$F^{-n}=\oplusop_{m\Le n}(V_{\DR})_m$ la {\it
filtration de Hodge}.

\medskip\noindent
{\sans (b)}\enspace
pour tout plongement $\sigma$ de $k$ dans une
cl\^{o}ture alg\'ebrique $C$ de $\dbR$, un
$\dbQ$-espace vectoriel $V_\sigma$ muni d'une
filtration croissante $W$ index\'ee par les entiers
pairs, fonctoriel en $C$.

\medskip\noindent
{\sans (c)}\enspace
pour $\sigma$ comme in (b), un isomorphisme de
comparaison $\comp_{\sigma,\DR}\colon\,
V_{\DR}\otimes_{k,\sigma}C\simover
V_\sigma\otimes C$, respectant la filtration
par le poids et fonctoriel en $C$.
Son inverse est not\'e $\comp_{\DR,\sigma}$.

On suppose ces donn\'ees telles que le
sous-espace
$$
(V_\omega)_n:=\frac{1}{(2\pi
i)^n}\,\comp_{\DR,\sigma}(\Gr_{-2n}^W(V_\sigma))
\Subset (V_{\DR})_n\otimes_{k,\sigma} C
$$
soit contenu dans $(V_{\DR})_n$ et
ind\'ependant de $\sigma$.

La cat\'egorie $\scrR_k^H$, ou simplement
$\scrR^H$, de ces syst\`emes, munie du produit
tensoriel \'evident, est une cat\'egorie
tannakienne, et 2.9 \`a 2.11 fournissent un
$\otimes$-foncteur $\real^H$ de $\MT(k)$ dans
$\scrR_k^H$.

La remarque suivante permettra de v\'erifier
\`a peu de frais le caract\`ere motivique
d'objets qui nous int\'eressent.

\proclaim{Proposition 2.14}
Le foncteur $\real^H\colon\, \MT(k)\to
\scrR_k^H$ est pleinement fid\`ele, d'image
essentielle stable par sous-objet.
\endproclaim

Le point essentiel est que le foncteur
$\real^H$ induit une injection de
$\Ext^1(\dbQ(0),\dbQ(n))$ dans
$\Ext^1(\real^H(\dbQ(0)),\real^H(\dbQ(n)))$.
Pour $n=1$, cette injectivit\'e est celle de 
$$
(\log\,\sigma)\colon\, k^*\otimes\dbQ\to
\prod\limits_{\sigma\colon\, k\to\dbC}
\dbC/2\pi i\dbQ
$$
(et un $\sigma$ suffit \`a l'injectivit\'e).
Pour $n>1$, elle est une cons\'equence de
l'injectivit\'e de l'application r\'egulateur
sur $K_{2n-1}(k)$.

\demo{Preuve}
Le foncteur fibre $\omega$ identifie la
cat\'egorie $\scrR_k^H$ \`a celle des
repr\'esentations d'un groupe alg\'ebrique
$\dbG_m\ltimes U_\omega^H$, et le foncteur
$\real^H$ au foncteur
$\Rep(G_m\ltimes U_\omega)\to\Rep(\dbG_m\ltimes
U_\omega^H)$ d\'efini par un homomorphisme
$$
\dbG_m\ltimes U_\omega^H\to \dbG_m\ltimes U_\omega.
$$
L'injectivit\'e sur $\Ext^1$ assure que ce
morphisme est surjectif, et 2.14 en
r\'esulte.
\enddemo

\subhead
2.15 Variantes
\endsubhead
(i) Puisque $\MT(\scrO_S)$ est une
sous-cat\'egorie pleine de $\MT(k)$ stable par
sous-objet, 2.14 vaut \'egalement pour le
foncteur
$$
\real^H\colon\, \MT(\scrO_S)\to \scrR_k^H
\tag2.15.1
$$

\noindent
(ii)
Supposons donn\'ee une factorisation du
foncteur 2.15.1 par une cat\'egorie
ab\'elienne $\scrR:$
$$
\real^H=F\circ \real^\scrR\colon\,
\MT(\scrO_S)\to\scrR\to\scrR_k^H,
$$
avec $\real^\scrR$ et $F$ exacts et $F$
fid\`ele.
Le foncteur $\real^\scrR$ est alors pleinement
fid\`ele d'image essentielle stable par
sous-objet.
``Pleinement fid\`ele'' r\'esulte formellement
de ``$F$ fid\`ele'' et ``$F\circ\real^\scrR$
pleinement fid\`ele''.
Que $F$ soit exact et fid\`ele implique que
pour $R$ dans $\scrR$, $F$ induit une
injection
$$
\{\text{sous-objets de $R$}\}\longrightarrow
\{\text{sous-objets de $F(R)$}\}
$$
et ``image essentielle stable par
sous-objet'' r\'esulte formellement de la
m\^{e}me propri\'et\'e pour
$F\circ\real^\scrR$.

On peut par exemple prendre pour $\scrR$ la
cat\'egorie $\scrR^{H+\ell}$ suivante de
syst\`emes de r\'ealisations: un objet est la
donn\'ee de (a) (b) (c) comme en 2.13, plus
les donn\'ees $\ell$-adiques suivantes:

\medskip\noindent
{\sans (d)}\enspace
pour $\ell$ un nombre premier et $\kbar$ une
cl\^{o}ture alg\'ebrique de $k$, un
$\dbQ_\ell$-espace vectoriel $V_\ell$, muni
d'une filtration finie croissante $W$ index\'ee par
les entiers pairs, la filtration par le poids,
d\'ependant fonctoriellement de $\kbar$.
La fonctorialit\'e en $\kbar$ doit d\'efinir une
action continue de $\Gal(\kbar/k)$, et on
demande que la repr\'esentation
$\Gr_{-2n}^W(V_\ell)$ de $\Gal(\kbar/k)$ soit
isomorphe \`a une somme de copies de
$\dbQ_\ell(n)$.

\medskip\noindent
{\sans (e)}\enspace
pour $\sigma$ un plongement de $k$ dans une
cl\^{o}ture alg\'ebrique $C$ de $\dbR$, et
pour $\kbar$ la cl\^{o}ture alg\'ebrique de
$k$ dans $C$, un isomorphisme de comparaison
$$
\comp_{\ell,\sigma}\colon\,
V_\sigma\otimes\dbQ_\ell\simover V_\ell
$$
fonctoriel en $C$.

On suppose l'existence d'un syst\`eme de
r\'eseaux $M_{\sigma,\dbZ}\Subset M_\sigma$
tels que les
$$
\comp_{\ell,\sigma}(M_{\sigma,\dbZ}\otimes \dbZ_l)
\Subset V_\ell
$$
soient stables sous $\Gal(\kbar/k)$ et
ind\'ependants de $\sigma$.

Tant dans $\scrR^H$ que dans $\scrR^{H+\ell}$,
les seuls objets simples sont les
r\'ealisations des $\dbQ(n)$.

\subhead
2.16
\endsubhead
Pour $k'$ une extension finie de $k$, on
dispose d'un $\otimes$-foncteur ``extension du
corps de base \`a $k'$'': $M\mapsto M_{(k')}$
de $\MT(k)$ dans $\MT(k')$.
L'existence et les propri\'et\'es du morphisme
``norme'' en $K$-th\'eorie assurent par
application de (0.4) que 
$$
\Ext_{\MT(k)}^1(\dbQ(n),\dbQ(m))\to
\Ext_{\MT(k')}^1(\dbQ(n),\dbQ(m))
\tag2.16.1
$$
est injectif et que, pour $k'/k$ une extension
galoisienne, l'image est form\'ee des
invariants sous Galois:
$$
\Ext_{\MT(k)}^1(\dbQ(n),\dbQ(m))\simover
\Ext_{\MT(k')}^1(\dbQ(n),\dbQ(m))^{\Gal(k'/k)}.
\tag2.16.2
$$

Par les arguments de la preuve de 2.14,
l'injectivit\'e de (2.16.1) implique que le
foncteur $M\mapsto M_{(k')}$ est pleinement
fid\`ele, d'image essentielle stable par
sous-objet.
En particulier, si un motif de Tate mixte
$M'$ sur $k'$ provient d'un motif de Tate
mixte $M$ sur $k$, ce dernier est
d\'etermin\'e \`a isomorphisme unique pr\`es
par $M'$.
Si la place finie $v'$ de $k'$ domine la
place $v$ de $k$, $M'$ est non ramifi\'e en
$v'$ si et seulement si $M$ l'est en $v$.

\subhead
2.17
\endsubhead
Par 2.16, les cat\'egories tannakiennes
$\MT(k)$, pour $k$ variable, forment un
pr\'echamp sur le site \'etale de
$\Spec(\dbQ)$.
On note $\MAT$ le champ associ\'e et on appelle
$\MAT(k)$ la cat\'egorie des motifs
d'Artin-Tate mixte sur $k$.
Explicitons cette d\'efinition.

Soit $k_1$ une extension galoisienne de $k$.
Notons $\scrC(k_1/k)$ la cat\'egorie des
extensions finies $k'$ de $k$ dans lesquelles
peut se plonger $k_1$.
Soit $\MAT(k_1/k)$ la cat\'egorie des sections
de $\MT$ au-dessus de $\scrC(k_1/k)$: la
cat\'egorie des syst\`emes d'objets
$M_{(k')}\in \MT(k')$ pour $k'\in\scrC(k_1/k)$,
et d'isomorphismes $M_{(k')(k'')}\simover
M_{(k'')}$ pour $k'\to k''$ dans
$\scrC(k_1/k)$, v\'erifiant une
compatibilit\'e pour $k'\to k''\to k'''$.
La cat\'egorie $\MAT(k)$ est la limite
inductive des $\MAT(k_1/k)$, pour
$\scrC(k_1/k)$ de plus en plus petit.

La filtration par le poids des objets des
$\MT(k')$ \'etant fonctorielle et compatible
aux extensions de corps de base, elle fournit
une filtration par le poids sur chaque objet
de $\MAT(k_1/k)$, ou de $\MAT(k)$.

La cat\'egorie $\MAT(k)$ h\'erite de $\MT(k')$
de foncteurs fibre $\omega$, $\DR$, $\sigma$
et $\ell$, et d'isomor-\break
phismes de comparaison.

\medskip\noindent
{\it Foncteur fibre $\omega$}:
le foncteur $M=(M_{(k')})_{k'\in\scrC(k_1/k)}
\mapsto M_{(k_1)}$ est une \'equivalence de
$\MAT(k_1/k)$ avec la cat\'egorie des objets de
$\MT(k_1)$ munis d'une action semi-lin\'eaire
de $\Gal(k_1/k)$
(donn\'ee de descente galoisienne).
Cette donn\'ee de descente induit une action
de $\Gal(k_1/k)$ sur
$\omega(M):=\omega(M_{(k_1)})$ et le groupe de
Galois motivique $\Aut(\omega)$ de
$\MAT(k_1/k)$ est le produit semi-direct de
$\Gal(k_1/k)$ par le groupe de Galois
motivique de $\MT(k_1)$.

Soit $\kbar$ une cl\^{o}ture alg\'ebrique de
$k$ et passons \`a la limite sur $k_1\Subset
\kbar$.
On d\'efinit
$$
\omega(M):=\omega(M_{(k_1)})\qquad
\text{pour $k_1\Subset \kbar$ assez grand}.
$$
Ce $\dbQ$-espace vectoriel gradu\'e est muni
d'une action continue de $\Gal(\kbar/k)$.

\medskip\noindent
{\it Foncteur fibre $\DR$}:
d\'efini par descente galoisienne \`a partir
des foncteurs fibres $\DR$ des $\MT(k')$.
Pour $M$ dans $\MAT(k_1/k)$ et
$k'\in\scrC(k_1/k)$, on a fonctoriellement en
$k'$
$$
M_{\DR}\otimes_k k'\simover (M_{(k')})_{\DR},
$$
et l'isomorphisme 2.10 fournit un isomorphisme
de comparaison 
$$
M_{\DR}=(\omega(M_{(k_1)})\otimes
k_1)^{\Gal(k_1/k)}.
$$

\medskip\noindent
{\it Foncteur fibre $\sigma$}:
pour $\sigma$ un plongement de $k$ dans une
cl\^{o}ture alg\'ebrique $C$ de $\dbR$.
C'est
$$
M_\sigma:=(M_{(k')})_{\sigma'}
$$
pour $\sigma$ le compos\'e $k\Subset k'
\,{\overset{\smash\sigma'}\to \longhook\,} C$ et $k'$
assez grand pour que $M_{(k')}$ soit dans
$\MT_{(k')}$.

\medskip\noindent
{\it Foncteur fibre $\ell$}:
il d\'epend du choix d'une cl\^{o}ture
alg\'ebrique $\kbar$ de $k$, et est d\'efini
comme \'etant $(M_{(k')})_\ell$ pour
$k'\Subset\kbar$ assez grand. 

\medskip\noindent
{\it Cas particulier.}
Les objets purement de poids $0$ de $\MT(k')$,
i.e. tels que $\Gr_{-2n}W=0$ pour $n\not=0$,
sont les sommes de copies de $\dbQ(0)$, et
$\omega_0$ identifie leur cat\'egorie \`a
celle des $\dbQ$-espaces vectoriels.
En cons\'equence, pour $\kbar$ une cl\^{o}ture
alg\'ebrique de $k$, le foncteur $\omega$
identifie la cat\'egorie des objets purement
de poids $0$ de $\MAT(k)$ avec celle des
$\dbQ$-espaces vectoriels munis d'une action
continue de $\Gal(\kbar/k)$ (motifs d'Artin).

Soit $k_1$ une extension galoisienne finie de
$k$.
Il r\'esulte de 2.16 que le foncteur $\MT(k)\to
\MAT(k_1/k)$ est pleinement fid\`ele d'image
essentielle stable par sous-objet.
L'image essentielle est aussi stable par
extensions:

\proclaim{Proposition 2.18}
Pour qu'un objet $M$ de $\MAT(k_1/k)$ soit dans
$\MT(k)$, il faut et il suffit que l'action de
$\Gal(k_1/k)$ sur $\omega(M)$ soit triviale.
\endproclaim

Si $M$ est purement de poids $0$, la
proposition r\'esulte du ``cas particulier''
ci-dessus.
Le cas o\`u $M$ est pur, i.e. purement d'un
seul poids, en r\'esulte par tensorisation avec
un $\dbQ(n)$.
Il reste \`a v\'erifier une stabilit\'e par
extensions.

\demo{$1^{\text{ere}}$ preuve}
Une suite spectrale d'Hochschild-Serre
fournit, pour $M$ et $N$ dans $\MT(k)$, des
isomorphismes
$$
\Ext_{\MAT(k_1/k)}^p(M,N)\simover \Ext^p
(M_{(k_1)},N_{(k_1)})^{\Gal(k_1/k)}.
$$

Pour v\'erifier que
$$
\Ext_{\MT(k)}^p(M,N)\to\Ext_{\MAT(k'/k)}^p
(M,N)=\Ext_{\MT(k_1)}^p(M_{(k_1)},
N_{(k_1)})^{\Gal(k_1/k)}
$$
est bijectif pour $p=0,1$ et injectif pour
$p=2$, il suffit par d\'evissage de le
v\'erifier pour $M$ et $N$ de la forme
$\dbQ(m)$ et $\dbQ(n)$.
Le cas $p=2$ r\'esulte de la nullit\'e des
$\Ext^2$, le cas $p=0$ r\'esulte de ce que
$\Hom=0$, sauf pour $m=n$, o\`u
$\End(\dbQ(n))=\dbQ$, et $p=1$ est 2.16.
\enddemo

Plutot que de d\'efinir la suite spectrale
d'Hochschild-Serre requise, on peut la
paraphraser:

\demo{$2^{\text{\`eme}}$ preuve}
Si $M$ est une extension de $\dbQ(n)$ par
$\dbQ(n+r)$, que $M_{(k_1)}$ admette une
donn\'ee de descente de $k_1$ \`a $k$,
respectant la structure d'extension de
$\dbQ(n)$ par $\dbQ(n+r)$, implique que la
classe de $M_{(k_1)}$ dans
$\Ext_{\MT(k_1)}^1(\dbQ(n),\dbQ(n+r))$ est fixe
sous $\Gal(k_1/k)$.
Elle provient donc de la classe d'une
extension $\Mtil$ de $\dbQ(n)$ par
$\dbQ(n+r)$ dans $\MT(k)$.
On v\'erifie, utilisant que
$\Hom_{\MT(k_1)}(\dbQ(n),\dbQ(n+r))=0$, que $M$
est isomorphe \`a l'image de $\Mtil$.

Prouvons par r\'ecurrence sur $r\Ge 0$ que 2.18
est vrai sous l'hypoth\`ese additionnelle que
pour un $n$ on ait $M=W_{-2n}M$ et
$W_{-2(n+r+1)}M=0$.

Le cas o\`u $r=0$ ($M$ pur) a d\'ej\`a \'et\'e
trait\'e.
Le cas $r=1$ se r\'eduit \`a celui d'une
extension de $\dbQ(n)$ par $\dbQ(n+1)$.

Pour $r\Ge 1$, l'hypoth\`ese de r\'ecurrence
assure que $M'=M/W_{-2(n+r)}M$ et
$M''=W_{-2(n+1)}M$ proviennent de $\Mtil'$ et
$\Mtil''$ dans $\MT(k)$.
Ce sont des extensions de la forme:
$\Gr_{-2n}^WM$ par $N$ et $N$ par
$\Gr_{-2(n+r)}^W M$.
Par nullit\'e des $\Ext^2$ dans $\MT(k)$, il
existe dans $\MT(k)$ un objet $\Mtil$,
extension it\'er\'ee de $\Gr_{-2n}^WM$ par $N$
par $\Gr_{-2(n+r)}^WM$, redonnant $\Mtil'$
et $\Mtil''$.
L'image de $\Mtil$ dans $\MAT(k_1/k)$ diff\`ere de
$M$ par une extension de $\Gr_{-2n}^WM$ par
$\Gr_{-2(n+r)}^W M$.
\enddemo

Cette derni\`ere provient d'une extension dans
$\MT(k)$, par r\'eduction au cas d'une
extension de $\dbQ(n)$ par $\dbQ(n+r)$, 
et ceci permet de corriger $\Mtil$ pour que son
image devienne isomorphe \`a $M$.

\newpage


\dspace
\subhead
3. Cohomologie et groupe fondamental 
\endsubhead

\subhead
3.1
\endsubhead
Soit $X$ un espace topologique s\'epar\'e et
connexe.
On note ${}_bP_a$ l'ensemble des classes
d'homotopie de chemins de $a$ \`a $b$, et
$\delta\gamma$ la composition des chemins
$_cP_b\times{} _bP_a\to {}_cP_a$.
Noter l'ordre des facteurs, choisi pour que si
$V$ est un syst\`eme local sur $X$ et que
$\gamma v\in V_b$ denote le transport de $v\in
V_a$ le long de $\gamma\in {}_bP_a$, on ait
$(\delta\gamma)v=\delta(\gamma v)$.

Nous ferons sur $X$ les hypoth\`eses (a) (b)
(c) suivantes.
(a) La th\'eorie classique du $\pi_1$
s'applique.
Par exemple: $X$ localement connexe par arcs
et localement simplement connexe.
(b) Le groupe fondamental est de type fini.
(c) Pour une projection $p\colon\, X^n\times
X^m\to X^n$ et un faisceau constant, $Rp_*$
commute au passage aux fibres.
Ces conditions sont remplies pour $X$
hom\'eomorphe au compl\'ement dans un poly\`edre
fini d'un sous-poly\`edre, par exemple pour
$X$ une vari\'et\'e alg\'ebrique complexe.

Fixons un corps $k$, le corps des
coefficients.
Soient $a\in X$, $\Lambda$ l'alg\`ebre de
groupe $k[\pi_1(X,a)]$ et $I$ son id\'eal
d'augmentation.
La donn\'ee d'un $\Lambda$-module \`a gauche
$V$, i.e. d'une repr\'esentation lin\'eaire de
$\pi_1(X,a)$ sur un $k$-espace vectoriel,
\'equivaut \`a celle d'un syst\`eme local 
$V\uptil$ sur $X$ de fibre $V$ en $a$.
Si on prend pour $V$ le $\Lambda$-module
$\Lambda$, on obtient le syst\`eme local de
fibre en $x$ l'espace vectoriel $k[{}_xP_a]$
de base ${}_xP_a$.
La structure de $\Lambda$-module \`a droite de
$\Lambda$ en fournit une sur ce syst\`eme
local $\Lamuptil$, i.e. une action \`a
droite de $\pi_1(X,a)$.
L'\'el\'ement $\gamma$ de $\pi_1(X,a)$ agit
sur ${}_xP_a$ par composition \`a droite.
Le syst\`eme local $(\Lambda/I^n)\uptil$ est
le quotient $\Lamuptil \otimes_\Lambda
\Lambda/I^n$ de $\Lamuptil$.

Le syst\`eme local $(\Lambda/I^n)\uptil$ est
extension it\'er\'ee des $n$ syst\`emes locaux
triviaux $I^p/I^{p+1}$ ($0\Le p<n$).
Sa fibre $(\Lambda/I^n)_a\!\!\!\uptil$ en $a$
est $\Lambda/I^n$ et l'unit\'e $1$ d\'efinit
$$
k\to (\Lambda/I^n)_a\!\!\!\uptil\colon\,
\lambda\longmapsto \lambda.1.
\tag3.1.1
$$
Le syst\`eme local $(\Lambda/I^n)\uptil$
est universel parmi les syst\`emes locaux $E$
extension it\'er\'ee de $n$ syst\`emes locaux
triviaux et munis de $k\to E_a$, car
$(\Lambda/I^n,1)$ est universel parmi les
$\Lambda$-modules $M$ tels que $I^nM=0$ munis
de $m_0\in M$.

\subhead
3.2
\endsubhead
Dans ce qui suit, nous laisserons de nombreux
signes ambigus.
Plusieurs conventions sont possibles, et
aucune ne nous semble clairement meilleure.
Elles n'importent gu\`ere, car elles
conduisent \`a des complexes isomorphe, avec
pour l'isomorphisme au pis l'ambigu{\umi}t\'e
d'un signe global.
Typiquement, changer de convention remplace un
complexe par un autre ayant les m\^{e}mes
composantes, et l'isomorphisme est donn\'e par
des signes d\'ependant de multidegr\'es.

\subhead
3.3
\endsubhead
Fixons $a,b\in X$ et consid\'erons dans $X^n$
les sous-espaces $b=t_1$,
$t_1=t_2,\dotsc,t_n=a$, not\'es $Y_0,\dots,Y_n$.
Pour $J\Subset\Delta_n:=\{0,1,\dotsc,n\}$,
soient $Y_J$ l'intersection des $Y_j$ ($j\in
J$) et $k_J$ le faisceau constant $k$ sur
$Y_J$, prolong\'e par z\'ero sur $X^n$.
Pour $J\Subset K$, on a $Y_K\Subset Y_J$ et on
dispose d'un morphisme de restriction $k_J\to
k_K$.
Si, pour $K=\{k_0,\dotsc,k_p\}$, les $k_\ell$
\'etant pris dans l'ordre croissant, et
$J=\{k_0,\dotsc,\khat_i,\dotsc,k_p\}$, on
l'affecte du signe $(-1)^i$, on obtient les
composantes de la diff\'erentielle d'un
complexe de faisceaux sur $X^n$:
$$
\underbark\to\oplus k_{\{j\}}\to\cdots\to
\oplusop_{\vert J\vert=p}
k_J\to\cdots\to k_{\Delta_n}.\tag3.3.1
$$
Noter que le dernier terme $k_{\Delta_n}$ est
nul sauf pour $n=0$ ou $a=b$, auxquels cas il
est r\'eduit \`a $k$ en $t_1=\ldots=t_n=a$.
Si on supprime le premier terme et qu'on
convient que les degr\'es vont de $0$ \`a $n$,
on obtient la r\'esolution \v{C}echiste
altern\'ee du faisceau constant $k$ sur la
r\'eunion $\cup Y_i$ des $Y_i$ pour le
recouvrement ferm\'e fini par les $Y_i$,
prolong\'ee par $0$ sur $X^n$.

D\'efinissons ${}_b\scrK_a\left<n\right>$,
ou simplement ${}_b\scrK_a$, comme \'etant le
complexe (3.3.1),
avec $k_{\Delta_n}$ omis, et les degr\'es
allant de $0$ \`a $n$.
Si $n>0$ et que $a\not=b$, c'est une
r\'esolution du prolongement par $0$ du
faisceau constant $k$ sur $X^n-\cup Y_i$, et
$$
\dbH^*(X^n,{}_b\scrK_a)=H^*(X^n\mod\,\cup
Y_i,k)\quad\text{(pour $n\not=0$,
$a\not=b$)}\tag3.3.2
$$

Si $a=b$, la derni\`ere diff\'erentielle
(3.3.1) fournit un morphisme de complexes
$$
\align
&{}_a\scrK_a\to k_{t_1=\cdots=
t_n=a}[-n],\tag3.3.3\\
\intertext{qui induit}
&\dbH^n(X^n,{}_a\scrK_a)\to k.\tag3.3.4
\endalign
$$

Si on r\'ep\`ete ces constructions en prenant
$a$ et $b$ comme param\`etres, on obtient un
complexe ${}_*\scrK_a$ sur $X\times X^n\times
X$, et, pour $p$ la projection $X\times
X^n\times X\to X\times X$, la restriction de
${}_*\scrK_*$ \`a $X^n\simover p^{-1}(b,a)$
s'identifie \`a ${}_b\scrK_a$.
Traitant seulement $b$ comme param\`etre, on
obtient de m\^{e}me ${}_*\scrK_a$ sur $X\times
X^n$.
Les $\dbH^i(X^n,{}_b\scrK_a)$ sont les fibres
du syst\`eme local $R^ip_*({}_*\scrK_*)$ sur
$X\times X$.

Nous allons, d'apr\`es Beilinson, v\'erifier la

\proclaim{Proposition 3.4}
(i) Pour $i<n$, $\dbH^i(X^n,{}_b\scrK_a)=0$.
\ (ii) Pour $b$ variable, le syst\`eme local
sur $X$ des $\dbH^n(X^n,{}_b\scrK_a)$, muni en
$b=a$ de (3.3.4), est le dual du syst\`eme
local $(\Lambda/I^{n+1})\uptil$ de 3.1, muni
du transpos\'e de (3.1.1).
\endproclaim

\demo{Preuve}
Pour $n=0$, $X^n$ est r\'eduit \`a un point,
${}_b\scrK_a$ est le faisceau constant $k$ en
degr\'e $0$ et (i) (ii) sont clairs.
Prouvons 3.4 par r\'ecurrence sur $n\Ge 1$;
l'hypoth\`ese de r\'ecurrence est la
validit\'e de 3.4 pour $m$ si $m<n$.

Soit $q$ la premi\`ere projection: $X^n\to X$.
Nous utiliserons la suite spectrale de Leray
pour $q$.
Ceci revient \`a \'ecrire que
$R\Gamma(X^n,\quad)=R\Gamma(X,Rq_*\quad)$ et
\`a filtrer $Rq_*$ par la filtration canonique
$\tau$.
Plus pr\'ecis\'ement, nous utiliserons une
filtration ayant une description uniforme pour
$a=b$ ou $a\not=b$, et qui revient \`a la
filtration canonique si $a\not=b$.

Un d\'evissage montre que la formation de
$Rq_*{}_b\scrK_a$ commute au passage aux
fibres.
Si on \'ecrit $X^n$ comme $X\times X^{n-1}$,
${}_b\scrK_a$ sur $X^n$ s'identifie au
c\^{o}ne$[-1]$ de
$$
\align
{}_*\scrK_a\left<n-1\right>\to{}_b\scrK_a\left<
n-1\right> &\quad\text{sur}\quad \{b\}\times
  X^{n-1}\tag3.4.1\\
\oplus~k_{\b{$\ssize t$}=a} &\quad\text{sur}\quad
\{a\}\times X^{n-1}
\endalign
$$
(complexe simple associ\'e au complexe double
de premi\`ere (resp. deuxi\`eme) colonne la
source (resp. but) de (3.4.1)).
\enddemo

On reconna{\hati}t \`a gauche les $k_J$ pour
$0\notin J$ et $J\not=\{1,\dotsc,m\}$, \`a
droite les $k_J$ d'une part pour $0\in J$, de
l'autre pour $J=\{1,\dotsc,n\}$.

Si $a\not=b$, la restriction de ${}_b\scrK_a$
\`a $q^{-1}(t)=X^{n-1}$ est donc: pour
$t\not=a,b$: ${}_t\scrK_a\left<n-1\right>$;
pour $t=b$: le c\^{o}ne$[-1]$ de
l'application identique de
${}_b\scrK_a\left<n-1\right>$; pour $t=a$: le
c\^{o}ne$[-1]$ de (3.3.3).
D'apr\`es l'hypoth\`ese de r\'ecurrence,
$R^iq_*{}_b\scrK_a=0$ si $i<n-1$.
Pour $i=n-1$, c'est un sous-faisceau du
syst\`eme local $(\Lambda/I^n)\uptil\upcheck$.
Il co{\umi}ncide avec
$(\Lambda/I^n)\uptil\upcheck$ si $t\not=a,b$,
est nul en $t=b$ et la fibre en a est le noyau
du transpos\'e de (3.3.1).
Pour $i\Ge n$, $R^iq_*{}_b\scrK_a$ est un
syst\`eme local sur $X-\{b\}$ prolong\'e par
$0$ sur $X$.
La suite spectrale de Leray v\'erifie donc
$E_2^{pq}=0$ pour $q\Le n-2$ ou $p=0$, et
$$
\align
&\dbH^i(X^n,{}_b\scrK_a)=0\quad\text{pour}\quad
 i<n\tag3.4.2\\
&\dbH^n(X^n,{}_b\scrK_a)=\dbH^1(X,R^{n-1}q_*
  {}_b\scrK_a)\quad\text{(pour $a\not=b$).}
\endalign
$$
Ce $\dbH^1$ est le groupe des extensions du
faisceau constant $k$ par
$(\Lambda/I^n)\uptil\upcheck$, trivialis\'ees
en $t=b$, et trivialis\'ees apr\`es avoir
pouss\'e par $(\Lambda/I^n)\upcheck\to k$ en
$t=a$.
Dualement, c'est le groupe des extensions
$\Etil$ de $(\Lambda/I^n)\uptil$ par $k$,
trivialis\'ees en $t=b$, et trivialis\'ees
au-dessus de $k\Subset\Lambda/I^n$ en $t=a$.

Traduisons en termes de $\Lambda$-modules.
A $\Etil$ correspond un $\Lambda$-module $E$,
extension de  $\Lambda/I^n$ par le module
d'augmentation $k$.
Il est annul\'e par $I^{n+1}$.
La donn\'ee additionelle en a est celle d'un
rel\`evement $\tilde{1}$ du g\'en\'erateur $1$
de $\Lambda/I^n$ dans $E$.
Soit $\alpha\colon\, I^n/I^{n+1}\to k$ la
restriction \`a $I^n/I^{n+1}$ du morphisme de
$\Lambda/I^{n+1}$ dans $E\colon\,
\lambda\mapsto\lambda \tilde{1}$.
L'extension $E$ est l'extension $E_\alpha$ de
l'extension $\Lambda/I^{n+1}$ de $\Lambda/I^n$
par $I^n/I^{n+1}$ en poussant par $\alpha$:
$$
\minCDarrowwidth{13pt}
\CD
0 @>>> I^n/I^{n+1} @>>> \Lambda/I^{n+1}
  @>>> \Lambda/I^n @>>> 0\\
@. @VV{\alpha}V @VVV @VVV\\
0 @>>> \dbQ @>>> E @>>> \Lambda/I^{n+1} @>>>
0.
\endCD
\tag3.4.3
$$

Pour $\alpha\colon\, I^n/I^{n+1}\to k$.
L'extension $E_\alpha\!\!\!\uptil$ de
$(\Lambda/I^n)\uptil$ par le syst\`eme local
trivial $k$ se d\'eduit de
$(\Lambda/I^{n+1})\uptil$ en poussant par
$$
\alpha\colon\text{ (syst\`eme local trivial
$I^n/I^{n+1})$}\to k.
$$
La trivialisation $\beta$ en $b$ fournit un
diagramme
$$
\minCDarrowwidth{13pt}
\CD
0 @>>> I^n/I^{n+1} @>>>
(\Lambda/I^{n+1})_b\!\!\!\uptil @>>>
(\Lambda/I^n)_b\!\!\!\uptil @>>> 0\\
@. @VV{\alpha}V @VVV @|\\
0 @>>> \dbQ @>>> (E_\alpha)_b\!\!\!\uptil
  &{\overset{\sssize\beta}\to\leftrightarrows} 
&(\Lambda/I^n)_b\!\!\!\uptil @>>> 0.
\endCD
\tag3.4.4
$$
De $\alpha$ et $\beta$ on d\'eduit une forme
lin\'eaire $\gamma$ sur
$(\Lambda/I^{n+1})_b\!\!\!\uptil$: \`a
$x\in(\Lambda/I^{n+1})_b\!\!\!\uptil$, d'image
$\xbar$ dans $(\Lambda/I^n)_b\!\!\!\uptil$,
attacher
$$
\gamma(x):=\image(x)-\beta(\xbar)\in\dbQ
\Subset(E_\alpha\!\!\!\uptil)_b.
$$
On a $\alpha=\gamma\vert(I^n/I^{n+1})$, et
$\gamma$ d\'etermine $\beta$.
Ceci fournit l'isomorphisme
$$
\dbH^n(X^n,{}_b\scrK_a)=
((\Lambda/I^{n+1})_b\!\!\!\uptil)\upvee
$$
annonc\'e.

Pr\'esent\'es un peu diff\'erement, ces
arguments continuent \`a s'appliquer pour
$a=b$.
L'image par $Rq_*$ du complexe ${}_b\scrK_a$,
c\^{o}ne$[-1]$ de (3.4.1) est encore un
c\^{o}ne$[-1]$.
On d\'efinit sa filtration croissante $\Fil$
en filtrant source et but par la filtration
canonique $\tau_{\Le}$.
Le gradu\'e $\Gr_i^{\Fil}$ s'envoie par un
quasi-isomorphisme sur le complexe
$$
R^iq_*{}_*\scrK_a\left<n-1\right>\to
\text{ (sa fibre en a)}_a\oplus
(k\text{ si }i=n-1)\tag3.4.5
$$
en degr\'es $0$ et $1$, translat\'e par
$[-i]$.
Pour $i<n-1$, (3.4.5) est nul.
Pour tout $i$, $\dbH^0((3.4.5))=0$.
La suite spectrale pour $R\Gamma(X,\quad)$ et
$\Fil$ donne donc
$$
\align
&\dbH^i(X^n,{}_b\scrK_a)=0\quad\text{pour
$i<n$}\tag3.4.6\\
&H^n(X^n,{}_b\scrK_a)=\dbH^1(X,(3.4.5)_{i=n-1}).
\endalign
$$
Ce $\dbH^1$ classifie les extensions du
faisceaux constant $k$ par
$(\Lambda/I^n)\uptil\upvee=R^{n-1}q_*{}_*\scrK_a
\left<n-1\right>$, trivialis\'ees apr\`es
avoir pouss\'e par le morphisme (3.4.5), et
les arguments donn\'es pour $a\not=b$ peuvent
\^{e}tre r\'ep\'et\'e.

Il reste, pour $a=b$, \`a identifier (3.3.4).
Le complexe (3.3.1) (degr\'es $0$ \`a $n+1$)
est le complexe simple associ\'e au complexe
triple de faisceaux sur $X^n=X\times X^{n-1}$:
$$
\CD
{}_*\scrK_a\left<n-1\right> @>>>{}_a\scrK_a
  \left<n-1\right>\text{ sur }q^{-1}(a)\\
@VVV @VVV\\
k[-(n-1)]_{{\b{$\ssize t$}}=a} @>{\sim}>> 
  k[-(n-1)]_{{\b{$\ssize t$}}=a}
\endCD
$$
et (3.3.4) est induit par le morphisme
$$
\CD
{}_*\scrK_a\left<n-1\right> @>>> 
  {}_a\scrK_a
  \left<n-1\right>\text{sur $\pr^{-1}(a)$}\\
@VVV \null\kern-3.0 cm @VVV\\
k[-(n-1)]_{{\b{$\ssize t$}}=a}
  @>>> \null\kern-3.0 cm 0
\endCD
\quad
\longrightarrow
\quad
\CD
0 @>>> \null\kern-25pt0\\
@VVV \null\kern-25pt@VVV\\
0 @>>> k[-(n-1)].
\endCD
$$

Filtrant comme pr\'ec\'edemment, on en
d\'eduit que (3.2.1) est le $\dbH^1$ du
morphisme de complexes
$$
\CD
(\Lambda/I^n)\upvee\uptil @>>>
(\Lambda/I^n)_a\!\!\!\uptil \oplus k_a\\
@VVV @VVV\\
0 @>>> k_a
\endCD
$$
(La fl\`eche verticale droite donnant le
morphisme est la diff\'erence entre
l'\'evaluation sur $k\Subset\Lambda/I^n$ et
l'identit\'e de $k_a$.)

Une extension de $k$ par
$(\Lambda/I^n)\upvee\uptil$, trivialis\'ee
apr\`es avoir poun\'e par (3.3.4), se dualise
en une extension de $(\Lambda/I^n)\uptil$ par
$k$ qui en $a$ est et trivialis\'ee, et
trivialis\'ee au-desus de $k\Subset\Lambda/I^n
=(\Lambda/I^n)_a\!\!\!\uptil$.
On lui attache la diff\'erence au-desus de $k$
entre ces deux trivialisations.
On laisse au lecteur le soin de v\'erifier que
cela revient \`a \'evaluer en $1$ une forme
lin\'eaire sur $\Lambda/I^{n+1}$.

\remark\nofrills{\bf 3.5 Remarque.}\enspace
Dans cette identification de 
$(k[{}_bP_a]\otimes_{\Lambda}\Lambda/I^{n+1})\upvee$
avec le groupe des extensions de
$(\Lambda/I^n)\uptil$ par $k$, munies de
trivialisations convenables en $a$ et $b$, le
sous-espace $(k[{}_bP_a]\otimes_\Lambda
\Lambda/I^n)\upvee$ est celui des extensions
triviales, avec la trivialisation \'evidente en
$a$ et une trivialisation en $b$.
L'injection de $(k[{}_bP_a]/I^n)\upvee\to
(k{}_aP_b]/I^{n+1})\upvee$ est donc le
$\dbH^n$ de
$$
\left({}_b\scrK_a\left<n-1\right>\text{ sur }
\{b\}\times X^{n+1}\right)[-1]\to
{}_b\scrK_a\text{ sur }X^n\tag3.5.1
$$
d\'eduit de la description du second membre
comme c\^{o}ne$[-1]$ de (3.4.1).
\endremark

\subhead
3.6
\endsubhead
Nous regarderons un objet simplicial
tant\^{o}t comme un foncteur contravariant sur la
cat\'egorie de tous les ensembles totalement
ordonn\'es finis non vides, tant\^{o}t comme un
foncteur sur la sous-cat\'egorie des
$\Delta_n=\{0,\dotsc,n\}$.
Cela revient au m\^{e}me (\'equivalence de
cat\'egories).
Si la valeur sur $\Delta_n$ d'un objet
simplicial $S$ est not\'ee $S_n$, on notera
$S_I$ sa valeur sur $I$ totalement ordonn\'e.
Elle est fonctorielle en $I$ et si $I$ a $n+1$
\'el\'ement, l'unique isomorphisme de $I$ avec
$\Delta_n$ fournit un isomorphisme de $S_I$
avec $S_n$.

Le morphismes $\Delta_n\to\Delta_1$ sont les
$$
\varphi_i\colon\, 0,\dotsc,i-1\longmapsto
0;\,\,
i,\dotsc,n\longmapsto 1\qquad
(0\Le i\Le n+1).\tag3.6.1
$$
Nous identifierons par cette num\'erotation
$X^{\Hom(\Delta_n,\Delta_1)}$ \`a $X\times
X^n\times X$.
L'ensemble simplicial
$\Delta_n\mapsto\Hom(\Delta_n,\Delta_1)$
contient l'ensemble simplicial
$\Delta_n\mapsto\Hom(\Delta_n,\partial\Delta_1)=
\{\varphi_0,\varphi_{n+1}\}$.
L'espace cosimplicial $\Delta_n\mapsto
X^{\Hom(\Delta_n,\Delta_1)}$ se projette donc sur
l'espace cosimplicial $\Delta_n\mapsto
X^{\Hom(\Delta_n,\partial\Delta_1)}=
X\times X$.
Prenant la fibre en $(b,a)$, on obtient un
espace cosimplicial $\Delta_n\mapsto X^n$.
Les $\partial_i\colon\, X^{n-1}\to X^n$
($0\Le i\Le n$) sont les
$$
\alignat2
t_1,\dotsc,t_{n-1}\longmapsto
  &~b,t_1,\dotsc,t_{n-1} &\quad &(i=0)\\
&~t_1,\dotsc,t_i,t_i,\dotsc,t_{n-1} 
  &\quad &(i\not=0,n)\\
&~t_1,\dotsc,t_{n-1},a &\quad &(i=n)\,\,,
\endalignat
$$
d'images les $Y_i$ de 3.3.
Plus g\'en\'eralement, pour $J\Subset\{
0,\dotsc,n\}$ de compl\'ement $I$, suppos\'e
non vide, $p=\vert I\vert-1$ et $\varphi_I$
l'inclusion de $I$, le morphisme cosimplicial
$\varphi_I\colon\, X^p\to X^n$ identifie $X^p$
avec $Y_J\Subset X^n$.
Pour $J$, $J'$ de compl\'ements $I$, $I'$ et
$p=\vert I\vert-1$, $p'=\vert I'\vert-1$, si
$J\Subset J'$ et que $\varphi_{II'}$ est
l'inclusion de $I'$ dans $I$, le morphisme
cosimplicial $\varphi_{II'}\colon\, X^{p'}\to
X^p$ s'identifie \`a l'inclusion de $Y_{J'}$
dans $Y_J$.
Noter que pour $J=\{j_0,\dotsc,\Hatj_k,\dotsc,j_p
\}\Subset J'=\{j_0,\dotsc,j_p\}$ donnant lieu
\`a $I'=\{i_0,\dotsc,{\Hati}_\ell,\dotsc,i_q\}
\Subset I=\{i_0,\dotsc,i_q\}$ (les $i$ et $j$
sont pris dans l'ordre croissant), on a
$p+q=n$, $k+\ell=i_\ell=j_k$, et le signe
$(-1)^k$ que nous avons employ\'e en 3.3 peut
se r\'ecrire
$$
(-1)^k=(-1)^\ell\cdot
(-1)^{j_k}=\eps(J)(-1)^\ell\eps(J')
$$
pour
$$
\eps(J)=\prod\limits_{j\in J}(-1)^j.\tag3.6.2
$$

\subhead
3.7
\endsubhead
Choisissons des complexes $S_m^*$ qui
calculent la cohomologie des $X^m$, et soient
fonctoriels pour les morphismes cosimpliciaux
entre les $X^m$.
Plus pr\'ecis\'ement: on prend pour $S_m^*$ le
complexe des sections globales d'une
r\'esolution fonctorielle $\scrS_m^*$ du
faisceau constant $k$ sur $X^m$, telle que
$\Gamma(X^m,\scrS_m^*)\to
R\Gamma(X^m,\scrS_m^*)$ soit un
quasi-isomorphisme.
A $m$ variable, les $S_m^*$ forment un
syst\`eme simplicial de complexes.
Les $\scrS_m^*$ fournissent une r\'esolution
de ${}_b\scrK_a\left<n\right>$ sur $X^n$:
avec les notations de 3.6, r\'esoudre $k_J$
par $\varphi_{I*}\scrS_I^*$ et utiliser la
fonctorialit\'e de $\scrS_I^*$ en $I$.

Un {\it syst\`eme de coefficients} $c$ sur le
simplexe standard $\Delta_n$, \`a valeurs dans
une cat\'egorie additive, attache \`a chaque
facette $F\Subset\Delta_m$ un objet $c(F)$,
contravariant en $F$ pour les morphismes
d'inclusion.
Il d\'efinit un complexe de cha{\hati}nes
$$
C_p(\Delta_n,c)=\bigoplus\limits_{\vert F\vert
=p+1}c(F).
$$
Dualement, un cosyst\`eme: $c(F)$ covariant en
$F$, d\'efinit un complexe de cocha{\hati}nes.

Si $s_{\smbul}$ est un objet simplicial, on
note encore $C_*(\Delta_n,s_{\smbul})$ le
complexe $C_*(\Delta_n,c)$ pour $c(F):=s_F$.
Dualement pour $s^{\smbul}$ cosimplicial.

Le complexe
$$
\Gamma(X^n,\text{ r\'esolution de
${}_b\scrK_a$ d\'efinie par $\scrS_{\smbul}$}),
$$
calcule l'hypercohomologie de ${}_b\scrK_a$.
Il s'identifie (avec les signes (3.6.2)) au
complexe simple associ\'e au complexe double
($=$ complexe de complexes)
$$
C_*(\Delta_n,S_{\smbul}^*),
$$
translat\'e par $[-n]$, et 3.4 (ii) se
r\'ecrit
$$
(k[{}_b P_a]\otimes_\Lambda
\Lambda/I^{n+1})\upvee=H^0C_*(\Delta_n, S_{\smbul}^*).
\tag3.7.1
$$

\example\nofrills{3.8 Exemple:}\enspace
prenons pour $S_m^*$ le complexe des cocha{\hati}nes
singuli\`eres de $X^m$.
Il ne rentre pas tout \`a fait dans le cadre
pr\'ec\'edent, mais on peut l'y ramener: si
$\scrS_m^*$ est le complexe des cocha{\hati}nes
singuli\`eres localis\'ees sur $X^m$ (le
faisceau associ\'e au pr\'efaisceau des
cocha{\hati}nes singuli\`eres),
$S_m^*\to\Gamma(X^m,\scrS_m^*)$ est un
quasi-isomorphisme.
Le complexe $S_m^*$ des cocha{\hati}nes
singuli\`eres est le dual du complexe $S_*^m$
des cha{\hati}nes singuli\`eres.
A $m$ variable, $S_*^m$ est cosimplicial en $m$,
et $C_*(\Delta_n,S_{\smbul}^*)$ est le dual de
$C^*(\Delta_n,S_*^{\smbul})$.
On a donc
$$
k[{}_b P_a]\otimes_\Lambda
\Lambda/I^{n+1}=H^0C^*(\Delta_n,S_*^{\smbul})
\tag3.8.1
$$

Pour $X=[0,1]$, $a=0$ et $b=1$, ${}_b P_a$ est
r\'eduit \`a un \'el\'ement et chaque
$k[{}_bP_a]\otimes_\Lambda \Lambda/I^{n+1}$
est r\'eduit \`a $k$.
Au membre de gauche, les
$C^*(\Delta_n,S_*^{\smbul})$ forment un
syst\`eme projectif, de morphismes de
transition d\'eduit, apr\`es dualisation, de
(3.5.1\quad):
$$
C^*(\Delta_n,S_*^{\smbul})\to
C^*(\Delta_{n-1},
S_*^{\smbul})
$$
est induit par $\partial_0\colon\,\Delta_{n-1}
\hookrightarrow\Delta_n$.
\endexample

Pour qu'un cocycle de degr\'e zero de
$C^*(\Delta_n,S_*^{\smbul})$ corresponde au
chemin de $0$ \`a $1$, il faut et il suffit
que sa valeur en $\{n\}=\partial_0^n(\Delta_0)\Subset
\Delta_n$ soit $1\in S^0$, le complexe des
cha{\hati}nes singuli\`eres du point.
On a $[0,1]=\vert\Delta_1\vert$.
Soit $\sigma_m$ le simplexe
$\vert\Delta_m\vert
\to\vert\Delta_1\vert^m=[0,1]^m$ de
coordonn\'ees les d\'eg\'en\'erescences
$\varphi_i$ ($1\Le i\le m$) de (3.6.1).
Son image est $\{\underbart\vert 1\Ge
t_1\Ge\ldots\Ge t_m\Ge 0\}$.
Comme cocycle, on peut prendre celui qui \`a
chaque face de dimension $p$ de $\Delta_n$
attache $\pm\sigma_p$, le signe \'etant $+$
pour $p=0$ et d\'ependant des conventiones de
signes d\'efinissant
$C^*(\Delta_n,S_*^{\smbul})$ pour $p>0$.

Par fonctorialit\'e en $X$, on en d\'eduit
que, pour $X$, l'image dans
$H^0C^*(\Delta_n,S_*)$ d'un chemin
$\gamma\colon\, [0,1]\to X$ de $a$ \`a $b$ est
repr\'esent\'ee par le cocycle qui \`a chaque
face de dimension $p$ de $\Delta_n$ attache
$\pm\gamma\circ\sigma_p$:
$\vert\Delta_p\vert\to[0,1]^p\to X^p$, le
signe \'etant comme ci-desssus.

\subhead
3.9
\endsubhead
La litt\'erature donne une autre recette pour
d\'eduire $(k[{}_bP_a]\otimes_\Lambda
\Lambda/I^{n+1})\upvee$ de l'objet simplicial
$S_{\smbul}^*$ de la cat\'egorie des complexes
consid\'er\'e en 3.7:

\medskip\noindent
{\sans (a)}\enspace
prendre le complexe normalis\'e
$NS_{\smbul}^*$ de cet objet simplicial;

\smallskip\noindent
{\sans (b)}\enspace
prendre le tronqu\'e $\sigma_{\Ge -n}
NS_{\smbul}^*$ qui ne garde que les
$NS_m^*$ pour $m\Le n$;

\smallskip\noindent
{\sans (c)}\enspace
prendre le complexe simple associ\'e au complexe
double $\sigma_{\Ge -n}NS_{\smbul}^*$ et son
$H^0$.

\medskip
L'\'equivalence entre les deux constructions
r\'esulte du lemme suivant.
Rappelons qu'une {\it cat\'egorie karoubienne}
est une cat\'egorie additive dans laquelle
tout endomorphisme idempotent admet un noyau,
et que l'\'equivalence de Dold-Puppe:
(objet simplicial
$S_{\smbul}$)$\to$ (complexe $NS_{\smbul}$ \`a
degr\'es $\Le 0$) est disponible dans une
telle cat\'egorie ($NS_k$ est facteur
direct de $S_k$)

\proclaim{Proposition 3.10}
Si $S_{\smbul}$ est un objet simplicial d'une
cat\'egorie karoubienne, le complexe 
$C_*(\Delta_n,S_{\smbul})$ et le
tronqu\'e normalis\'e $\sigma_{\Ge
-n}NS_{\smbul}$ sont fonctoriellement homotopes.
\endproclaim

Nous v\'erifierons par la m\'ethode des
mod\`eles acycliques l'\'enonc\'e dual, pour
un objet cosimplicial $S^{\smbul}$,
$\sigma_{\Le n}NS^{\smbul}$ et
$C^*(\Delta_n,S^{\smbul})$.

\demo{Preuve}
Nous noterons $C_{\dbZ}$ la cat\'egorie
pr\'eadditive ayant les m\^{e}mes objets
qu'une cat\'egorie $\scrC$, pour laquelle
$\Hom(A,B)$ est le groupe ab\'elien librement
engendr\'e par $\Hom_{\scrC}(A,B)$.
Nous noterons $\scrC_{\kar}$ son enveloppe
karoubienne, obtenue par adjonction formelle
de sommes finies et de facteurs directs
noyaux d'endomorphismes idempotents.
Le foncteur contravariant
$$
X\longmapsto h_{\dbZ}^X\colon\,
\scrC_{\kar}\to\HOM(\scrC,Ab),\quad 
h_{\dbZ}^X(Y)=\Hom_{\scrC_{\kar}}(X,Y)
$$
est pleinement fid\`ele.
Il induit une \'equivalence de $\scrC_{\kar}$
avec la sous-cat\'egorie de la cat\'egorie
ab\'elienne $\HOM(\scrC,Ab)$ d'objets les
facteurs directs de sommes finies de
foncteurs
$$
Y\longmapsto\text{ groupe ab\'elien librement
engendr\'e par $\Hom_{\scrC}(X,Y)$},
$$
pour $X$ dans $\scrC$.
Ces foncteurs sont des objets projectifs de la
cat\'egorie ab\'elienne $\HOM(\scrC,Ab)$, car
pour $C$ dans $\scrC$
$$
\Hom(h_{\dbZ}^C,F)=F(C)
$$
(variante du lemme de Yoneda).
\enddemo

Si $\Delta$ est la cat\'egorie des $\Delta_n$,
l'objet cosimplicial $(\Delta_n)$ de
$\Delta_{kar}$ est universel en ce sens que,
pour $\scrA$ une cat\'egorie karoubienne,
$F\mapsto F((\Delta_n))$ est un \'equivalence
$$
\text{(foncteurs $\Delta_{\kar}\to\scrA)\to$
(objets cosimpliciaux de $\scrA$)}
$$
Pour v\'erifier l'\'enonc\'e $(3.10)^*$ dual de
3.10, et montrer que l'\'equivalence
d'homotopie peut \^{e}tre donn\'ee par des
formules universelles, il suffit de v\'erifier
$(3.10)^*$ dans le cas universel: pour l'objet
$(\Delta_n)$ de $\Delta_{\kar}$.
Si $\Delta_{\Le n}$ est la sous-cat\'egorie de
$\Delta$ d'objets les $\Delta_m$ pour $m\Le
n$, tant $\sigma_{\Ge-n}N\Delta_{\smbul}$ que
$C_*(\Delta_n,\Delta_{\smbul})$ sont dans
$(\Delta_{\Le n})_{\kar}$, et c'est dans
$(\Delta_{\Le n})_{\kar}$ que nous
travaillerons.

Appliquons le foncteur $X\mapsto h_{\dbZ}^X$
\`a $\sigma_{\Le n}N\Delta^{\smbul}$ et \`a
$C^*(\Delta_n,\Delta^{\smbul})$.
On obtient deux complexes d'objets projectifs
de la cat\'egorie ab\'elienne
$\HOM(\Delta_{\Le n},Ab)$.
Nous montrerons qu'il sont l'un et l'autre une
r\'esolution projective du foncteur constant
$\dbZ$, donc sont homotopes, et on conclut par la
pleine fid\'elit\'e de $X\mapsto h_{\dbZ}^X$.
Si on \'evalue en $\Delta_p$ ($p\Le n$), on
obtient les complexes $\sigma_{\Ge
-n}NS_{\smbul}$ et $C_*(\Delta_n,S_{\smbul})$
de 3.10 pour $S_{\smbul}$ le complexe
cha{\hati}nes simpliciales de $\Delta_p$.
Ce sont

\medskip\noindent
{\sans (a)}\enspace
le complexe des cha{\hati}nes non
d\'eg\'en\'er\'ees de $\Delta_p$, \'egal \`a
son tronqu\'e car $p\le n$.

\smallskip\noindent
{\sans (b)}\enspace
le complexe des cha{\hati}nes combinaisons
lin\'eaires de
$\Delta_q\to\Delta_n\times\Delta_p$, avec la
premi\`ere projection injective.

\medskip
Ces deux complexes sont augment\'es vers $\dbZ$,
fonctoriellement en $\Delta_p$.
Puisque $\vert\Delta_p\vert$ est contractile,
et que le complexe (a) calcule son homologie,
il est une r\'esolution de $\dbZ$.
Le complexe (b) est le complexe des
cha{\hati}nes non d\'eg\'en\'er\'ees du
sous-polytope simplicial suivant de
$\vert\Delta_n\times\Delta_p\vert$:

\parindent=40pt
\medskip\noindent
\Item{(3.10.1)}
la r\'eunion, pour $\varphi\colon\,
\Delta_n\to\Delta_p$, des graphes
$\Gamma_\varphi$ des morphismes
$\vert\varphi\vert\colon\,\vert\Delta_n\vert\to
\vert\Delta_p\vert$.

\parindent=20pt
\medskip
Pour v\'erifier que le complexe (b) est une
r\'esolution de $\dbZ$, il suffit donc de
v\'erifier le

\proclaim{Lemme 3.10.2}
L'espace (3.10.1) est contractile.
\endproclaim

\demo{Preuve}
Ordonnons lexicographiquement l'ensemble des
applications croissantes
$$
\varphi\colon\,\{0,\dotsc,n\}\to
\{0,\dotsc,p\}.
$$
Montrons que pour tout $\varphi\not=0$,
$$
\Gamma_\varphi\CAP\bigcup\limits_{\varphi' < \varphi}
\Gamma_{\varphi'}
$$
est une union non vide de faces du simplexe
$\Gamma_{\varphi}$, mais non de toutes les faces
de $\Gamma_\varphi$, de sorte que
$\bigcup\limits_{\varphi'<\varphi}\Gamma_{\varphi'}$
est r\'etracte par d\'eformation de
$\bigcup\limits_{\varphi'\Le
\varphi}\Gamma_{\varphi'}$.

Quels que soient $\varphi'$ et $\varphi$, si
$F\Subset\{0,\dotsc,n\}$ est l'ensemble des
$i$ tels que $\varphi(i)=\varphi'(i)$, et
$\vert F\vert$ la facette correspondante
$\Delta_m$, $\Gamma_\varphi\CAP
\Gamma_{\varphi'}$ est le graphe de la
restriction de $\vert\varphi\vert$ ou
$\vert\varphi'\vert$ \`a $\vert F\vert$.

Identifions $\Gamma_\varphi$ \`a
$\vert\Delta_m\vert$ par la projection sur
$\Delta_m$.
Si $\varphi'<\varphi$, soit $i$ minimal tel
que $\varphi'(i)<\varphi(i)$.
Si $\varphi''$ est d\'efini par
$\varphi''(j)=\varphi'(j)$ (resp.
$\varphi(j)$) pour $j\le i$ (resp. $j>i$), on
a $\varphi''<\varphi$ et
$\Gamma_\varphi\CAP\Gamma_{\varphi'}\Subset
\Gamma_\varphi\CAP\Gamma_{\varphi''}$, qui
s'identifie \`a la face d'indice $i$ de
$\vert\Delta_m\vert$.
L'intersection
$\Gamma_\varphi\CAP\bigcup\limits_{\varphi'<\varphi}
\Gamma_{\varphi'}$ est donc r\'eunion non vide
de faces.

Pour que la face d'indice $i$ soit obtenue, il
faut et il suffit que $\varphi(i)\not=0$, et
que, si $i\not=0$, on ait
$\varphi(i-1)<\varphi(i)$.
Parce que $n\Ge p$, toute les faces ne sont
pas obtenues: soit $\varphi(i)=\varphi(i+1)$
pour un $i$, soir $n=m$, $\varphi$ est
l'identit\'e, et $\varphi(0)=0$.
\enddemo

\subhead
3.11
\endsubhead
Si $X$ est une vari\'et\'e diff\'erentiable,
que $k=\dbR$ et que, dans 3.7, on calcule la
cohomologie de $X^n$ \`a l'aide du complexe de
de Rham, on obtient un objet simplicial
$\Omega^p(X^n)$ de la cat\'egorie des
complexes (degr\'e diff\'erentiel gradu\'e $p$,
degr\'e simplicial $n$).
Soit $N\Omega^*(X^*)$ son normalis\'e, et
$\bfs N\Omega^*X^*$ le complexe simple
associ\'e, dont la composante de degr\'e $k$
est $\oplusop_{p-n=k}N\Omega^p(X^k)$.
C'est un proche cousin du complexe des
int\'egrales it\'er\'ees de Chen.
Nous nous proposons d'expliquer la relation
entre les deux.

Une $q$-forme $\beta$ sur l'espace
${}_b\Omega_a(X)$ des chemins $C^\infty$ de
$a$ \`a $b$ est pour Chen la donn\'ee, pour
toute vari\'et\'e $U$ et toute famille
$C^\infty$ $^\colon\, U\to {}_b\Omega_a(X)$
de chemins de $a$ \`a $b$ param\'etr\'ee par
$U$, d'une $q$-forme (not\'ee
$\gamma^*\beta$) sur $U$, cette donn\'ee
v\'erifiant que pour $f\colon\, V\to U$ on
a $(\gamma f)^*\beta=f^*(\gamma^*\beta)$.
A une $p$-forme sur $X^n$, Chen attache comme
suit une $(p-n)$-forme sur ${}_b\Omega_a(X)$:

\medskip\noindent
{\sans (a)}\enspace
r\'ep\'etant avec param\`etres la construction
de 3.8, on attache \`a une famille $\gamma$ de
chemins param\'etr\'ee par $U$ un morphisme
$$
\gamma_n=\gamma\circ\sigma_n\colon\,
U\times\vert\Delta_n\vert\to
U\times[0,1]^n\to X^n;
$$

\noindent
{\sans (b)}\enspace
\`a la forme $\alpha$ sur $X^n$ on attache
l'int\'egrale de $\gamma_n^*\alpha$ le long
des fibres de $U\times\vert\Delta_n\vert\to U$.
Plus pr\'ecis\'ement, on attache
$$
\left<\alpha\right>=(-1)^{n(n+1)/2}\int_{U\times
\vert\Delta_n\vert/U}\gamma_n^*\alpha,
$$
o\`u l'int\'egrale le long des fibres est
d\'efinie de sorte que pour $\beta'$ sur
$U$ et $\beta''$ sur $\vert\Delta_n\vert$,
on ait 
$$
\int_{U\times\vert\Delta_n\vert/U}
\pr_1^*\beta'\w \pr_2^*\beta''=\beta'.
\int_{\vert\Delta_n\vert}\beta''.
$$

Le signe est choisi pour avoir, si $\alpha$
est de degr\'e $p$
$$
d\left<\alpha\right>=\left<d\alpha\right>+
(-1)^p\left<\sum(-1)^ij_i^*\alpha\right>\,\,:
$$
compatibilit\'e de la diff\'erentielle
ext\'erieure sur ${}_b\Omega_a(X)$ avec la
diff\'erentielle $D$ du complexe simple
associ\'e au complexe double des
$\Omega^*(X^n)$ (pour * le premier degr\'e).

Pour la diff\'erentielle transpos\'ee de
$\Sing_*X^*$, si $\gamma$ est un chemin de $a$
\`a $b$, c'est $\sum(-1)^p\gamma\circ\sigma_p$
qui est un cycle.

Chen remplace les $\Omega^*X^n$ par le
syst\`eme simplicial de sous-complexes
$\otimesop^{n}\Omega^*(X)$.
Le {\it complexe des int\'egrale it\'er\'ees}
de Chen est l'image par $\left<\quad\right>$
de la somme des $\otimesop^{n} \Omega^*(X)$.
Cette image est form\'ee des sommes de 
$$
\left<\alpha_1\vert\ldots\vert\alpha_n\right>:=
\left<\pr_1^*\alpha_1\w\ldots\w\pr_n^*\alpha_n
\right>.
$$

Les d\'eg\'en\'erescences $s_i\colon\, X^n\to
X^{n-1}$ sont les projections $(x_1\ldots
x_n)\to
(x_1\ldots\xhat_i\ldots x_n)$ et le complexe
normalis\'e de $\otimesop^{n}\Omega^*X$, dans
sa version ``diviser par les $\im(s_i)$'' est
$$
N\otimesop^{n}\Omega^*X=\otimesop^{n}(\Omega^*X/
\text{constantes}).\tag3.11.1
$$
Par K\"unneth, les $\otimesop^{n}\Omega^*X$
sont quasi-isomorphes aux $\Omega^*X^n$, et
donc les $N\otimesop^{n}\Omega^*X$ aux
$N\Omega^*X^n$.

Si un des $\alpha_i$ est de degr\'e $0$,
$\left<\alpha_1\vert\ldots\vert\alpha_n\right>$
est nul et $\left<\quad\right>$ induit donc
$$
\left<\quad\right>\colon\,
N\otimesop^{n}\Omega^*X\to\text{complexe des
int\'egrales it\'erees.}
$$
Les deux membres sont filtr\'es par $n$; par
les $\oplusop_{m\Le n}N\otimesop^{m}\Omega^*X$
et leurs images.
Chen (1973) (Lemme 4.3.1) prouve 
que, pour cette filtration $\Fil$ du complexe
des int\'egrales it\'eries, $\Gr_n^{\Fil}$ est
le quotient $\otimesop^{n}(\tau_{\Ge 1}\Omega^*X)$
de $\otimesop^{n}\Omega^*X$.
Le complexe $\tau_{\Ge 1}\Omega^*X$ est le
quotient de $\Omega^*X$ par $\Omega^0X$ et son
image par $d$.
La vari\'et\'e $X$ \'etant connexe, le
morphisme de $(\Omega^*X/\text{constantes})$
dans $\tau_{\Ge 1}\Omega^*X$ est un
quasi-isomorphisme, et les
$$
\Fil_nN\otimesop^{n}\Omega^*X
\to\Fil_n(\text{int\'egrales it\'er\'ees})
$$
sont donc des quasi-isomorphismes.

\remark{Remarque}
On peut donner du complexe des int\'egrales
it\'er\'ees une description directe, en terme
de l'alg\`ebre gradu\'ee $\Omega^*(X)$ et des
morphismes d'\'evaluation en $a$ et $b$, sans
invoquet
les formes diff\'erentielles sur
${}_b\Omega_a$: l'image de $\oplusop_{m\Le
n}\otimesop^{m}\Omega_X^*$ est son quotient par
les $\alpha_1\otimes\ldots\otimes\alpha_m$
tels qu'un $\alpha_i$ soit de degr\'e $0$, et
leurs diff\'erentielles (construction bar
r\'eduite de Chen 1976).
\endremark

\subhead
3.12
\endsubhead
Soit $X$ une vari\'et\'e alg\'ebrique lisse
connexe sur un corps $K$ et $a,b\in X(K)$.
Comme dans le cas topologique 3.6, le sch\'ema
cosimplicial $X^{\Hom(\Delta_n,\Delta_1)}$
s'envoie dans le sch\'ema cosimplicial constant
$X\times
X=X^{\Hom(\Delta_n,\partial\Delta_1)}$ et,
prenant la fibre en $(b,a)\in X\times X$, on
obtient un sch\'ema cosimplicial
$\Delta_n\mapsto X^n$.
Regardons-le comme un objet cosimplicial de la
cat\'egorie additive $\SmCor(K)$.
Quel que soit $n$, on peut lui attacher le
complexe
$$
C^*(\Delta_n,X^*)[n]
$$
(3.7) de $\SmCor(K)$, ou, dans l'enveloppe
karoubienne $\SmCor(K)_{\kar}$ de $\SmCor(K)$,
le complexe homotope $\sigma_{\Le n}NX^*$.
Ces complexes fournissent des objets,
naturellement isomorphes d'apr\`es 3.10, de la
cat\'egorie triangul\'ee motivique $\DM(K)$.
Notons les ${}_b\Omega_a^{[n]}(X)$.

Appliquons \`a ${}_b\Omega_a^{[n]}(X)$ le
foncteur ``r\'ealisation''  cohomologique,
suivi de $H^0$.
D'apr\`es 3.4, pour tout plongement complexe
$\sigma$ de $K$, dans la r\'ealisation de
Betti $H_\sigma$ correspondante (coefficients
rationnels), on obtient
$(\dbQ[{}_bP_a]\otimes_\Lambda
\Lambda/I^{n+1})\upvee$ pour $X(\dbC)$.

Supposons que $X$, comme objet de $\DM(K)$,
soit de Tate mixte, i.e. dans $\DMT(K)$.
Les ${}_b\Omega_a^{[n]}(X)$ sont alors
\'egalement dans $\DMT(K)$.
Si $K$ est un corps de nombres, ou plus
g\'en\'eralement quand on dispose de la
conjecture d'annulation de Beilinson-Soul\'e,
on peut, apr\`es tensorisation avec $\dbQ$,
appliquer $H^0$ et dualiser pour obtenir un
objet de la cat\'egorie $\MT(K)$ des motifs de
Tate mixtes sur $K$.

Si on travaille avec $\sigma_{\Le n}NX^*$, les
formules classiques en topologie (voir 
Wojtkowiak (1993)) fournissent, apr\`es
passage \`a la limite inductive en $n$, une
structure d'alg\`ebre commutative dans la
cat\'egorie des $\Ind$-objets de $\MT(K)$, et
$$
{}_bA_a(X):=\colim_n(H^0\sigma_{\Le
n}NX^*)\upvee
$$
m\'erite le nom d'alg\`ebre affine de l'espace
motivique des chemins de $a$ \`a $b$.
Pour trois points $a$, $b$, $c$, on dispose
aussi (loc. cit.)d'un coproduit
$$
{}_c A_a(X)(X)\to {}_cA_b(X)\otimes
{}_bA_a(X)
$$
correspondant \`a la composition des chemins.
Pour $a=b$, $\Spec({}_aA_a(X))$ est le {\it
groupe fondamental rendu unipotent motivique}.

\subhead
3.13
\endsubhead
Prenons pour $X$ le compl\'ement, dans une
droite projective $P$ sur $K$, d'un ensemble
fini de points rationnels, et montrons comment
comparer ces constructions \`a celles de
Deligne (1989).
En r\'ealisation de de Rham, $X$ fournit le
complexe (\`a diff\'erentielle nulle) des
formes diff\'erentielles \`a p\^{o}les
logarithmiques:
$$
K\mapright{0}\Omega_P^1(\log\,S).\tag3.13.1
$$
Pour $X^n$, on obtient la uissance tensorielle
$n^{\ieme}$ de (3.13.1).
Ces puissances tensorielles forment un
syst\`eme simplicial de complexes.
Normalisant, ce qui revient \`a supprimer $K$
dans (3.13.1) (cf. (3.11.1)), on obtient
apr\`es troncation en $n$
$$
\oplusop_{m\Le n}\Omega_p^1(\log\,S)^{\otimes m}.
\tag3.13.2
$$

La comparaison avec la r\'ealisation de Betti
est par les int\'egrales it\'er\'ees de Chen,
et ceci est compatible aux constructions de
Deligne (1989).
Quant \`a la r\'ealisation en cohomologie
$\ell$-adique, les arguments prouvant 3.4
s'appliquent aussi bien dans le cadre
$\ell$-adique.

\newpage


\dspace
\subhead
4. Le groupe fondamental unipotent motivique
d'une vari\'et\'e rationnelle
\endsubhead

\subhead
4.1
\endsubhead
Soient $k$ un corps de nombres, $P$ une droite
projective sur $k$, $X$ le compl\'ement d'un
ensemble fini $S$ de points rationnels et $x,y\in
X(k)$.
Dans Deligne (1989) \S{13}, nous avons d\'efini un
``espace de (classes d'homotopie de) 
chemins motivique'' not\'e $P_{y,x}$.
Sous les hypoth\`eses ci-dessus, et apr\`es
oubli de l'aspect cristallin de loc. cit.,
c'est un sch\'ema affine en la cat\'egorie
tannakienne $\scrR^{H+\ell}$ de 2.15.
Nous le noterons ici $P_{y,x}^{H+\ell}$, ou
$P_{y,x}^{H+\ell}(X)$ s'il y a lieu de
pr\'eciser $X$.
On dispose d'une ``composition des chemins''
$P_{z,y}^{H+\ell}\times P_{y,x}^{H+\ell}
\to P_{z,x}^{H+\ell}$.
Posons $\pi_1^{H+\ell}(X,x):=P_{x,x}^{H+\ell}$.
Pour la composition des chemins, c'est un
$\scrR^{H+\ell}$-sch\'ema en groupes affine et
$P_{y,x}^{H+\ell}$ est un espace principal
homog\`ene \`a droite sous
$\pi_1^{H+\ell}(X,x)$, \`a gauche sous
$\pi_1^{H+\ell}(X,y)$: un
$(\pi_1^{H+\ell}(X,x),
\pi_1^{H+\ell}(X,y))$-bitorseur.
Pour $\sigma$ un plongement de $k$ dans une
cl\^{o}ture alg\'ebrique $C$ de $\dbR$, la
r\'ealisation $\sigma$ de
$\pi_1^{H+\ell}(X,x)$ est l'enveloppe
alg\'ebrique prounipotente de $\pi_1(X(C),x)$.
Si $I$ est l'id\'eal d'augmentation de
l'alg\`ebre de groupe $\dbQ[\pi_1(X(C),x)]$, c'est
le spectre de l'alg\`ebre de Hopf commutative
$$
\colim(\dbQ[\pi_1(X(C),x)]/I^N)\upvee.
$$
i.e. c'est l'enveloppe alg\'ebrique
pro-unipotente de $\pi_1(X(C),x)$.

\subhead
4.2
\endsubhead
Soit $A_{y,x}^{H+\ell}$ l'alg\`ebre affine
(2.6) de $P_{y,x}^{H+\ell}$.
D'apr\`es 3.12, 3.13, le 
Ind-objet $A_{y,x}^{H+\ell}$ de
$\scrR^{H+\ell}$ est l'image par le foncteur
r\'ealisation d'un Ind-objet $A_{y,x}$ de
$MT(k)$.
D'apr\`es 2.14, 2.15 (ii), le produit
$A_{y,x}^{H+\ell}\otimes A_{y,x}^{H+\ell}\to
A_{y,x}^{H+\ell}$ provient d'un produit sur
$A_{y,x}$, et $P_{y,x}^{H+\ell}$ est la
r\'ealisation du $MT(k)$-sch\'ema affine
$P_{y,x}:=\Spec(A_{y,x})$.
De m\^{e}me, la composition des chemins est
d\'efinie par un coproduit $A_{z,x}^{H+\ell}\to
A_{z,y}^{H+\ell}\otimes A_{y,x}^{H+\ell}$, ce
coproduit se rel\`eve \`a $MT(k)$ et d\'efinit
un morphisme $P_{z,y}\times P_{y,x}\to
P_{z,x}$.

Nous appellerons ici $P_{y,x}$ l'{\it espace
de chemins de $x$ \`a $y$ motivique}.
Nous appelerons $P_{x,x}$ le groupe
fondamental motivique de $X$ en $x$ et le
noterons $\pi_1^{\mot}(X,x)$ ou simplement
$\pi_1(X,x)$.
Mise en garde: bien que la notation ne
l'indique pas, il s'agit de $\pi_1$ rendus
unipotents.

\subhead
4.3
\endsubhead
Appelons ``point base de $X$'' soit un point
rationnel de $X=P-S$, soit un vecteur tangent
non nul en un point de $S$ (``points base \`a
l'infini'').
Pour $x$ et $y$ deux points-base de $X$, 
Deligne (1989) d\'efinit \S{15} un
``espace de chemins motivique de $x$ \`a
$y$''.
Comme en 4.1, c'est, apr\`es oubli d'un
aspect cristallin, un
$\scrR^{H+\ell}$-sch\'ema affine.
Nous le noterons ici $P_{y,x}^{H+\ell}$, et,
pour $x=y$, $\pi_1^{H+\ell}(X,x)$.

\proclaim{Th\'eor\`eme 4.4}
Avec les notations de 4.1 et 4.3, si $x$ et
$y$ sont deux points-base de $X$, le
$\scrR^{H+\ell}$-sch\'ema $P_{y,x}^{H+\ell}$
est motivique, i.e. la r\'ealisation d'un
$MT(k)$-sch\'ema $P_{y,x}$.
\endproclaim

Si $S$ est vide ou r\'eduit \`a un point,
$X=P-S$ est simplement connexe, les
$P_{y,x}^{H+\ell}$ sont r\'eduits \`a un point
et 4.4 est trivial.
Nous pouvons donc supposer,
et supposerons, que $\vert S\vert\Ge
2$.
Le cas o\`u $x$ et
$y$ sont \`a distance finie ayant d\'ej\`a
\'et\'e trait\'e, nous supposerons
au moins l'un de $x$ et $y$ \`a l'infini.

\medskip\noindent
{\sans Cas 1.}
$x\in X(k)$ et $y$ est un vecteur tangent non
nul en $\ybar\in S$.

Soit $Z$ le $\MT(k)$-sch\'ema
des homomorphismes de $\dbQ(1)$
dans $\pi_1^{\mot}(X,x)$.
Identifiant homomorphismes de sch\'emas en
groupe unipotents, et homomorphismes entre
leurs alg\`ebres de Lie, on peut en donner la
description suivante: Lie $\pi_1^{\mot}(X,x)$
est une pro-alg\`ebre de Lie dans $MT(k)$, en
particulier un pro-objet de $MT(k)$ (c'est la
limite projective des Lie~$(\pi_1^{\mot}(X,x)/\grZ^N)$, pour $\grZ$ la
s\'erie centrale descendante),
et on d\'efinit $Z$ comme \'etant le
$MT(k)$-sch\'ema pro-vectoriel (2.6) Lie
$\pi_1^{\mot}(X,x)(-1)$.

Dans $Z^{H+\ell}:=\real^{H+\ell}(Z)$,
on dispose du sous-sch\'ema ferm\'e ``classe de
conjugaison de la monodromie locale autour de
$\ybar\,$''.
Voici sa d\'efinition.
Ecrivons $\dbQ(1)^{H+\ell}$ pour
$\real^{H+\ell}(\dbQ(1))$.
On dispose d'un morphisme ``monodromie locale
autour de $\ybar\,$'':
$$
\dbQ(1)^{H+\ell}\to\pi_1^{H+\ell}(X,y).
\tag{4.4.1}
$$
Dans (4.4.1), $\dbQ(1)^{H+\ell}$ est vu comme
un $\scrR^{H+\ell}$-sch\'ema vectoriel, avec sa
structure de groupe.
Si $T_{\ybar}$ est l'espace tangent en
$\ybar$, et que
$T_{\ybar}^*:=T_{\ybar}-\{0\}$, il y a lieu de
voir (4.4.1) comme un morphisme
$$
\pi_1^{H+\ell}(T_{\ybar}^*,y)\to\pi_1^{H+\ell}
(X,y).
$$
La donn\'ee de (4.4.1) \'equivaut \`a celle de
$$
\dbQ(1)^{H+\ell}\to\Lie\,\pi_1^{H+\ell}
(X,y).
\tag4.4.2
$$
Dans (4.4.2), $\dbQ(1)^{H+\ell}$ est un objet
de $\scrR^{H+\ell}$, muni de sa structure
d'alg\`ebre de Lie commutative.

On dispose aussi de l'espace de chemins
$P_{y,x}^{H+\ell}$ et, par composition des
chemins, de
$$
(p,\gamma)\longmapsto p^{-1}\gamma p\colon\,
P_{y,x}^{H+\ell}\times\pi_1^{H+\ell}(X,y)\to
\pi_1^{H+\ell}(X,x).
$$
Composant avec (4.4.1), on en d\'eduit
$$
P_{y,x}^{H+\ell}\times\dbQ(1)^{H+\ell}\to
\pi_1^{H+\ell}(X,x).
\tag4.4.2
$$

Le groupe $\pi_1^{H+\ell}(X,y)$ agit par
composition \`a gauche sur
$P_{y,x}^{H+\ell}$.
Divisant \`a gauche par
$\dbQ(1)^{H+\ell}$, vu par (4.4.1) comme
sous-groupe de $\pi_1^{H+\ell}(X,y)$, on d\'efinit
$$
P_{\ybar,x}^{H+\ell}:=\dbQ(1)^{H+\ell}\setminus
P_{y,x}^{H+\ell}
$$
et (4.4.2) se factorise par
$P_{\ybar,x}^{H+\ell}\times\dbQ(1)^{H+\ell}$,
d\'efinissant
$$
P_{\ybar,x}^{H+\ell}\to Z^{H+\ell}\tag4.4.3
$$

\proclaim{Lemme 4.5}
Le morphisme (4.4.3) est un plongement ferm\'e
de $\scrR^{H+\ell}$-sch\'emas.
\endproclaim

\demo{Preuve}
Il suffit de le v\'erifier en r\'ealisation de
de Rham, o\`u l'assertion se r\'eduit \`a la
suivante.
\enddemo

\proclaim{Lemme 4.6}
Soit $\scrL$ l'alg\`ebre de Lie sur $k$
engendr\'ee  par des \'el\'ements $(e_s)_{s\in
S}$ soumis \`a la seule relation que $\sum
e_s=0$.
On suppose $\vert S\vert\Ge 2$, on fixe
$\ybar\in S$ et on pose $e:=e_{\ybar}$.
Soit $G$ le sch\'ema en groupe prounipotent
limite projective des
$G_N:=\exp(\scrL/\grZ^N)$.
Le morphisme de $G/\exp(ke)$ dans $\Lie(G)$:
$g\mapsto \ad\,g(e)$, identifie $G/\exp(ke)$
\`a un sous-sch\'ema ferm\'e de $\Lie(G)$,
l'orbite adjointe de $e$.
\endproclaim

Dans 4.6, $\Lie(G)$ est le sch\'ema pro-vectoriel
limite projective des sch\'emas vectoriels
$\Lie(G_N)$.

\demo{Preuve}
Les orbites de l'action d'un groupe unipotent
sur un sch\'ema affine sont ferm\'ees Demazure,
Gabriel (1970) IV 2.2.7).
Si $e_N$ est l'image de $e$ dans
$\scrL/\grZ^N$ et $Z_N\Subset G_N$ son
fixateur, l'application $g\mapsto \ad\,g(e_N)$
de $G_N/Z_N$ dans $\Lie(G_N)$ est donc un
plongement ferm\'e.
Par passage \`a la limite projective,
l'application $g\mapsto\ad\,g(e)$ de
$\lim\,G_N/Z_N$ dans $\Lie(G)=\lim\,\Lie(G_N)$
est \'egalement un plongement ferm\'e.

L'alg\`ebre de Lie du 
sous-groupe $Z_N$ de $G_N$ est le
centralisateur de $e_N$ dans
$\scrL/\grZ^N$.
Fixons $b$ dans $S$ distinct de $\ybar$.
L'alg\`ebre de Lie $\scrL$ est librement
engendr\'ee par les $e_s$ ($s\not=b$).
D'apr\`es 4.7 ci-dessous, le centralisateur de
$e$ est r\'eduit \`a $ke$.
Munissons $\scrL$ de la graduation pour
laquelle les $e_s$ sont de degr\'e $1$.
La sous-alg\`ebre $\grZ^N$ est la
sous-alg\`ebre ``degr\'e $\Ge N$''.
L'application lin\'eaire
$\ad\,e\colon\,\scrL\to\scrL$ \'etant de
degr\'e $1$ et de noyau r\'eduit \`a $ke$, on a
$$
\Lie\,Z_N=\Ker(\ad\,e_N\colon\,\scrL/\grZ^N\to
\scrL/\grZ^N)=ke_N+\grZ^{N-1}/\grZ^N.
$$

La projection de $G_{N+1}/Z_{N+1}$ sur 
$G_N/Z_N$ se factorise donc par
$G_N/\exp(ke_N)$ et
$$
G/\exp(ke)=
\lim\,G_N/\exp(ke_N)\simover \lim\,G_N/Z_N(e).
$$
\enddemo

\proclaim{Lemme 4.7}
Dans une alg\`ebre de Lie libre $\scrL$ sur
$k$, le centralisateur d'un g\'en\'erateur $e$
est r\'eduit aux multiples de ce
g\'en\'erateur.
\endproclaim

D'apr\`es Reutenauer (1993) 2.10, 
dans une alg\`ebre de Lie
libre, deux \'el\'ements non nuls quelconques
qui commutent sont multiples l'un de l'autre.
Le cas particulier 4.7 est, plus simplement,
cons\'equence de ce que l'alg\`ebre de Lie
libre $\Lib(\{e\}\CUP A)$ est le produit
semi-direct de $ke$ et de
$\Lib((\ad\,e^n(a))_{a\in A,n\Ge 0}$ et de ce
que la d\'erivation $\ad\,e^n(a)\mapsto
\ad\,e^{n+1}(a)$ de l'alg\`ebre associative
libre engendr\'ee par les $\ad\,e^n(a)$ est
injective.

Il r\'esulte de 4.5 que $P_{\ybar,x}^{H+\ell}$
est motivique.
Soient en effet $A$ l'alg\`ebre de $Z$ et
$A^{H+\ell}:=\real^{H+\ell}(A)$ celle de
$Z^{H+\ell}$.
D'apr\`es 4.5, $P_{\ybar,x}^{H+\ell}$ est
isomorphe \`a un sous-sch\'ema ferm\'e de
$Z^{H+\ell}$, d\'efini par un id\'eal
$a^{H+\ell}$ de $A^{H+\ell}$.
D'apr\`es  2.14, 2.15, cet id\'eal est la
r\'ealisation d'un id\'eal $a$ de $A$, et
$P_{\ybar,x}^{H+\ell}$ est isomorphe \`a la
r\'ealisation de $\Spec(A/a)$.

\subhead
4.8 Fin du cas $1$ de la preuve de 4.4
\endsubhead
Soit $b\in S$ distinct de $\ybar$, et
$X':=P-\{\ybar,b\}$.
Identifions $\dbP$ et $\dbP^1$, par un
isomorphisme envoyant $\ybar$ sur $0$ et $b$
sur $\infty$, donc $X'$ sur
$A^*=\Spec\,k[u,u^{-1}]$.
Soit $t:=$ du $(y)$ la coordonn\'ee du vecteur
tangent $y$ en $u=0$, et notons encore $t$ le
point de $\dbG_m$ de coordonn\'ee $t$.
Les espaces de chemins $P_{y,x}^{H+\ell}(X')$
et $P_{t,x}^{H+\ell}(\dbG_m)$ sont
canoniquement isomorphes (Deligne (1989) \S{14}).
Puisque $t$ et $x$ sont des points-base \`a
distance finie de $\dbG_m$, $P_{y,x}^{H+\ell}$
est motivique.
C'est d'ailleurs simplement un torseur de
Kummer: on a $\pi_1^{\mot}(\dbG_m,x)=\dbQ(1)$,
$t/x\in k^*$ d\'efinit une extension de
$\dbQ(0)$ par $\dbQ(1)$, donc un
$\dbQ(1)$-torseur, et $P_{y,x}^{H+\ell}$ 
est la r\'ealisation de ce torseur.

\proclaim{Lemme 4.9}
Le morphisme
$$
P_{y,x}^{H+\ell}(X)\to
P_{y,x}^{H+\ell}(X')\times
P_{\ybar,x}^{H+\ell}(X)
$$
$(1^{\text{\rm en}}$ facteur: fonctorialit\'e
en $X$; $2^{\text{\rm e}}$ facteur: passage ou
quotient$)$ est un isomorphisme.
\endproclaim

Puisque tant $P_{y,x}^{H+\ell}(X')$ que
$P_{\ybar,x}^{H+\ell}(X)$ sont motiviques 4.9 
implique que $P_{y,x}^{H+\ell}(X)$ est motivique.

\demo{Preuve de 4.9}
Il suffit de v\'erifier 4.9 en r\'ealisation
de de Rham.
Avec les notations de 3.6, l'assertion se
r\'eduit \`a la suivante.
\enddemo

\proclaim{Lemme 4.10}
Regardons $\dbG_a=\exp(ke)$ comme un quotient
de $G$ par $e\mapsto e$, $e_a\mapsto 0$ pour
$a\not=\ybar,b$.
L'application
$$
G\to \exp(ke)\times\exp(ke)\setminus G
$$
est un isomorphisme.
\endproclaim

\demo{Preuve}
$G$ est le produit semi-direct de $\exp(ke)$
par le noyau de la projection sur $\exp(ke)$,
cf. la preuve de 4.7.
\enddemo

\subhead
4.11
\endsubhead
Pour finir la preuve de 4.4 il reste, compte
tenu de ce que $P_{y,x}^{H+\ell}$, est
isomorphe comme $\scrR^{H+\ell}$-sch\'ema \`a
$P_{x,y}^{H+\ell}$, \`a traiter le

\medskip\noindent
{\sans Cas 2.}
$x$ et $y$ sont des points base \`a l'infini.

Choisissons un point base $z\in X(k)$.
On a une expression de $P_{y,x}^{H+\ell}$
comme compos\'e de bitorseurs
$$
P_{y,x}^{H+\ell}=P_{y,z}^{H+\ell}
\times_{\pi_1^{H+\ell}(X,z)}  P_{z,x}^{H+\ell}
$$
o\`u au second membre chaque
$\scrR^{H+\ell}$-sch\'ema est motivique.
Que $P_{y,x}^{H+\ell}$ soit motivique en
r\'esulte.

\subhead
4.12
\endsubhead
Soit $X$ une vari\'et\'e alg\'ebrique lisse
sur un corps de nombres $k$, et $x,y\in X(k)$.
On suppose $X$ s\'epar\'e.
Sous l'hypoth\`ese que $H^1(\Xbar,\scrO)=0$
pour $\Xbar$ une compactification lisse de $X$,
nous avons d\'efini dans Deligne (1989) \S{13} un
``espace de chemins'' $P_{y,x}$ qui est un
sch\'ema en une cat\'egorie de syst\`emes de
r\'ealisations.
Nous le noterons ici,
apr\`es omission de l'aspect cristallin,
$P_{y,x}^{\real}$.
Comme pr\'ec\'edemment, les $P_{y,x}^{\real}$
forment un groupo{\umi}de.

\proclaim{Th\'eor\`eme 4.13}
Si la vari\'et\'e $X$ est unirationnelle,
$P_{y,x}^{\real}$ est motivique,
r\'ealisation d'un $\MAT(k)$-sch\'ema.
\endproclaim

Les cat\'egories de syst\`emes de
r\'ealisations consid\'er\'es sont des champs
sur $\Spec(k)_{\et}$.
Le th\'eor\`eme \'equivaut donc \`a dire qu'il
existe une extension galoisienne $k_1$ de $k$
telle qu'apr\`es extension du corps de base
\`a $k'$ dans $\scrC(k_1/k)$,
$P_{y,x}^{\real}(X_{k'})$ est motivique,
r\'ealisation d'un $\MT(k')$-sch\'emas.

\demo{\bf Preuve pour $X$ un ouvert de $\dbP^n$}
Soit $D\Subset\dbP^n$ une droite, assez
g\'en\'erale pour que le morphisme de
sch\'emas en groupe unipotents
$$
\pi_1^{\real}(X\CAP D,x)\to
\pi_1^{\real}(X,x)
$$
($x$ quelconque dans $X\CAP D$) soit surjectif.
\enddemo

Si $X\CAP D$ est le compl\'ement dans $D$ d'un
ensemble fini de points rationnels,\break
$\pi_1^{\real}
(X\CAP D,x)$ est dans $\scrR^{H+\ell}$ et,
d'apr\`es 4.4, motivique, r\'ealisation
de $\pi_1^{\mot}(X\CAP D,x)$ dans $MT(k)$.
L'alg\`ebre affine de $\pi_1^{\real}(X,x)$ est
une sous-alg\`ebre de celle de $\pi_1^{\real}(X
\CAP D,x)$, donc est \'egalement motivique,
dans $MT(k)$.
En g\'en\'eral, pour $k'$ assez grand (i.e.,
dans un $\scrC(k_1/k)$ convenable), $X\CAP D$
deviendra apr\`es extension du corps de base
de $k$ \`a $k'$ le compl\'ement d'un ensemble de
points rationnels sur $k'$, et la construction
ci-dessus nous fournit
$\pi_1^{\mot}(X_{k'},x)$ dans $MT(k')$,
fonctoriel en $k'$, i.e. $\pi_1^{\mot}(X,x)$
dans $\MAT(k)$.

Soit $y\in X(k)$, distinct de $x$.
Prouvons que $P_{y,x}^{\real}$ est motivique.
Soit $D$ la droite de $\dbP^n$ passant par $x$
et $y$.
Comme pr\'ec\'edemment, on se ram\`ene \`a
supposer que $X\CAP D$ est le compl\'ement, dans
$D$, d'un ensemble fini de points rationnels.
Le $\scrR^{H+\ell}$-sch\'ema
$P_{y,x}^{\real}(X\CAP D)$ est alors la
r\'ealisation de $P_{y,x}^{\mot}(X\CAP D)$
dans $\MT(k)$.
On obtient $P_{y,x}^{\mot}(X)$ en poussant le
$\pi_1(D\CAP X,x)$-torseur
$P_{y,x}^{\mot}(D\CAP X)$ par
$$
\pi_1(D\CAP X,x)\to\pi_1(X,x).
$$

Pour $y,z\in X(k)$, $P_{y,z}^{\mot}$ est
d\'eduit
des $\pi_1^{\mot}(X,x)$-torseurs
$P_{y,x}^{\mot}$ et $P_{z,x}^{\mot}$.

\demo{\bf Preuve de $4.13$ {\rm (cas g\'en\'eral)}}
La vari\'et\'e $X$, \'etant unirationelle, est
domin\'ee par un ouvert $U$ de $\dbP^n$.
Fixons $x\in U(k)$.
Pour $\sigma$ un plongement complexe de $k$,
le morphisme
$\pi_1(U(\dbC),x)\to\pi_1(X(\dbC),x)$ a une
image d'indice fini.
Le morphisme induit sur les enveloppes
alg\'ebriques unipotentes est donc surjectif,
et l'alg\`ebre affine de $\pi_1^{\real}(X,x)$
est motivique en tout que sous-alg\`ebre de
l'alg\`ebre affine de $\pi_1^{\real}(X,x)$,
que nous savons d\'ej\`a \^{e}tre motivique.
Le cas des $P_{y,x}$ se traite ensuite comme
plus haut, en utilisant que deux points de $X$
peuvent \^{e}tre connect\'es par une 
cha\^{\dotless\char'020}ne de courbes
rationnelles.
\enddemo

\remark{\bf Remarque 4.14}
Utilisant 4.4, on pourrait aussi dans 4.13
prendre pour $x$ et $y$ des points-base \`a
l'infini, au sens de Deligne (1989) \S{15}.
\endremark

\proclaim{Proposition 4.15}
Supposons que $\Xbar$ soit une
compactification normale de la vari\'et\'e
unirationelle $X$, et soient $x,y\in X$.
Si les composantes irr\'eductibles de
codimension un dans $X$ de $Y:=\Xbar-X$ sont
absolument irr\'eductibles, alors
$P_{y,x}^{\mot}$ est un $\MT(k)$-sch\'ema.
\endproclaim

\demo{Preuve}
Supposons d'abord que $x=y$.
Pour que $\pi_1(X,x)$ soit un
$\MT(k)$-sch\'ema il suffit que
$\Lie\,\pi_1(X,x)$ soit dans $\MT(k)$.
D'apr\`es 2.18, il suffit m\^{e}me que son
gradu\'e par le poids le soit.
Ce gradu\'e \'etant engendr\'e par $H_1(X)$,
purement de poids $-2$, il suffit que l'action
de $\Gal(\kbar/k)$ sur $\omega H_1(X)$ soit
triviale.
Ce groupe est
$$
\scrO^*(X_{\kbar})/\kbar^*\otimes\dbQ
$$
et ses \'el\'ements sont d\'etermin\'es par
leurs valuations le long des composantes
irr\'eductibles de codimension un de
$Y_{\kbar}$.
Si le groupe $\Gal(\kbar/k)$ ne les permute
pas, il agit donc trivialement sur
$\scrO^*(X_{\kbar})/\kbar^*\otimes\dbQ$.
\enddemo

Le cas g\'en\'eral r\'esulte du

\proclaim{Lemme 4.16}
Si un $\MAT(k)$-sch\'ema $P$ est un torseur
sous un $\MT(k)$-sch\'ema en groupes unipotent
$G$, il est lui-m\^{e}me un $\MT(k)$-sch\'ema.
\endproclaim

\demo{Preuve}
Le groupe $G$ \'etant unipotent, si on le fait
agir sur son alg\`ebre affine $A$ par
translations \`a droite, celle-ci admet une
filtrations exhaustive Fil par des sous-motifs
stables sous $G$, et tels que $G$ agisse
trivialement sur le gradu\'e associ\'e.
L'alg\`ebre affine $A_P$ de $P$, d\'eduite de
$A$ en poussant par le $G$-torseur $P$, admet
donc une filtration, encore not\'ee Fil, telle
que $\Gr^{\Fil}(A)\simover \Gr^{\Fil}(A_P)$.
Il reste \`a appliquer 2.18.
\enddemo

\proclaim{Proposition 4.17}
Soit $S$ un ensemble fini de places finies de
$k$ et supposons que la vari\'et\'e
unirationelle $X$ est la fibre g\'en\'erale du
compl\'ement $X_{\scrO}$, dans $\Xbar_{\scrO}$
propre et lisse sur $\Spec(\scrO_S)$, d'un
diviseur \`a croisements normaux relatif,
r\'eunion de diviseurs lisses, chacun de fibre
g\'en\'erale absolument irr\'eductible.
Si $x$ et $y$ dans $X(k)$ proviennent de
$x_{\scrO}$, $y_{\scrO}$ dans
$X_{\scrO}(\scrO_S)$, alors $P_{y,x}^{\mot}$
est dans $\MT(\scrO_S)$.
\endproclaim

\demo{Preuve}
Pour tout nombre premier $\ell$, l'hypoth\`ese
assure que la r\'ealisation $\ell$-adique de
$P_{y,x}^{\mot}$ est non ramifi\'ee en dehors
de $S$ et des places au-dessus de $\ell$, et on
applique 1.5.
\enddemo

\remark{\bf 4.18 Remarque}
Un \'enonc\'e analogue vaut pour $x$ ou $y$ \`a
l'infini.
\endremark

\subhead
4.19
\endsubhead
Soit $X$ une vari\'et\'e lisse unirationnelle
sur un corps de nombres $k$.
Le th\'eor\`eme 4.13 permet de d\'efini la
cat\'egorie 
$\MAT(X/k)$ des {\it syst\`emes locaux
unipotents de motifs d'Artin-Tate mixtes} sur
$X$.
On choisit $x\in X(k)$ et on d\'efinit un
objet $M$ de $\MAT(X/k)$ comme \'etant la
donn\'ee de $M_x$ dans $\MAT(k)$ muni d'une
action de $\pi_1(X,x)$.
Il y a lieu de voir $M_x$ comme \'etant la
fibre de $M$ en $x$.
Le choix de $x$ n'importe pas: pour $y$ dans
$X(k)$, on d\'efinit $M_y$ comme \'etant le
tordu de $M_x$ par le $\pi_1(X,x)$-torseur
$P_{y,x}$.
Il reviendrait au m\^{e}me de d\'efinir
$\MAT(X/k)$ comme la cat\'egorie des 
repr\'esentations, dans $\MAT(k)$, du
groupo{\umi}de des $P_{y,x}$ pour $x,y\in
X(k)$.
On pourrait aussi prendre des points base \`a
l'infini.

Puisque $\Lie\,\pi_1(X,x)$ est en poids $<0$,
$\pi_1(X,x)$ agit trivialement sur
$\Gr^W(M_x)$ et, pour $x$ et $y$ deux  points
base de $X$, $\Gr^W(M_x)$ et $\Gr^W(M_y)$ sont
canoniquement isomorphes.
Le foncteur fibre $\omega(M):=\omega(M_x)$ est
donc 
ind\'ependant de $x$.
On d\'efinit la sous-cat\'egorie $\MT(X)$ de
$\MAT(X/k)$ par la condition ``$M_x$ est dans
$\MT(k)$''.
Le choix de $x$ n'importe pas (2.18).

\subhead
4.20
\endsubhead
Soit $\sigma$ un plongement de $k$ dans une
cl\^{o}ture alg\'ebrique $C$ de $\dbR$.
Pour $M$ dans $\MAT(X/k)$, $(M_x)_\sigma$ est
une structure de Hodge-Tate mixte, munie d'une
action de Hodge mixte de $\pi_1(X,x)_\sigma$.
D'apr\`es Hain and Zucker (1987) (5.21),
ces donn\'ees d\'efinissent
une variation admissible $M_\sigma$ de
structures de Hodge mixtes sur $X(C)$.

Le torseur $(P_{y,x})_\sigma$ tord $(M_\sigma)_x$ 
en $(M_\sigma)_y$ et la variation $M_\sigma$ ne 
d\'epend donc pas du choix de $x$.

De m\^{e}me, $M$ d\'efinit un
$\dbQ_\ell$-faisceaux lisse sur $X$.

\subhead
4.21
\endsubhead
R\'eciproquement, supposons donn\'e sur $X$ un
syst\`eme de r\'ealisations $M$, candidat \`a
\^{e}tre la r\'ealisation d'un syst\`eme local
unipotent de motifs de Tate mixtes, de
gradu\'e par le poids constant.
Alors, si en un point $x$ de $X$ ce syst\`eme
est motivique, il est lui-m\^{e}me motivique.
En effet, dans les diverses r\'ealisations, il
d\'efinit une repr\'esentation de
$\pi_1(X,x)$.
Appliquant 2.14, on obtient une
repr\'esentation de $\pi_1^{\mot}(X,x)$ sur
le motif de r\'ealisation $M_x$, et cette
repr\'esentation d\'efinit le syst\`eme local
voulu.

\subhead
4.22
\endsubhead
Si $K$ est le corps des fonctions rationnelles
sur $X$, il est raisonnable de d\'efini
$\MT(K)$ comme la limite inductive des
cat\'egorie $\MT(U)$, pour $U$ un ouvert non
vide de plus en plus petit de $X$.
La relation entre cette cat\'egorie et la
cat\'egorie d\'eriv\'ee motivique $\DM(K)$ de
Voevodsky n'est pas claire.
Faute de disposer de la conjecture
d'annulation de Beilinson-Soul\'e pour $K$, on
ne sait pas extraire de $\DM(K)$ une
cat\'egorie ab\'elienne de motifs de Tate
mixtes sur $K$.
Pour $X$ rationnelle de dimension un, la
conjecture d'annulation est vraie, mais
nous n'avons pas prouv\'e l'\'equivalence des
deux constructions.
Pour la cat\'egorie $\MT(K)$ d\'efinie ici, on
a
$$
\alignat2
&\Ext^1(\dbQ(0),\dbQ(1))=K^*\otimes\dbQ
  &\qquad&\text{\rm et}\\
&\Ext^1(\dbQ(0),\dbQ(n))=K_{2n-1}(k) 
  &\qquad&\text{\rm pour $n\Ge 2$}.
\endalignat
$$

\subhead
4.23 Remarque
\endsubhead
Soit $\dbQbar$ une cl\^{o}ture alg\'ebrique de
$\dbQ$ et d\'efinissons $\MT(\dbQbar)$ 
comme la limite inductive des $\MT(F)$, pour
$F$ une extension finie de $\dbQ$ dans
$\dbQbar$.
Goncharov (1994) conjecture que, dans
cette limite, les $P_{y,x}(\dbP^1-S)$
engendrent tous les motifs de Tate
mixte.
En d'autres termes: que tout motif de Tate mixte
sur $\dbQbar$ est un sous-quotient d'un
produit tensoriel de sous-motifs d'alg\`ebres
affines $\scrO(P_{y,x}(\dbP^1-S))$ et de duaux
de tels sous-motifs.

Les remarques suivantes montrent qu'il est en
tout cas difficile de sortir de la
sous-cat\'egorie tannakienne $\scrT$ de
$\MT(\dbQbar)$ engendr\'ee par les
$P_{y,x}(\dbP^1-S)$.

\medskip\noindent
{\sans (a)}\enspace
Pour $X$ lisse unirationnelle, les
$P_{y,x}(X)$ sont dans $\scrT$.
Cela r\'esulte de la preuve de 4.13.

\smallskip\noindent
{\sans (b)}\enspace
Si $M$ est un syst\`eme local unipotent de
motifs de Tate mixte sur $X$ lisse
unirationnelle, et si en un point $M$ est dans
$\scrT$, alors $M$ est en chaque point dans
$\scrT$.
En effet, $M_y$ est tordu de $M_x$ par
$P_{y,x}$.
Il se plonge donc dans
$M_x\otimes\scrO(P_{y,x})$, et est dans
$\scrT$ puisque $M_x$ et $\scrO(P_{y,x})$ le
sont.

\newpage


\dspace
\subhead
5. Le groupe fondamental motivique du
compl\'ement, dans $\dbP^1$, de $0$, $\infty$,
et $\mu_N$
\endsubhead

\subhead
5.1
\endsubhead
Soient $N$ un entier $\Ge 1$,
$k$ engendr\'e sur $\dbQ$ par une racine
primitive $N^{\ieme}$ de l'unit\'e
et $\scrO$ l'anneau de ses entiers.
Il est commode de ne pas supposer que $k$ soit
le sous-corps $\dbQ(\exp(2\pi i/N))$ de
$\dbC$.
Notons $A_k$, ou simplement $A$ la droite
affine type sur $k$ et posons
$A^*=A-\{0\}=\Spec\,k[u,u^{-1}]$.
Soit $X$ le compl\'ement dans $A^*$ de
l'ensemble $\mu_N:=\mu_N(k)$ des racines
$N^{\iemes}$ de l'unit\'e.
S'il y a ambigu{\umi}t\'e sur $N$, on \'ecrira
plut\^{o}t $k_N$, $\scrO_N$ et $X_N$.

Le groupe dih\'edral $\dbZ/2\ltimes\mu_N$ agit
sur $X\Subset \dbP_k^1$: le g\'en\'erateur de
$\dbZ/2$ agit par $u\mapsto u^{-1}$, et
$\xi\in\mu_N$ agit par $u\mapsto\xi u$.
Pour $N=1,2$ ou $4$, le groupe des
projectivit\'es transformant $X$ en
lui-m\^{e}me est plus grand.
Pour $N=1$, il s'identifie au groupe des
permutations de $\{0,1,\infty\}$.
Pour $N=2$, $(0,\infty,1,-1)$ est un quaterne
harmonique, et on obtient le groupe des
automorphismes d'un carr\'e de sommets
cons\'ecutifs $0$, $1$, $\infty$, $-1$.
Pour $N=4$, les points $0$, $1$, $i$, $-1$, $-i$,
$\infty$ forment sur la sph\`ere de Riemann les
sommets d'un octa\`edre, et on obtient le
groupe des d\'eplacements de cet octa\`edre.

\subhead
5.2
\endsubhead
Si $x$ et $y$ sont deux points-base (4.3) de
$X$, l'espace de chemins motivique $P_{y,x}$
est un $\MT(k)$-sch\'ema affine (4.4).
La r\'ealisation $\omega$ de $P_{y,x}$ est
ind\'ependante de $x$ et $y$.
C'est un $\dbQ$-sch\'ema en groupe pro-unipotent, 
et la composition des chemins est
la loi de groupe.
Autre fa\c{c}on de dire: en r\'ealisation
$\omega$, il y a un chemin canonique de $x$
\`a $y$, et le compos\'e des chemins
canoniques de $x$ \`a $y$ et de $y$ \`a 
$z$ est le chemin canonique de $x$ \`a $z$.
Bien s\^{u}r, cette trivialisation du
$\omega(\pi_1(X,x))$-torseur $\omega(P_{y,x})$
n'est en g\'en\'eral pas motivique, i.e. ne
provient pas d'une trivialisation du
$\pi_1(X,x)$-torseur $P_{y,x}$.
L'action sur $\omega(P_{y,x})$
du groupe des automorphismes du foncteur fibre
$\omega$ d\'epend de $x$ et $y$.

\subhead
5.3
\endsubhead
Pour $z$ dans $\dbP_k^1-\{\infty\}$ et
$\lambda$ dans $k$ nous noterons
$\lambda_z$ le vecteur tangent en $z$ de
coordonn\'ee
$du(\lambda_z)=\lambda$.
Nous noterons $\lambda_\infty$ le vecteur
tangent en $\infty$ de coordonn\'ee
$du^{-1}(\lambda_\infty)=\lambda$.

Nous prendrons comme points-base les
$\lambda_z$ pour $z=0$, $\infty$ ou dans
$\mu_N$ et $\lambda$ une racine de l'unit\'e
dans $k$.
Ce syst\`eme de points-base est stable sous
l'action du groupe dih\'edral
$\dbZ/2\ltimes\mu_N$.
Une application de 1.7 montre comme en 4.17,
4.18 que les $P_{y,x}$ correspondants sont non
ramifi\'es en dehors de 
l'ensemble des places de $k$ au-dessus
d'un nombre premier divisant $N$: ce sont
des $\MT(\scrO[1/N])$-sch\'emas.
Nous d\'ecrivons en 5.4 et 5.5 des structures
sur le syst\`eme des $P_{y,x}$.
Pour toutes ces structures, il suffit pour les
construire de les construire dans
$\scrR^{H+\ell}$, ce qui est fait dans
Deligne (1989), et d'appliquer 2.14.

\subhead
5.4
\endsubhead
Soient $T$ une droite (sch\'ema vectoriel de
dimension un) et $T^*=T-\{0\}$.
Le $\pi_1$ motivique de $T^*$ est ab\'elien,
ind\'ependant du point-base.
C'est $\dbQ(1)$.
Pour $x\in T^*(k)$, l'application
$\lambda\mapsto\lambda x$ est un isomorphisme
de $A^*$ avec $T^*$.
Pour $y\in T^*(k)$, il induit un isomorphisme
de $\dbQ(1)$-torseurs de
$P_{y/x,1}(A^*)$ avec $P_{y,x}(T^*)$.
Le $\dbQ(1)$ torseur $P_{t,1}(A^*)$
est le {\it torseur de Kummer} $K(t)$, et la
composition des chemins fournit un
isomorphisme de torseurs $K(tu)=K(t)+K(u)$.

En particulier, si $\alpha\in k^*$ est une
racine $n^{\ieme}$ de l'unit\'e, elle fournit
une trivialisation de $nK(\alpha)$ et, puisque la
multiplication par $n$ est un automorphisme de
$\dbQ(1)$, par division, une trivialisation de
$K(\alpha)$.

Pour $x\in T^*(k)$ et $\alpha$ une racine de
l'unit\'e, le $\dbQ(1)$-torseur $P_{\alpha
x,x}$ est ainsi trivialis\'e: on dispose d'un
chemin motivique de $x$ \`a $\alpha x$.

Prenons pour $T$ l'espace tangent \`a
$\dbP^1$ en $0$, $\infty$ ou en
$\zeta\in\mu_N$.
Le groupo{\umi}de des $P_{y,x}(T^*)$, pour
$x,y\in T^*(k)$, s'envoie dans le 
groupo{\umi}de
des $P_{y,x}(X)$.
Pour $x\in T^*(k)$ et $\alpha$ une racine de
l'unit\'e on dispose donc d'un chemin
motivique de $x$ \`a $\alpha x$ dans $X$: une
trivialisation du $\pi_1(X,x)$-torseur
$P_{\alpha x,x}(X)$.

Pour le syst\`eme de points base $\lambda_z$
de 5.3, $P_{\lambda_z,\mu_t}$ ne d\'epend donc
que de $z$ et $t$, et on \'ecrira simplement
$P_{z,t}$ pour $P_{\lambda_z,\mu_t}$.
Quand d'autres points-base \`a l'infini auront \`a
\^{e}tre utilis\'es, nous retournerons \`a la
notation $\lambda_z$.

\subhead
5.5
\endsubhead
On dispose sur le syst\`eme des
$\MT(\scrO[1/N])$-sch\'emas $P_{y,x}$ ($x,y\in
S:=\{0,\infty\}\Cup \mu_N$) des structures
suivantes.

\medskip\noindent
{\sans (A)}\enspace
La composition des chemins.

\smallskip\noindent
{\sans (B)}\enspace
Pour chaque $x\in S$, un 
morphisme de sch\'emas en groupe, la {\it
monodromie locale autour de $x$}:
$$
\dbQ(1)\to P_{x,x}.
$$

\smallskip\noindent
{\sans (C)}\enspace
une \'equivariance sous le groupe dih\'edral
$\dbZ/2\ltimes\mu_N$ de 5.1.

\subhead
5.6
\endsubhead
Si on applique un foncteur fibre $F$, \`a
valeurs dans les espaces vectoriels sur un
corps $K$, \`a $\dbQ(1)$, aux $P_{y,x}$ et 
aux structures 5.5, on obtient une structure
de l'esp\`ece suivante:

\medskip\noindent
(1)\enspace
Un espace vectoriel de dimension un $K(1)$.

\noindent
(2)\enspace
Pour chaque $x,y\in S$, un sch\'ema $P_{y,x}$
sur $K$.

\noindent
$(3)_{\text{A}}$\enspace
Un syst\`eme de morphismes de sch\'emas
$$
P_{z,y}\times P_{y,x}\to P_{z,x}\tag5.6.1
$$
faisant des $P_{y,x}$ un groupo{\umi}de.
Les sch\'emas en groupe $P_{x,x}$ sont
pro-unipotents.

\noindent
$(3)_{\text{B}}$\enspace
Pour chaque $x\in S$, un homomorphisme
$$
\text{(groupe additif }
K(1))\to P_{x,x}.\tag5.6.2
$$
L'alg\`ebre de Lie de $P_{x,x}$ est une
pro-alg\`ebre de Lie.
Parce que $P_{x,x}$ est pro-unipotent, la
donn\'ee (5.6.2) \'equivaut \`a celle de
$$
K(1)\to \Lie\,P_{x,x}.\tag5.6.3
$$

\noindent
$(3)_{\text{C}}$\enspace
Une $\dbZ/2\ltimes\mu_N$ \'equivariance.

\subhead
5.7
\endsubhead
Il nous sera utile d'affaiblir une structure
5.6 en restreignant les points-base \`a
\^{e}tre $0$ ou dans $\mu_N$.
Plus pr\'ecis\'ement, \`a une structure 5.6,
on peut attacher une structure affaiblie du
type suivant.

\medskip\noindent
(1)\enspace
Comme en 5.6 (1).

\noindent
(2)\enspace
Un sch\'ema en groupe pro-unipotent $P_{0,0}$
et, pour $\zeta\in\mu_N$, un $P_{00}$-torseur
$P_{\zeta,0}$.
On note $P_{\zeta\zeta}$ le tordu de $P_{00}$
par ce torseur, i.e. le sch\'ema en groupe des
automorphismes du $P_{0,0}$-torseur
$P_{\zeta,0}$.

\noindent
$(3)_{\text{B}}$\enspace
Pour $x=0$ ou $x\in\mu_N$, un morphisme
(5.6.2) ou, ce qui revient au m\^{e}me,
(5.6.3).

\noindent
$(3)_{\text{C}}$\enspace
Une $\mu_N$-\'equivariance.

\remark{\bf Remarques} (i)
La donn\'ee d'une structure 5.7 \'equivaut \`a
celle de (1) et de

\noindent
$(2')$\enspace
Un sch\'ema en groupe pro-unipotent $P_{00}$
et un $P_{00}$-torseur $P_{10}$.
Comme en (2), $P_{10}$ d\'efinit par torsion
$P_{11}$.

\noindent
$(3')_{\text{B}}$\enspace
Comme en $(3)_{\text{B}}$, pour $x=0$ ou $1$.

\noindent
$(3')_{\text{C}}$\enspace
Une action sur $P_{00}$ de $\mu_N$, respectant
la monodromie locale autour de $0$.

On r\'ecup\`ere (2) en d\'efinissant
$P_{\zeta,0}$ comme une copie de $P_{1,0}$,
l'\'equivariance $(3)_{\text{C}}\,:\zeta: P_{1,0}\to
P_{\zeta,0}$ \'etant l'identit\'e de $P_{1,0}$
sur cette copie.

\noindent
(ii)\enspace
Cette structure est un peu plus pr\'ecise que
la donn\'ee de $P_{00}$, de l'action de
$\mu_N$, de $(3)_{\text{B}}$ pour $x=0$ et de
la classe de conjugaison ``monodromie locale
autour de $1$'' de morphismes $\dbQ(1)\to
P_{1,1}$: cette classe de conjugaison est un
espace homog\`ene, qu'on peut identifier au
quotient $P_{0,1}/\dbQ(1)$.
\endremark

\subhead
5.8
\endsubhead
Prenons pour foncteur fibre le foncteur fibre
$\omega$, et d\'ecrivons la structure
d'esp\`ece (5.6) obtenue.
Puisque $\omega$ est \`a valeurs dans les
$\dbQ$-espaces vectoriels gradu\'es, les
alg\`ebres affines des $\omega(P_{y,x})$ seront
des $\dbQ$-alg\`ebres gradu\'ees.

Soit $\scrL$ l'alg\`ebre de Lie engendr\'ee
par des g\'en\'erateurs $e_x$, 
$x\in\{0,\infty\}\Cup\mu_N$, soumis \`a la seule
relation que $\sum e_x=0$.
En cas d'ambigu{\umi}t\'e, on l'\'ecrira plus
pr\'ecisement $\scrL_N$.
On la gradue en donnant aux $e_x$ le degr\'e un.
Elle est librement engendr\'ee par les $e_x$,
$x\not=\infty$.

Soit $\prod$ le groupe pro-unipotent
$$
\prod:=\lim\,\exp(\scrL/\text{degr\'e $\Ge
n$)}.
$$
On munit son alg\`ebre affine de la graduation
d\'eduite de celle de $\scrL$.
L'alg\`ebre de Lie de $\prod$ est le
compl\'et\'e $\scrL^{\uphat}
:=\prod\limits_{n}\scrL_n$ de $\scrL$.

Avec ces notations, la structure d'esp\`ece
5.6, et les graduations, sont comme suit.

\medskip\noindent
$(1)^\omega$\enspace
L'espace vectoriel $\dbQ$, en degr\'e un.

\noindent
$(2)^\omega$\enspace
Pour chaque $x$, $y$, une copie
$\prod\nolimits_{y,x}$ de $\prod$.

\noindent
$(3)_{\text{A}}^\omega$\enspace
Morphismes (5.6.1): la loi de groupe de
$\prod$.

\noindent
$(3)_{\text{B}}^\omega$\enspace
Le morphisme (5.6.3) est
$$
\dbQ\to\scrL^{\uphat}\colon\, 1\longmapsto e_x
$$

\noindent
$(3)_{\text{C}}^\omega$\enspace
Pour $\sigma\in\dbZ/2\ltimes\mu_N$, agissant
sur $\{0,\infty\}\Cup\mu_N$ comme en 5.1, $\sigma$
envoie $\prod\nolimits_{y,x}$ sur
$\prod\nolimits_{\sigma y,\sigma x}$.
C'est le morphisme de groupes
$[\sigma]\colon\,\prod\to\prod$
qui induit sur l'alg\`ebre de Lie $e_x\mapsto
e_{\sigma x}$.

D'apr\`es 2.10, la m\^{e}me description vaut,
sur $k$, pour le foncteur fibre $\DR$.

Le groupe $G_\omega\left<\MT(\scrO[1/N])\right>$ 
(2.1) agit sur cette structure d'esp\`ece 5.8.
Nous nous proposons d'\'etudier cette action.

\subhead
5.9
\endsubhead
La structure 5.8 d'esp\`ece 5.6 fournit par
oubli une structure d'esp\`ece 5.7.
Pr\'esent\'ee comme en 5.7 Remarque (i), c'est
la suivante.

\medskip\noindent
$(1)^\omega$\enspace
L'espace vectoriel $\dbQ$.

\noindent
$(2)^\omega$\enspace
Une copie $\prod\nolimits_{0,0}$ de $\prod$ et
le $\prod\nolimits_{0,0}$ torseur trivial
$\prod$, not\'e $\prod\nolimits_{1,0}$.
Le tordu de $\prod\nolimits_{0,0}
$ par ce torseur est une nouvelle copie de
$\prod$, not\'ee $\prod\nolimits_{1,1}$.

\noindent
$(3)_{\text{B}}^\omega$\enspace
Les morphismes (5.6.3), pour $x=0$ ou $1$ sont
$\dbQ\to\scrL^{\uphat}\colon\, 1\to e_x$.

\noindent
$(3)_{\text{C}}^\omega$\enspace
L'action de $\mu_N$ sur
$\prod\nolimits_{0,0}$, donn\'ee sur
l'alg\`ebre de Lie $\scrL^{\uphat}$  par
$e_x\mapsto e_{\xi x}$.

\subhead
5.10
\endsubhead
Soit $H_\omega$ le sch\'ema en groupe des
automorphismes de la structure 5.9 d'esp\`ece
5.7.
L'action de $H_\omega$ sur l'espace vectoriel
de dimension un $(1)^\omega$ est un morphisme
$$
H_\omega\to \dbG_m.\tag5.10.1
$$
Soit $V_\omega$ son noyau.
Les graduations fournissent un morphisme de
$\dbG_m$ dans $H_\omega$.
Puisque $(1)^\omega$ est en degr\'e un, c'est
une section de (5.10.1).
Elle fait de $H_\omega$ un produit
semi-direct:
$$
H_\omega=\dbG_m\ltimes V_\omega.\tag5.10.2
$$

L'action de $G_\omega$ sur la structure 5.9 se
factorise par $H_\omega$, puisque $H_\omega$ a
\'et\'e d\'efini comme le groupe
d'automorphismes de r\'ealisations $\omega$
respectant des structures motiviques.
Comparant \`a 2.1, on v\'erifie que le
morphisme
$$
G_\omega=\dbG_m\ltimes U_\omega\to H_\omega
=\dbG_m\ltimes V_\omega\tag5.10.3
$$
est compatible aux d\'ecompositions en
produits semi-directs.

\proclaim{Proposition 5.11}
Notons $1_{y,x}$ l'\'el\'ement neutre de la
copie $\prod\nolimits_{y,x}$ de $\prod$.
Avec les notations de 5.9 et 5.10, le
morphisme $v\mapsto v(1_{1,0})$ est un
isomorphisme de sch\'emas de $V_\omega$ avec
$\prod\nolimits_{1,0}$.
\endproclaim

Rappelons que pour $\xi\in\mu_N$, 
$[\xi]$ est l'automorphisme de $\prod$
qui induit sur l'alg\`ebre de Lie
$e_x\mapsto e_{\xi x}$.
On notera $\prod\nolimits_{0,1}$ le
$\left(\prod\nolimits_{1,1},
\prod\nolimits_{0,0}\right)$-bitorseur inverse
du $\left(\prod\nolimits_{0,0},
\prod\nolimits_{1,1}\right)$-bitorseur
$\prod\nolimits_{1,0}$.
C'est une nouvelle copie $\prod$.
Pour $x,y,z\in\{0,1\}$, on dispose d'une
``composition des chemins''
$\prod\nolimits_{z,y}\times\prod\nolimits_{y,x}\to
\prod\nolimits_{z,x}$, simplement donn\'ee par
la loi de groupe de $\prod$.

\demo{Preuve}
Soit $v$ dans $V_\omega$ et posons
$a:=v(1_{1,0})$.
\enddemo

Puisque $V_\omega$ agit trivialement sur
$(1)^\omega$ et respecte
$(3)_{\text{B}}^\omega$, on a dans
$\Lie\,\prod\nolimits_{00}$, identifi\'e \`a
$\scrL^{\uphat}$
$$
v(e_0)=e_0\qquad\qquad
\text{(dans $\Lie\,{\tsize \prod\nolimits_{0,0})}$}
\tag5.11.1
$$
et dans $\Lie\,\prod\nolimits_{1,1}$,
\'egalement identifi\'e \`a $\scrL^{\uphat}$
$$
v(e_1)=e_1\qquad\qquad\text{(dans
$\Lie\,{\tsize \prod\nolimits_{1,1}}$)}.
$$
Transformant par $v$ l'identit\'e
$$
\left(\exp(e_1)\text{ dans }{\tsize
\prod\nolimits_{0,0}}\right)=1_{0,1} \cdot \left(\exp(e_1)
\text{ dans }{\tsize\prod\nolimits_{11}}\right)\cdot 
1_{1,0},
$$
on obtient
$$
v(\exp(e_1)\text{ dans }{\tsize \prod\nolimits_{0,0}})
=a^{-1}\exp(e_1)a,
$$
qui \'equivaut \`a
$$
v(e_1\text{ dans }\Lie\,
{\tsize\prod\nolimits_{0,0}})=
\ad_{a^{-1}}(e_1).\tag5.11.2
$$

Par \'equivariance, (5.11.2) implique que pour
$\zeta\in\mu_N$,
$$
v(e_\zeta\text{ dans }
\Lie\,{\tsize
\prod\nolimits_{0,0}})=\ad_{[\zeta]a^{-1}}
(e_\zeta).\tag5.11.3
$$

Pour $a$ dans $\prod$, notons
$\left<a\right>_0$ l'automorphisme de $\prod$
induisant sur l'alg\`ebre de
$\Lie\,\scrL^{\uphat}$ l'automorphisme, encore
not\'e $\left<a\right>_0$:
$$
\left<a\right>_0\colon\, e_0\longmapsto
  e_0,\quad
e_\zeta\longmapsto
\ad_{[\zeta]a^{-1}}(e_\zeta).
\tag5.11.4
$$
Notons $\left<a\right>_{1,0}$ l'automorphisme du
sch\'ema $\prod$
$$
\left<a\right>_{1,0}\colon\, g\longmapsto 
a \cdot \left<a\right>_0(g).\tag5.11.5
$$
Les formules (5.11.1) et (5.11.3) montrent que
$v$ agit par $\left<a\right>_0$ sur
$\prod\nolimits_{0,0}$. 
Sur le $\prod\nolimits_{0,0}$-torseur
$\prod\nolimits_{1,0}$, $v$ agit par
$\left<a\right>_{1,0}$: transformer par $v$
l'identit\'e ($g$ dans $\prod\nolimits_{1,0}$)
$=1_{1,0} \cdot (g\text{ dans
}\prod\nolimits_{0,0})$.

R\'eciproquement, quelque soit $a$ dans $\prod$,
ces formules d\'efinissent un automorphisme
$\left<a\right>$ de la structure 5.9
d'esp\`ece 5.7.

\subhead
5.12
\endsubhead
Si on transporte \`a la copie
$\prod\nolimits_{1,0}$ de $\prod$ la loi de
groupe de $V_\omega$, on obtient une nouvelle
loi de groupe $\circ$ sur le sch\'ema $\prod$.
Elle est donn\'ee par
$$
a \circ b=a \cdot \left<a\right>_0(b).\tag5.12.1
$$
En effet, $\left<b\right>_{1,0}$ envoie
$1_{1,0}\in\prod\nolimits_{1,0}$ sur $b$ et
$\left<a\right>_{1,0}$ envoie
$b=1_{1,0} \cdot b\in\prod\nolimits_{1,0}$ sur
$a \cdot \left<a\right>_0(b)$.

Prenant l'espace tangent \`a l'origine, on
d\'eduit de 5.11 un isomorphisme de pro-espaces
vectoriels de $\Lie\,V_\omega$ avec
$\scrL^{\uphat}$.
Transportons \`a $\scrL^{\uphat}$ l'action de
$\Lie\,V_\omega$ d\'eduite de l'action de
$V_\omega$ sur la structure 5.9.
D\'erivant (5.11.4), on trouve que 
$a\in \scrL^{\uphat}$ agit sur
$\scrL^{\uphat}=\Lie
\left(\prod\nolimits_{0,0}\right)$ par
$$
\partial_a\colon\, e_0\longmapsto 0,\quad
e_\zeta\longmapsto \left[-[\zeta](a),e_\zeta\right].
\tag5.12.2
$$
L'action sur $\prod\nolimits_{00}$ s'en
d\'eduit.
Il est clair sur (5.12.2) que
$\partial_{e_1}=0$.
Puisque le centralisateur de $e_1$ dans
$\scrL^{\uphat}$ est r\'eduit aux multiples de
$e_1$ (4.7), le noyau de $a\mapsto\partial_a$ est
r\'eduit \`a $\dbQ e_1$.

Pour \'ecrire l'action de
$\Lie\,V_\omega=\scrL^{\uphat}$ tant sur
$\prod\nolimits_{00}$ que sur le
$\prod\nolimits_{00}$-torseur
$\prod\nolimits_{1,0}$, il est commode de
plonger ces sch\'emas affines dans le dual de
leur alg\`ebre affine, un sch\'ema
provectoriel dont ils sont le sous-sch\'ema
des \'el\'ements groupaux.
A la structure de groupe de
$\prod\nolimits_{0,0}$ correspond une
structure d'alg\`ebre sur ce dual, \`a celle
de $\prod\nolimits_{0,0}$-torseur de
$\prod\nolimits_{0,1}$ une structure de module
\`a droite sous cette alg\`ebre:
$$
\spreadmatrixlines{2\jot}
\matrix
{\tsize\prod\nolimits_{0,1}}  &:\text{torseur
sous}\hfill &{\tsize\prod\nolimits_{0,0}} 
  &,\text{\'el\'ements groupaux de}\\
\null\kern-8pt\CAP &&\null\kern-8pt\CAP &\\ 
\dbQ\ll e_0,(e_\zeta)\gg 
   &:\text{module \`a droite sur}\hfill 
  &\dbQ\ll e_0,(e_\zeta)\gg &.
\endmatrix
$$

\medskip
Si $\dbQ[\eps]$ est l'alg\`ebre des
nombres duaux $(\eps^2=0)$, l'action de
$a\in\scrL^{\uphat}=\Lie(V_\omega)$ est le
coefficient de $\eps$ dans l'action du
$\dbQ[\eps]$-point $1+a\eps$ de $V_\omega$.
L'action sur l'alg\`ebre $\dbQ\ll
e_0,(e_\zeta)\gg$, duale de
$\scrO(\prod\nolimits_{0,0})$, est $e_0\mapsto
e_0$, $e_\zeta\mapsto
(1+[\zeta]a\eps)^{-1}e_\zeta(1+[\zeta]a\eps)$,
et correspond \`a l'action de a par la
d\'erivation $\partial_a\colon\, e_0\to 0$,
$e_\zeta\mapsto-[[\zeta]a,e_\zeta]$.
L'action sur le module \`a droite $\dbQ\ll
e_0,(e_\zeta)\gg$, dual de
$\scrO(\prod\nolimits_{0,1})$, est par
$x\mapsto (1+a\eps)(x+\partial_a x \cdot \eps)$ et
correspond \`a l'action de a par
$$
x\longmapsto ax+\partial_a x.\tag5.12.3
$$

Au crochet de Lie de $\Lie\,V_\omega$
correspond un nouveau crochet de Lie
$\{\,\,,\,\,\}$ sur $\scrL^{\uphat}$.

\proclaim{Proposition 5.13}
On a
$$
\{a,b\}=[a,b]+\partial_a(b)-\partial_b(a).
\tag5.13.1
$$
\endproclaim

\demo{Preuve}
V\'erifions d'abord que
$$
[\partial_a,\partial_b]=\partial_{[a,b]+\partial_a
(b)-\partial_b(a)}.\tag5.13.2
$$
Puisque, comme $\partial_a$ et $\partial_b$,
la d\'erivation $[\partial_a,\partial_b]$ de
$\scrL^{\uphat}$ fixe $e_0$ et est 
$\mu_N$-\'equivariante, il suffit de tester son action
sur $e_1$.
On a
$$
\align
\partial_a\partial_b &\colon\, e_1\longmapsto
  -[b,e_1]\longmapsto-[\partial_a(b),e_1]+
  [b,[a,e_1]]\\
\partial_b\partial_a &\colon\, e_1\longmapsto
  -[a,e_1]\longmapsto-[\partial_b(a),e_1]
  +[a,[b,e_1]]\\
[\partial_a,\partial_b] &\colon\,
e_1\longmapsto
-[\partial_a(b)-\partial_b(a)+[a,b],e_1].
\endalign
$$
\enddemo

Deux m\'ethodes pour d\'eduire (5.13.1) de
(5.13.2):

\medskip\noindent
{\sans a.}\enspace
Le crochet $\{\,\,,\,\,\}$ est compatible aux
graduations, toutes les constructions faites
\'etant compatibles \`a l'action de $\dbG_m$.
Le noyau de l'action de
$\scrL^{\uphat}=\Lie(V_\omega)$ sur $\scrL^{\uphat}$
\'etant r\'eduit aux multiples de $e_1$, si $a$
et  $b$ sont homog\`enes (de degr\'e $\Ge 1$),
il r\'esulte de (5.13.2) que les deux membres
de (5.13.1)
co{\umi}ncident \`a un multiple de $e_1$
pr\`es.
Etant homog\`ene de degr\'e $\Ge 2$, ils
co{\umi}ncident. 

\medskip\noindent
{\sans b.}\enspace
Testons sur (5.12.3).
Si $\mu_a$ est la multiplication \`a gauche
par $a$, on a
$$
\align
[\mu_a+\partial_a, &\mu_b+\partial_b] =
\mu_{[a,b]+\partial_ab-\partial_ba}+
  [\partial_a,\partial_b]\\
&=\mu_{[a,b]+\partial_ab-\partial_b a}+
\partial_{[a,b]+\partial_ab-\partial_ba}.
\endalign
$$

\remark{\bf 5.14 Remarque}\enspace
(i) On a
$$
\{a,e_1\}=[a,e_1]+\partial_a(e_1)-\partial_{e_1}(a)
=[a,e_1]-[a,e_1]+0=0:
$$
$e_1$ est dans le centre de
$(\scrL^{\uphat},\{\,\,,\,\,\})$.

\medskip\noindent
(ii)
On prendra garde que l'application
exponentielle de Lie $V_\omega$ dans
$V_\omega$, identifi\'ee par 5.11 et 5.12 \`a
une application de $\scrL^{\uphat}$ dans $\prod$,
n'est pas l'application exponentielle de
$\scrL^{\uphat}$ dans $\prod$.
Il resulte de 5.12.4 que, si $\mu_a$ est la
multiplication \`a gauche par $a$, c'est
$$
\align
\exp &\colon\, \Lie\,V_\omega\to V_\omega
  \colon\,\scrL^{\uphat}\to\prod\Subset\dbQ
\ll e_0,(e_\zeta)\gg\colon\\
&a\longmapsto\sum (\mu_a+\partial_a)^n/n!
\quad\text{appliqu\'e \`a $1$}.
\endalign
$$
Premiers termes:
$$
a+(a^2+\partial_a(a))/2+(a^3+2a \cdot \partial_a(a)+
\partial_a(a) \cdot a+\partial_a^2(a))/
6+\cdots
$$
\endremark

\remark{\bf Remarque 5.15}
Ainsi que la notation le sugg\`ere, le groupe
$H_\omega$ (resp. $V_\omega$) est la
r\'ealisation $\omega$ d'un
$\MT(\scrO[1/N])$-sch\'ema en groupe $H$ (resp. $V$).
On peut d\'efinir $H$ par la propri\'et\'e que
pour tout foncteur fibre $F$, $F(H)$ est,
fonctoriellement en $F$, le sch\'ema en
groupe des automorphismes de la structure 5.7
correspondante.
Il agit sur $\dbQ(1)$, $P_{0,0}$, et $P_{0,1}$
et $V$ est le noyau de l'homomorphisme
$H\to\dbQ(1)$ donnant l'action de $H$ sur
$\dbQ(1)$.
Il reste vrai dans $\MT(\scrO[1/N])$ que $P_{1,0}$
est un espace principal homog\`ene sous $V$.
\endremark

Comme en 5.10, le groupe
$G:=\pi(\MT(\scrO[1/N]))$ s'envoie dans $H$, et
ce morphisme induit un morphisme
$$
U\to V\tag5.15.1
$$
de $\MT(\scrO[1/N])$-sch\'emas en groupes.
Sa r\'ealisation $\omega$ est induite par
(5.10.3).

\subhead
5.16
\endsubhead
Fixons un plongement $\sigma$ de $k$ dans
$\dbC$, i.e. un isomorphisme de $\mu_N(k)$
avec $\mu_N(\dbC)$.
Le {\it droit chemin} du point base $1_0$ de
$\dbC^*-\mu_N$ vers le point base $(-1)_1$ va
par segments de droite successifs \quad
(a) dans l'espace tangent en $0$, de $1$ \`a
$\eps>0$; (b) de $\eps$ \`a $1-\eps$ dans
$\dbC^*-\mu_N$; (c) de $-\eps$ \`a $-1$ dans
l'espace tangent \`a un.
Ce droit chemin a une image $\dch$ dans
$(P_{(-1)_1,1_0})_\sigma$.
Via l'isomorphisme de comparaison
$\comp_{\sigma,\omega}$, il correspond \`a un
\'el\'ement de 
$$
\dch(\sigma)\in (P_{(-1)_1,1_0})_\omega\otimes\dbC
=\prod(\dbC) \Subset \dbC\ll
e_0,(e_\zeta)_{\zeta\in\mu_N(k)}\gg.
$$

Les coefficients de $\dch(\sigma)$ sont
donn\'es par 
$$
\spreadlines{1\jot}
\leqalignno{
&\text{coefficient de }e_{x_1}\ldots 
     e_{x_n}= &(5.16.1)\cr
&=\text{ int\'egrale it\'er\'ee de $0$ \`a $1$ de }
\frac{dz}{z-\sigma(x_1)}\cdots
\frac{dz}{z-\sigma(x_n)}\cr
&=\int_0^1\left(\frac{dz_1}{z_1-\sigma(x_1)}
\int_0^{z_1}\left(\frac{dz_2}{z_2-\sigma(x_2)}
\cdots \int_0^{z_{n-1}}
\frac{dz_n}{z_n-\sigma(x_n)}\right)\ldots
\right).\cr}
$$
Si $x_n=0$ ou si $x_1=1$, l'int\'egrale diverge
et doit \^{e}tre r\'egularis\'ee:
remplacer les limites d'int\'egration $0$ et
$1$  par $\eps$ et $1-\eta$; l'int\'egrale
it\'er\'ee convergente obtenue a pour comportement
asymptotique pour $\eps,\eta\to 0$
$$
\align
\text{Int\'egrale it\'er\'ee}
&=\text{polyn\^{o}me
$(\log\,\eps,\log\,\eta)$}\\
&\qquad +O(\sup(\eps\vert\log\,\eps\vert^A,
\eta\vert\log\,\eta\vert^B))
\endalign
$$
pour $A$ and $B$ convenables, et on prend le
terme constant du polyn\^{o}me.

Les valeurs multizeta et leur analogue
cyclotomique apparaissent quand dans (5.16.1)
on d\'eveloppe $dz/z-\zeta$ en s\'erie
g\'eom\'etrique:
$$
\frac{dz}{z-\zeta}=\frac{-\zeta^{-1}z}{1-\zeta^{-1}z}
\,\,\frac{dz}{z}=-\sum\limits_{n\Ge 1}(\zeta^{-1}z)^n.
\frac{dz}{z}.
$$

\proclaim{Proposition 5.17}
{\rm (i)}\enspace Pour $\zeta_1,\dotsc,\zeta_m$ dans
$\mu_N(k)$ et des entiers $s_1,\dotsc,s_m\Ge
1$, si $s_1\not=1$, le coefficient de
$e_0^{s_1-1}e_{\zeta_1}\ldots
e_0^{s_m-1}e_{\zeta_m}$ dans $\dch(\sigma)$ est
$$
(-1)^m\sum\limits_{n_1>\ldots>n_m>0}
\frac{\sigma\left(\zeta_1^{n_2-n_1}\zeta_2^{n_3-n_2}
\ldots \zeta_m^{-n_m}\right)}
{n_1^{s_1}\ldots n_m^{s_m}}
$$

\noindent
{\rm (ii)}\enspace
$\dch(\sigma)$ est caract\'eris\'e par les
propri\'et\'es d'\^{e}tre un \'el\'ement
groupal (terme constant 1, et 
$\Delta\dch(\sigma)=\dch(\sigma)\otimes
\dch(\sigma)$) ayant les coefficients (i) et
dont les coefficients de $e_0$ et $e_1$ sont
nuls.
\endproclaim

Notons $\tau$ l'application de $\dbG_m$ dans
$H_\omega$ d\'efinie par les graduations
(5.10.2) et $\left<g\right>$ l'\'el\'ement de
$V_\omega(\dbC)$ correspondant par 5.11 \`a
$g\in\prod(\dbC)$.

\proclaim{Proposition 5.18}
L'\'el\'ement
$\left<\dch(\sigma)\right>\tau(2\pi i)$
de $H_\omega(\dbC)$ transforme les 
$\dbQ$-formes $(P_{0,0})_\omega$ et
$(P_{1,0})_\omega$ de
$(P_{0,0})_\omega\otimes\dbC$ et
$(P_{1,0})_\omega\otimes\dbC$ en les
$\dbQ$-formes $(P_{0,0})_\sigma$ et
$(P_{1,0})_\sigma$.
\endproclaim

\demo{Preuve}
Les \'el\'ements de $\pi_1(X(\dbC),1_0)$ ont
dans $\prod(\dbC)$ une image rationnelle pour
la $\dbQ$-forme $(P_{0,0})_\sigma$.
Ce $\pi_1$ contient la monodromie locale en
$0$ et le chemin allant droit pr\`es de $1$,
tournant autour de $1$, et revenant.
Images: $\exp(2\pi i\,e_0)$ et
$\dch(\sigma)^{-1}\exp(2\pi
i\,e_1)\dch(\sigma)$.
Ce sont les images de $\exp(e_0)$ et
$\exp(e_1)$ par
$\left<\dch(\sigma)\right)\tau(2\pi i)$.
La $\dbQ$-forme $(P_{0,0})_\sigma$ est
$\mu_N$-\'equivariante.
Les images des $\exp(e_\zeta)$ par
$\left<\dch(\sigma)\right>\tau(2\pi i)$ sont
donc elles aussi rationnelles.
La pro-alg\`ebre de Lie $\scrL^{\uphat}$ est
l'unique $\dbQ$-forme de
$\scrL_{\dbC}^{\uphat}$ pour laquelle les
$e_0$ et $e_\zeta$ sont rationnels, et $\prod$
est donc l'unique $\dbQ$-forme du groupe
$\prod_{\dbC}$ pour laquelle $\exp(e_0)$ et
les $\exp(e_\zeta)$ sont rationnels.
Transportant par
$\left<\dch(\sigma)\right>\tau(2\pi i)$, on
obtient l'assertion pour $(P_{0,0})_\omega$.

Montrons que la $\dbQ$-forme
$(P_{1,0})_\sigma$ de
$(P_{1,0})_\omega\otimes\dbC$ est l'image de
la $\dbQ$-forme $(P_{1,0})_\omega$ par
$\left<\dch(\sigma)\right>\tau(2\pi i)$.
Ces deux $\dbQ$-formes sont des
$(P_{0,0})_\sigma$-torseurs.
Il suffit donc de montrer qu'elles admettent
un point rationnel commun.
L'image $\dch(\sigma)$ de $1$ par
$\left<\dch(\sigma)\right>\tau(2\pi i)$ est en
effet dans $(P_{1,0})_\omega(\dbQ)$ car image
du droit chemin de $1_0$ \`a $(-1)_1$.
\enddemo

\subhead
5.19
\endsubhead
Notons encore $a_\sigma$ un \'el\'ement de
$G_\omega\left<\MT(\scrO[1/N])\right>$
image d'un \'el\'ement\break
 $a_\sigma\in G_\omega\left<\MT(k)\right>$ de 2.12, et
$\iota(a_\sigma)$ son image par 5.10.3 dans
$H_\omega$.
De m\^{e}me pour $a_\sigma^0$ (2.12.2).
L'\'el\'ement $\iota(a_\sigma)$ lui-aussi
transforme les $\dbQ$-formes $\omega$ de
$(P_{00})_\omega\otimes\dbC$ et
$(P_{1,0})_\omega\otimes\dbC$ en les
$\dbQ$-formes $\sigma$.
Il existe donc $v_\sigma\in H_\omega(\dbQ)$
tel que
$$
\left<\dch(\sigma)\right>\tau(2\pi
i)=\iota(a_\sigma)v_\sigma.\tag5.19.1
$$

Si comme en (2.12.2) on a choisit $a_\sigma$
de la forme $a_\sigma^0\tau(2\pi i)$ avec
$a_\sigma^0\in U_\omega(\dbC)$, tant
$\left<\dch(\sigma)\right>\tau(2\pi i)$ que
$\iota(a_\sigma)$ ont la m\^{e}me projection
$2\pi i$ dans le quotient $\dbG_m$ de
$H_\omega$, $v_\sigma$ est dans
$V_\omega(\dbQ)$ et on peut r\'ecrire (5.19.1)
sous les formes
$$
\gather
\left<\dch(\sigma)\right>=\iota(a_\sigma^0)\cdot
  (\tau(2\pi i)v_\sigma\tau(2\pi i)^{-1})\tag5.19.2\\
\tau(2\pi i)^{-1}\left<\dch(\sigma)\right>
   \tau(2\pi i)=
  (\tau(2\pi i)^{-1}\iota(a_\sigma^0)
  \tau(2\pi i))v_\sigma.\tag5.19.3
\endgather
$$

A ces formules correspondent deux fa\c{c}ons
de voir les relations entre les coefficients
(5.17.1) de $\dch(\sigma)$ impliqu\'ees par la
th\'eorie motivique.

La formule (5.19.2) implique que
$\dch(\sigma)\in\prod(\dbC)$ est contenu dans
la sous-vari\'et\'e alg\'ebrique d\'efinie sur
$\dbQ$ suivante.
Soit $\iota(U_\omega)\Subset\left(\prod,\circ\right)$
le sous-groupe unipotent image de $U_\omega$
par (5.10.3).
Soit $D_\sigma$ l'image de la droite de
coordonn\'ee $t$ par $t\mapsto
\tau(t)v_\sigma\tau(t)^{-1}$.
Cette application est d\'efinie en $t=0$ car
$\scrL^{\uphat}$ est \`a degr\'es $>0$.
On a $t=0\mapsto$ \'el\'ement neutre.
Par (5.19.2), $\dch(\sigma)$ est contenu dans
l'image par le produit $\circ$ de
$\iota(U_\omega)$ et $D_\sigma$.
Noter que $D_\sigma$ peut d\'ependre de
$\sigma$, mais $\iota(U_\omega)$ non.

La formule (5.19.3) sugg\`ere de consid\'erer
plut\^{o}t $\tau(2\pi i)^{-1}(\dch(\sigma))$,
de coefficients convergents donn\'es par la
formule 5.17 (i), divis\'ee par $(2\pi
i)^{\sum s_i}$.
Les autres coefficients sont d\'eduits de
ceux-l\`a par des formules rationnelles.
L'\'el\'ement $\tau(2\pi
i)^{-1}(\dch(\sigma))\in\prod(\dbC)$ est
contenue dans la sous-vari\'et\'e d\'efinie
sur $\dbQ$ $\iota(U_\omega)v_\sigma$, une
classe lat\'erale rationnelle pouvant
d\'ependre de $\sigma$ du sous-groupe
$\iota(U_\omega)$ de $\left(\prod, \circ \right)$

\subhead
5.20
\endsubhead
Supposons que $N=1$.
On a $k=\dbQ$; $\sigma$ est unique.
Choisissons $a_\sigma^0$ r\'eel (2.12).
Puisque $\dch(\sigma)$ est r\'eel, (5.12.2)
montre que $\tau(2\pi i)v_\sigma\tau(2\pi
i)^{-1}$ est r\'eel.
Le logarithme de $v_\sigma$ est donc purement
en degr\'e pair, $\tau(t)v_\sigma\tau(t)^{-1}$
ne d\'epend que de $t^2$ et ceci permet de
poser $d(t^2):=\tau(t)v_\sigma\tau(t)^{-1}$.

Regardons $\prod$ comme l'espace des
\'el\'ements groupaux de $\dbC\ll e_0,e_1\gg$.
L'alg\`ebre de Lie de $U_\omega$ \'etant en
degr\'e $\Ge 3$, les coefficients de
$\dch(\sigma)$ et de $v_\sigma$
co{\umi}ncident en degr\'ee $\Le 2$.
On a 
$$
\dch(\sigma)=1-\zeta(2)[e_0,e_1]+\text{degr\'e}\Ge 3
$$
et donc, puisque $\zeta(2)=-(2\pi i)^2/24$,
$$
v_\sigma=1+\tfrac{1}{24}[e_0,e_1]+\text{degr\'e}\Ge 3.
$$

Pour $g\in\iota(U_\omega(\dbC))$ et $d\in
D_\sigma$, $d=\tau(t)v_\sigma\tau(t)^{-1}$,
$t^2$ est donc $24$ fois le coefficient de
$e_0e_1$ dans $d$, \'egal \`a celui de $e_0e_1$
dans $g\circ d$.
Ceci d\'etermine $d$ \`a partir de $g\circ d$
et
$$
\iota(U_\omega)\times D_\sigma\to
\left(\prod, \circ \right)
$$
est un plongement ferm\'e: l'image est
l'ensemble des $x\in\prod\Subset\dbC\ll
e_0,e_1\gg$ tels que
$$
x\circ d\text{ (24.coefficient de $e_0e_1$
dans $x$)}^{-1}\in\iota(U_\omega).
\tag5.20.1
$$

Une variante d'une conjecture de Grothendieck
affirme que $a_\sigma$ est Zariski-dense (sur
$\dbQ$) dans $G_\omega$.
Si tel est le cas, toutes les relations
alg\'ebriques entre les coefficients de
$\dch_\sigma$ sont donn\'ees par 5.20.1.

On observera que les \'equations polynomiales
entre les coefficients de $\dch_\sigma$
d\'eduites de 5.20.1 ont les propri\'et\'es
suivantes.

\medskip\noindent
{\sans (a)}\enspace
Homog\'en\'eit\'e: toute \'equation est somme
d'\'equations isobares, pour le poids.
Ceci exprime que $\iota(U_\omega)\circ
D_\sigma\Subset V_\omega$ est stable sous
l'action int\'erieure de $\dbG_m\Subset
H_\omega$.
En effet, tant $\iota(U_\omega)$ que
$D_\sigma$ le sont.

\smallskip\noindent
{\sans (b)}\enspace
Si on ajoute l'\'equation ``$\pi^2=0$'',
plus pr\'ecis\'ement: ``nullit\'e du
coefficient de $e_0e_1$'', on obtient l'id\'eal 
des \'equations d\'efinissant le
sous-groupe $\iota(U_\omega)$ de
$\left(\prod,0\right)$.

\smallskip\noindent
{\sans (c)}\enspace
Si on ajoute plut\^{o}t l'\'equation
``$\pi^2=\alpha$'', on obtient une classe
lat\'erale de ce groupe.

\medskip
G. Racinet donne dans sa th\`ese un syst\`eme
d'\'equations polynomiales v\'erifi\'ees par
les coefficients de $\dch(\sigma)$ et prouve
que ce syst\`eme v\'erifie (a) et la forme
affaiblie suivante de (b) (c): si on ajoute
l'\'equation ``$\pi^2=0$'', on obtient un
sous-groupe de $\left(\prod,0\right)$, et si
on ajoute ``$\pi^2=\alpha$'', on obtient une
classe lat\'erale de ce sous-groupe.

\subhead
5.21
\endsubhead
La situation est similaire pour $N=2$,
$\zeta(2)$ continuant \`a jouer le m\^{e}me
r\^{o}le que pour $N=1$.
Il n'est plus vrai que $\Lie\,U_\omega$ est en
degr\'e $\Ge 3$, mais on est sauv\'e par la
compatibilit\'e entre les cas $N=1$ et $N=2$:
$$
\spreadlines{2\jot}
\CD
U_\omega\left<\MT(\dbZ[1/2])\right> @>>>
\left(\prod\text{ pour }N=2,o\right)\\
@VVV @VVV\\
U_\omega\left<\MT(\dbZ)\right> @>>>
  \left(\prod\text{ pour }N=1,o\right).
\endCD
$$

\subhead
5.22
\endsubhead
Pour $N\Ge 3$, $k$ est totalement complexe.
Le coefficient de $e_\zeta$ dans
$\dch(\sigma)$ est $-\log(1-\sigma(\zeta))$
(resp. $0$ pour $\zeta=1$), et le r\^{o}le
pr\'ec\'edemment jou\'e par le coefficient
$-\zeta(2)$ de $e_0e_1$ est jou\'e par la
diff\'erence des coefficients de $e_\zeta$
et de $e_{\zeta^{-1}}$, pour $\zeta\in\mu_N$
outre que $\pm1$.
Cette diff\'erence est un multiple rationnel
de $2\pi i$, d\'ependant de $\sigma$.
La courbe $D_\sigma$ est encore une droite,
param\'etr\'ee par la difference
$\delta_{\coeff}$ des coefficients de
$e_\zeta$ et $e_{\zeta^{-1}}$, 
$$
\iota(U_\omega)\times D_\sigma\to
\left(\prod, \circ\right)
$$
est un plongement ferm\'e, et l'image est
l'ensemble des $x\in\prod\Subset\dbC\ll
e_0,(e_\alpha)_{\alpha\in\mu_N}\gg$ tels que
$$
x\circ\text{(point de $D_\sigma$ donn\'e par
$\delta_{\coeff}(x))^{-1}\in\iota(U_\omega)$}.
\tag5.22.1
$$

La conjecture de Grothendieck affirme \`a
nouveau que toutes les relations polynomiales
\`a coefficients rationnels entre les
coefficients de $\dch(\sigma)$ proviennent de
(5.22.1).

Le syst\`eme d'\'equations polynomiales
d\'eduit de (5.22.1) a encore les
propri\'et\'es (a), (b), (c) de 5.20, avec
``$\pi^2$'' remplac\'e par ``$2\pi i\,$'' (plus
pr\'ecis\'ement, par la diff\'erence des
coefficients de $e_\zeta$ et
$e_{\zeta^{-1}}$), et Racinet a obtenu 
pour tout $N$ des r\'esultats semblables \`a ceux
mentionn\'es en 5.20.

\subhead
5.23
\endsubhead
Permettons \`a nouveau \`a $N$ d'\^{e}tre un
quelconque entier $\Ge 1$.
Puisque $\dch(\sigma)$ est un \'el\'ement
groupal de
$\dbC\ll e_0,(e_\alpha)_{\alpha\in\mu_N}\gg$, le
produit de deux coefficients de $\dch(\sigma)$
est combinaison lin\'eaire \`a coefficients
rationnels de coefficients, de sorte que
conna{\hati}tre les relations lin\'eaires entre les
coefficients suffit pour conna{\hati}tre  les
relations polynomiales entre eux.

\proclaim{Th\'eor\`eme 5.24}
Soit $d_n$ la dimension en degr\'e $n$ de
l'image de $\iota\colon\,\Lie\,U_\omega\to
\Lie\,V_\omega$ (5.10.3).
Soit $A$ l'alg\`ebre de polyn\^{o}mes
gradu\'ee engendr\'ee sur $\dbQ$ par $d_n$
g\'en\'erateurs degr\'e $n$, et par un
g\'en\'erateur additionnel $t_0$ en degr\'e
$1$ (resp. $2$ si $N=1$ ou $2$).

Il existe un homomorphisme $\varphi_\sigma$ de
$A$ dans $\dbC$, envoyant $t_0$ sur $2\pi i$
(resp. $\pi^2$ si $N=1$ ou $2$), tel que les
coefficients de mon\^{o}mes de degr\'e $d$
dans $\dch_\sigma$ soient contenus
dans\break
 $\varphi_\sigma$ (partie de degr\'e $d$ de $A$).
\endproclaim

\demo{Preuve}
On prend simplement pour $A$ l'alg\`ebre
affine de $\iota(U_\omega)\times D_\sigma$.
Elle a la structure dite d'alg\`ebre
gradu\'ee, la fonction ``coefficient d'un
mon\^{o}me de degr\'e $d$'' est de degr\'e
$d$, et $\dch(\sigma)$ est un point de
$\Spec(A)$ \`a valeurs dans $\dbC$, i.e. un
homomorphisme $\varphi_\sigma$ de $A$ dans
$\dbC$.
Le g\'en\'erateur $t_0$ est la coordonn\'ee de
$D_\sigma$, convenablement normalis\'ee.
\enddemo

\proclaim{Corollaire 5.25}
Soit $D_n$ la dimension du $\dbQ$-espace
vectoriel engendr\'e par les coefficients des
monomes de degr\'e $n$ dans $\dch(\sigma)$.
Le s\'erie g\'en\'eratrice $\sum D_nt^n$ est
terme \`a terme major\'ee par la suivante:

\roster
\item"{\sans (i)}"
Pour $N=1$: $1/(1-t^2-t^3)$.

\item"{\sans (ii)}"
$N=2$: $1/(1-t-t^2)$.

\item"{\sans (iii)}"
Pour $N\Ge 3$ ayant $\nu$ facteurs premiers:
$1/\left(1-\left(\frac{\varphi(N)}{2}+(\nu-1)\right)
t+(\nu-1)t^2\right)$, o\`u $\varphi$
est l'indicateur d'Euler
$\vert(\dbZ/N\dbZ)^*\vert$.
\endroster
\endproclaim

\demo{Preuve}
La s\'erie majorante est la s\'erie
g\'en\'eratrice des dimensions des composantes
homog\`enes de l'alg\`ebre gradu\'ee de
polyn\^{o}me construite comme l'alg\`ebre $A$
de 5.25, mais en partant de $\Lie(U_\omega)$
plut\^{o}t que d'un quotient.
En effet, cette s\'erie est le produit de
$\frac{1}{1-t}$ (resp $\frac{1}{1-t^2}$ si
$N=1$ ou $2$) par la s\'erie g\'en\'eratrice
pour l'alg\`ebre enveloppante de
$\Lie(U)_\omega^{\GR}$.
Cette derni\`ere est une alg\`ebre associative
libre.
S\'erie g\'en\'eratrice: $\sum
f(t)^n=1/(1-f(t))$, avec
$$
f(t)=\frac{\varphi(N)}{2}\,\,\frac{t}{1-t}+(\nu-1)t
$$
(resp. $f(t)=t^3/1-t^2$ si $N=1$, $t/1-t^2$ si
$N=2$).
Le calcul final est laiss\'e au lecteur.
\enddemo

\remark\nofrills{\bf Remarque 5.26.}\enspace
Pour $N=1$, nous retrouvons ainsi le
r\'esultat de Terasoma (2002).
Terasoma consid\`ere un groupe de cohomologie
relative, correspondant \`a un motif de Tate
mixte sur $\dbZ$ qui est sans doute le
suivant: dans l'alg\`ebre affine de ${}_1P_0$,
prendre la partie de poids $\Le 2n$, puis les
coinvariants de la monodromie locale en $0$
et $1$, agissant \`a droite (resp. gauche) sur
$P_{1,0}$.
Comme la n\^{o}tre, sa preuve repose sur 1.6.
\endremark

\subhead
5.27
\endsubhead
La conjonction de (a) toute relation
$\dbQ$-lin\'eaire entre les coefficients de
$\dch(\sigma)$ est somme de relations isobares
(b) la borne 5.25 est atteinte \'equivaut \`a
la conjonction de (c) la conjecture de
Grothendieck que $a_\sigma$ est $\dbQ$-Zariski
dense dans $U_\omega$, et (d) l'injectivit\'e
de $\iota\colon\, U_\omega\to V_\omega$.
La condition (d) est en g\'en\'eral fausse.
Elle est vraie pour $N=2$, $3$ ou $4$ et on
ignore si elle est vraie pour $N=1$.

Pour obtenir une estimation raisonable, il
faudrait conna{\hati}tre la dimension des
composantes gradu\'ees l'image de $\Lie\,U_\omega$
dans $\Lie\,V_\omega$.
On ne conna{\hati}t ces dimensions que pour
$N=2,3,4$.
La question plus int\'eressante serait de
d\'eterminer non seulement ces dimensions, mais
encore l'image elle-m\^{e}me de $\Lie\,U_\omega$.
On ne sait le faire dans aucun cas.

\subhead
5.28
\endsubhead
Il serait tr\`es int\'eressant de disposer
pour les motifs de Tate mixtes de
r\'ealisations cristallines qui, appliqu\'ees
aux espaces de chemins motiviques, fournissent
les r\'ealisations consid\'er\'ees en Deligne
(1989).

Dans le cas de motifs de Tate mixtes $M$ sur
$\dbZ[1/D]$, on veut, pour chaque nombre
premier $p$ premier \`a $D$, un automorphisme
dit de Frobenius
$$
F_p\colon\, M_\sigma\otimes\dbQ_p\to
M_\sigma\otimes\dbQ_p.
$$
Cet automorphisme est l'action d'un
\'el\'ement canonique de $G_\omega(\dbQ_p)$,
pour
$$
G_\omega=G_\omega\left<\MT(\dbZ[1/D])\right>.
$$
Sur $\dbQ(1)_\sigma\otimes\dbQ_p$, $F_p$ est
la multiplication par $1/p$.

Pour $N=1$ ou $2$, on a $k=\dbQ$ et, pour $p$
premier \`a $N$, on disposerait donc de
$F_p^{-1}$ de la forme
$$
F_p^{-1}=\varphi_p\tau(p)\in G_\omega(\dbQ_p)
$$
et de l'image $\iota(\varphi_p)$ de
$\varphi_p$ dans $V_\omega(\dbQ_p)$.
Plongeons comme d'habitude $V_\omega(\dbQ_p)$
dans\break
$\dbQ_p\ll e_0,e_1\gg$, pour $N=1$, ou
$\dbQ_p\ll e_0,e_1,e_{-1}\gg$, pour $N=2$.
Les coefficients de $\iota(\varphi_p)$
semblent \^{e}tre des analogues $p$-adiques des
valeurs multiz\^{e}ta (cf Deligne (1989) 3.4).
D'apr\`es 5.20, on s'attend \`a ce qu'ils
v\'erifient toutes les identit\'es satisfaites
par les coefficients de $\dch(\sigma)$, ainsi
que ``$\pi^2=0$'' (nullit\'e du coefficient de
$e_0e_1$).
Le formalisme cristallin esp\'er\'e prouverait
en tout cas les identit\'es d\'eduites de
(5.20.1), et de ``$\pi^2=0$''.

Il serait int\'eressant aussi de disposer pour
ces coefficients d'expressions $p$-adiques qui
rendent claires qu'ils v\'erifient des
identit\'es du type
$$
\align
&\coeff(e_0^{n-1}e_1)\coeff(e_0^{m-1}e_1)=\\
&\coeff(e_0^{n-1}e_1e_0^{m-1}e_1)+\coeff
  (e_0^{m-1}e_1e_0^{n-1})\\
&\qquad\qquad\qquad +\coeff(e_0^{n+m-1})e_1),
\endalign
$$
qui pour $\dch(\sigma)$ expriment que
$$
\sum \frac{1}{k^n}.\sum \frac{1}{\ell^m}=
\sum\limits_{k>\ell}\frac{1}{k^n}
\frac{1}{\ell^m}+\sum\limits_{\ell>k}
\frac{1}{k^n}\frac{1}{\ell^m}+\sum
\frac{1}{k^{n+m}}.
$$

Pour $N\Ge 3$, le plus commode est de se
rappeller que, m\^{e}me si les points de
$\mu_N$ ne sont d\'efinis que sur $k$,
$X=A^*-\mu_N$ est d\'efini sur $\dbQ$, de
sorte que $\pi_1(X,0)$ et $P_{1,0}(X)$ sont de
fa\c{c}on naturelle des
$\MAT(k/\dbQ)$-sch\'emas.
La r\'ealisation de de Rham de
$\Lie\,\pi_1(X,0)$ est une $\dbQ$-forme
$\scrL_{\DR}^{\uphat}$ de
$\scrL^{\uphat}\widehat{\otimes} k$: les
invariants de $\Gal(k/\dbQ)$ agissant sur $k$
et permutant les $e_\zeta$ ($\zeta\in\mu_N$).
Pour $p\nmid N$, la th\'eorie cristalline
fournit un automorphisme de
$\scrL_{\DR}^{\uphat}$ et du
$\exp(\scrL_{\DR}^{\uphat})$-torseur
$P_{1,0}(X)_{\DR}$.

\newpage


\dspace
\subhead
6. Profondeur 
\endsubhead

\subhead
6.1
\endsubhead
Avec les notations de 5.1, pour tout point
base $x$, l'inclusion de $X$ dans
$A^*$ induit un morphisme du
groupe fondamental motivique de $X$ dans
celui, $\dbQ(1)$, de $A^*$.
On d\'efinit la filtration par la profondeur
$D$ de $\pi_1(X,x)$ par
$$
\align
&D^0\pi_1(X,x):=\pi_1(X,x)\tag6.1.1\\
&D^1\pi_1(X,x):=\Ker(\pi_1(X,x)
\to \pi_1(A^*,x)
\endalign
$$
et, pour $n\Ge 1$, 
$$
\align
D^n\pi_1(X,x):=&\text{le $n^{\ieme}$ groupe de la}\\
&\text{s\'erie centrale descendante de
$D^1\pi_1(X,x)$}.
\endalign
$$
On filtre de m\^{e}me les alg\`ebres de Lie:
$$
D^n\Lie\,\pi_1(X,x)=\Lie\,D^n\pi_1(X,x).
\tag6.1.2
$$

Si on pousse le $\pi_1(X,x)$-torseur $P_{y,x}$
par $\pi_1\to \pi_1/D^n$, les
$P_{y,x}/D^n$ obtenus forment un
groupo{\umi}de quotient du groupo{\umi}de
des $P_{y,x}$.

Appliquons A.10, A.11.
La filtration par la profondeur induit une
filtration de l'alg\`ebre enveloppante
compl\'et\'ee $\scrU\uphat\Lie\,\pi_1(X,x)$
ainsi que de l'alg\`ebre affine
$\scrO(\pi_1(X,x))$, qui en est duale.
Elle induit aussi  une filtration sur les
alg\`ebres affines des torseurs $P_{y,x}$.
L'alg\`ebre affine du quotient
$\pi_1(X,x)/D^{n+1}$ de $\pi_1(X,x)$ est la
sous-alg\`ebre de l'alg\`ebre affine de
$\pi_1(X,x)$ engendr\'ee par $D^{-n}
\scrO(\pi_1(X,x))$.
De m\^{e}me pour les $P_{y,x}$.

\subhead
6.2
\endsubhead
En r\'ealisation $\omega$, l'alg\`ebre de Lie
de $\pi_1(X,x)$ est $\scrL^{\uphat}$ (5.8), le
compl\'et\'e de l'alg\`ebre de Lie libre
$\scrL$ engendr\'ee par $e_0$ et les $e_\zeta$
($\zeta\in\mu_N$).
La filtration de profondeur 
est d\'eduite 
la graduation par le degr\'e en les $e_\zeta$,
appel\'e le {\it $D$-degr\'e}.
Cette graduation 
est compatible \`a la
graduation provenant de celle de $\omega$ pour
laquelle $e_0$ et les $e_\zeta$ sont de
degr\'e un et est $\mu_N$-\'equivariante.
Elle se propage aux alg\`ebres
enveloppantes, aux alg\`ebres affines et aux
alg\`ebres affines des torseurs triviaux
$(P_{y,x})_\omega$.
L'action correspondante $\rho_D$ de $\dbG_m$
pour laquelle $\rho_D(\lambda)$ fixe $e_0$ et
transforme $e_\zeta$ en $\lambda e_\zeta$ sera
appell\'ee action ``$D$-degr\'e''.

On prendra garde que les filtrations $D$ sont
motiviques et stables par l'automorphisme
$u\mapsto u^{-1}$ de $X$, mais que le
$D$-degr\'e est d\'efini seulement en
r\'ealisation $\omega$ et n'est pas stable par
$u\mapsto u^{-1}$.
L'automorphisme de $\scrL\uphat$ induit par
$u\mapsto u^{-1}$  \'echange $e_0$, qui est
homog\`ene de $D$-degr\'e $0$, et
$e_\infty=-e_0-\sum e_\zeta$,  qui ne l'est
pas.

\subhead
6.3
\endsubhead
Le $\MT(\scrO[1/N])$-sch\'ema en groupe $V$
de 5.15 agit sur $P_{0,0}$ et $P_{1,0}$.
On d\'efinit la filtration $D$ de $\Lie\,V$
par la condition que $D^p\Lie\,V$ soit maximal
tel que l'action de $\Lie\,V$ sur
$\scrO(P_{0,0})$ et $\scrO(P_{0,1})$ envoie
$D^p\Lie\,V\otimes D^q\scrO(P_{0,0})$ dans
$D^{p+q}\scrO(P_{0,0})$ et $D^p\Lie\,V\otimes
D^q\scrO(P_{1,0})$ dans
$D^{p+q}\scrO(P_{1,0})$.

Passons \`a la r\'ealisation $\omega$, o\`u
$V_\omega$ s'identifie \`a
$\left(\prod,\circ\right)$ agissant sur les copies
$\prod_0$ et $\prod_{1,0}$ de $\prod$ par
$g\mapsto\left<g\right>_0$ et
$\left<g\right>_{1,0}$.

\proclaim{Proposition 6.4}
L'action ``$D$-degr\'e'' $\rho_D$ de $\dbG_m$
sur $\prod$ respecte les actions
$\left<g\right>_0$ et $\left<g\right>_{10}$ de
$\left(\prod,\circ\right)$ sur $\prod$.
Elle respecte donc la loi de groupe $\circ$ sur
$\prod$ et le crochet $\{\,\,,\,\,\}$ de
$\scrL\uphat$.
\endproclaim

L'action $\left<g\right>_0$ de $g\in\prod$ sur
le groupe $\prod$ est caract\'eris\'ee par les
propri\'et\'es de commuter \`a l'action de
$\mu_N$ sur $\prod$, de fixer $\exp(e_0)$ et
d'envoyer $\exp(e_1)$ sur $g^{-1}\exp(e_1)g$.
Son conjugu\'e par $\rho_D(\lambda)$ commute
encore \`a $\mu_N$, fixe $\exp(e_0)$ et envoie
$\exp(\lambda e_1)$ sur\break
$\rho_D(\lambda)(g)^{-1}\exp(\lambda e_1)
\rho_D(\lambda)(g)$.
Cette derni\`ere propri\'et\'e \'equivaut \`a
envoyer $\exp(e_1)$ sur
$\rho_D(\lambda)(g)^{-1}\exp(e_1)
\rho_D(\lambda)(g)$: $\rho_D(\lambda)$
conjugue $\left<g\right>_0$ en
$\left<\rho_D(\lambda)(g)\right>_0$.
Conjugant $\left<g\right>_{10}\colon\,x\mapsto
g\left<g\right>_0(x)$ par $\rho_D(\lambda)$,
on en d\'eduit que $\rho_D(\lambda)$
transforme $\left<g\right>_{10}$ en
$\left<\rho_D(\lambda)(g)\right>_{10}$.

\subhead
6.5
\endsubhead
Puisque l'action de $\Lie\,V_\omega=(\scrL\uphat,
\{\,\,,\,\,\})$ sur $\scrO(P_{1,0})_\omega$
est fid\`ele, il r\'esulte de 6.4 que la
filtration $D$ de $\Lie\,V_\omega$ co{\umi}ncide
avec la filtration $D$ de $\scrL\uphat$
d\'eduite de la graduation par le $D$-degr\'e
de $\scrL$.
En particulier, $D^0\Lie\,V_\omega=\Lie\,V_\omega$.

Transportons \`a $V_\omega$ l'action
``$D$-degr\'e'' $\rho_D$ de $\dbG_m$ sur
$\left(\prod,\circ\right)$.
Cette action gradue $\Lie(V_\omega)$, son
alg\`ebre enveloppante universelle
compl\'et\'ee et l'alg\`ebre affine
$\scrO(V_\omega)$.
Les filtrations correspondantes
co{\umi}ncident avec celles d\'efinie par la
filtration $D$ de $\Lie\,V_\omega$ par les
$D^p\Lie\,V_\omega$.
Il en r\'esulte que l'isomorphisme de 
sch\'emas de $V_\omega$ 
avec $\prod$ induit un isomorphisme de
$V_\omega/D^pV_\omega$ avec $\prod/D^p\prod$,
les alg\`ebres affines des quotients \'etant
dans les deux cas les sous-alg\`ebres
engendr\'ees par $D^{-p}$ des alg\`ebres
affines (A.10).
Plus pr\'ecis\'ement, l'action de $V_\omega$
sur $\prod_{10}$ fait de $\prod_{10}$ un
espace principal homog\`ene sous $V_\omega$,
trivialis\'e par $1_{10}$ fixe sous $\rho_D$,
et induit une action principale homog\`ene de
$V_\omega/D^pV_\omega$ sur $\prod_{10}/D^p$.
Ind\'ependamment du foncteur fibre, il reste
vrai que l'action de $V$ sur $P_{0,0}$ et
$P_{1,0}$ induit une action de $V/D^pV$ sur
$P_{0,0}/D^p$ et $P_{1,0}/D^p$, et que
$P_{1,0}/D^p$ est un espace principal
homog\`ene sous $V/D^pV$. 

\proclaim{Proposition 6.6}
Le morphisme 5.15 de $U$ dans $V$ se factorise
par $D^1V$.
\endproclaim

Le quotient $D^1V/D^2$ \'etant ab\'elien, ce
morphisme induit
$$
U^{\ab}\to D^1V/D^2V.\tag6.6.1
$$

\demo{Preuve}
D'apr\`es 6.5, il suffit de v\'erifier que $U$
agit trivialement sur $P_{1,0}/D^1$, le
$\dbQ(1)$-torseur des chemins, dans
$A^*$, de $1_0$ \`a $1$.
En effet, $P_{\mu,\lambda_0}$
($A^*$) est canoniquement
isomorphe \`a
$P_{\mu,\lambda}(A^*)$, et,
faisant $\lambda=\mu$, on voit que
$P_{1,0}/D^1$ est le $\dbQ(1)$-torseur trivial,
sur lequel $U$ agit trivialement.
\enddemo

\subhead
6.7
\endsubhead
Nous nous proposons de montrer que (6.6.1) est
injectif et de calculer son image.
Noter que l'action de $U$ sur $U^{\ab}$, ainsi
que l'action de $U$ sur $D^1V/D^2V$
(d\'eduite de l'action int\'erieure de $V$ sur
lui-m\^{e}me), sont triviales.
Les deux membres de (6.6.1) sont donc des
produits de sch\'emas en groupe vectoriels
$\dbQ(n)$, et le passage \`a la r\'ealisation
$\omega$ ne perd aucune information.

En r\'ealisation $\omega$, le noyau $D^1\scrL$
de $\Lib(e_0,(e_\zeta))\to\Lib(e_0)$:
$e_0\mapsto e_0$, $e_\zeta\mapsto 0$, est,
pour le crochet [\quad], l'alg\`ebre de $\Lie$
libre engendr\'ee par les $(\ad\,e_0)^n(e_\zeta)$
($n\Ge 0$, $\zeta\in\mu_N)$, de $D$-degr\'e
$1$.
Le quotient $D^1\scrL/D^2\scrL$ a donc pour
base les images des $(\ad\,e_0)^n(e_\zeta)$,
$\ad\,e_0$ \'etant pris au sens du crochet
$[\,\,,\,\,]$.

Notons $\scrE^n$ la composante de poids $-2n$
de $D^1\scrL/D^2\scrL$, de base les images des\break
$(\ad\,e_0)^{n-1}(e_\zeta)$.
Le quotient $\Gr_D^1(\scrL\uphat)$ est le
produit des $\scrE^n$.
Quand $n$ sera fix\'e, nous \'ecrisons
simplement $\scrE$ pour $\scrE^n$, et
$E_\zeta$ pour la base des
$(\ad\,e_0)^{n-1}(e_\zeta)$.

Fixons un plongement complexe $\sigma$ de $k$.
L'\'el\'ement $\dch(\sigma)$ de
$V_\omega(\dbC)\Subset\dbC\ll e_0,(e_\zeta)\gg$
a un coefficient de $e_0$ nul (5.17).
Il est donc dans $D^1V_\omega(\dbC)$.
Calculons sa projection dans
$D^1V_\omega/D^2V_\omega$, d'alg\`ebre de Lie
$D^1\scrL\uphat/D^2\scrL\uphat$.
Pour
$$
x=\sum\lambda_{n,\zeta}(\ad e_0)^{n-1}(e_\zeta)
+\text{\rm termes de $D$-degr\'e $\Ge 2$},
$$
l'exponentielle de $x$ dans $V_\omega=\prod$
est donn\'ee par 5.14 (ii).
On a
$$
\exp(x)=1+x+\text{\rm termes de $D$-degr\'e $\Ge
2$},
$$
et $\lambda_{n,\zeta}$ est donc le coefficient
de $e_0^{n-1}e_\zeta$ dans $x$.

D'apr\`es 5.17, la projection de
$\dch(\sigma)$ dans $D^1V_\omega/D^2V_\omega$
est donc l'exponentielle de l'\'el\'ement
$$
\dch^1(\sigma)=\sum\limits_{n,\zeta}
\biggl(\sum\limits_{m}
\frac{\sigma(\zeta^{-m})}{m^n}\biggr)
(\ad\,e_0)^{n-1}(e_\zeta)
$$
de son alg\`ebre de Lie $\prod\scrE^n$.
Pour $n=1$, $\zeta=1$, la somme sur $m$ diverge
et est \`a remplacer par $0$.

Identifions $U^{\ab}$ et $D^1V/D^2V$ \`a
leurs alg\`ebre de Lie par l'application
exponentielle et projetons sur $U^{\ab}$ et
$D^1V/D^2V$.
On obtient qu'il existe $a_\sigma$ dans 
$\gru_\omega^{\ab}
\widehat{\otimes}\dbC$ v\'erifiant
$$
\align
&\text{Image de $a_\sigma$ dans
$\left(\Lie\,D^1V_\omega/D^2V_\omega=\prod \scrE^n\right)\widehat{\otimes}\dbC $}
\tag6.7.1\\
&=\text{$\dch^1(\sigma)+$ \'el\'ement de
$\prod(2\pi i)^n\scrE^n$}.
\endalign
$$

Notons  $\scrI$ l'image de Lie $U_\omega$ dans Lie
$V_\omega$.
La filtration $D$ de Lie $V$ induit une
filtration $D$ de $\scrI$.
On prendra garde que l'action ``$D$-degr\'e''
de $\dbG_m$ sur $\Lie\,V_\omega$ n'a aucune
raison de respecter $\scrI$ (et en fait ne
respecte pas $\scrI$).

\proclaim{Th\'eor\`eme 6.8}\enspace
{\rm (i)} Le morphisme (6.6.1) est injectif; il
induit donc un isomorphisme de
$\gru_\omega^{\ab}$ avec $\scrI/D^2\scrI$.

\noindent
{\rm (ii)} En r\'ealisation $\omega$ et en poids
$-2n$, l'image dans $\scrE:=\scrE^n$ de
$\omega_n(\gru^{\ab})$ est le sous-espace
$\scrE_{\distr}^+$ des $\sum x_\zeta E_\zeta$
v\'erifiant la relation de parit\'e
$$
x_{\zeta^{-1}}=(-1)^{n-1}x_\zeta\tag6.8.1
$$
et la relation de distribution de poids $n-1$
(resp. de distribution trou\'ee si $n=1$): si
$N=Md$ avec $d\Ge 1$ et que $\eta^M=1$, on a
$$
x_\eta=\frac{1}{d^{n-1}}\sum\limits_{\zeta^d=\eta}
x_\zeta\tag6.8.2
$$
(resp. $x_1=0$ et (6.8.2) seulement pour
$\eta\not=1$ si $n=1$).
\endproclaim

\remark\nofrills{\bf Remarque 6.9 {\rm (i)}.}\enspace
La relation de parit\'e (6.8.1) est la
relation (6.8.2) pour $d=-1$.

(ii)\enspace
Le $\dbQ$-espace vectoriel $\scrE_{\distr}^+$
est muni de formes lin\'eaires $x\mapsto
x_\zeta$ v\'erifiant les relations de parit\'e
et de distribution de 6.8 (ii), et est
universel pour cette propri\'et\'e.
Son dual $\scrD_{\distr}^+$ est muni d'une
application $\mu_N\to\scrD_{\distr}^+$ v\'erifiant
les m\^{e}mes identit\'es, et est universel
pour cette propri\'et\'e: c'est l'espace dans
lequel prend ses valeurs la distribution de
poids $n-1$ (resp. trou\'ee si $n=1$) de
parit\'e $(-1)^{n-1}$ universelle.
\endremark

\demo{6.10 Preuve de 6.8}
Pour $\sigma$ un plongement complexe de $k$,
posons
$$
d_\sigma:=\sum\limits_{\zeta}\sum\limits_{m}
\frac{\sigma(\zeta^{-m})}{m^n}E_\zeta.
\tag6.10.1
$$
Pour $n=1$, le coefficient, divergent, de
$E_1$ est \`a remplacer par z\'ero.

D'apr\`es 6.7.1, il existe $a_\sigma$ dans
l'image de $\omega_n(\gru^{\ab})$ tel que
$$
d_\sigma\equiv a_\sigma\pmod{(2\pi i)^n\scrE}.
\tag6.10.2
$$
\enddemo

\medskip\noindent
{\sans Cas $n=1$.}\enspace
Fixons $\sigma$ et identifions $k$ \`a son
image dans $\dbC$ par $\sigma$.
En particulier, \'ecrivons $\zeta$ pour
$\sigma\zeta$.
Avec cette notation, (6.10.1) se r\'ecrit
$$
(d_\sigma)_\zeta=-\log(1-\zeta).
$$

D'apr\`es (2.3.1), $\omega_n(\gru^{\ab})$ est
le dual de $\scrO_N[1/N]^*\otimes\dbQ$.
Si $\scrO_N^+$ est la partie r\'eelle de
$\scrO_N$, c'est encore le dual de
$\scrO_N^+[1/N]^*\otimes\dbQ$.
Soit $\scrF$ son image dans $\scrE$.

Puisque $a_\sigma$ est dans $\scrF_\dbC$, sa
partie r\'eelle est dans $\scrF_\dbR$.
Prenant la partie r\'eelle de (6.10.2), on
obtient que
$$
-\sum \log\vert1-\zeta\vert~E_\zeta
$$
est dans $\scrF_\dbR$.

D'apr\`es un th\'eor\`eme de Bass (1966) (cf.
Washington (1997) 8.9 et 12.18),
$\scrO_N[1/N]^*\otimes\dbQ$, muni de
$\zeta\mapsto 1-\zeta$ (resp. $1$ si
$\zeta=1$) est la distribution trou\'ee paire
universelle sur $\mu_N$.
Puisque $\scrO_M^+[1/N]^*
\otimes\dbQ\simover\scrO_N[1/N]^*\otimes\dbQ$, 
l'application 
$x\mapsto \log\,\vert x\vert$ est injective sur
$\scrO_N[1/N]^*/\torsion$ et il n'y a donc
d'autres relations rationnelles entre les
$\log\,\vert1-\zeta\vert$ que celles
d\'eduites de (6.8.1) et (6.8.2).
L'image  $\scrF$ contient donc
$\scrE_{\distr}^+$.
Puisque $\scrE_{\distr}^+$ a le m\^{e}me rang
que $\scrO_N^*[1/N]^*\otimes\dbQ$, (i) et (ii)
en r\'esultent.

\medskip\noindent
{\sans Cas $n>1$.}\enspace
Soit $\scrF$ l'image de
$\omega_n*(\gru^{\ab})$ dans $\scrE$.
Prenant, selon la parit\'e de $n$, la partie
r\'eelle ou la partie imaginaire de (6.10.2),
on obtient que pour tout $\sigma$
$$
d_\sigma+(-1)^{n-1}\dbar_\sigma
$$
est dans $\scrF_\dbC$.
Il r\'esulte de (6.10.1) et des identit\'es
$$
\align
&\sum\limits_{\zeta^d=\eta}\zeta^m
=\cases
0 &\text{si $d\nmid m$}\\
d\eta^{m/d} &\text{si $d\mid m$}
\endcases\\
\vspace{5pt}
&\zeta^{-1}=\zetabar
\endalign
$$
que les $d_\sigma+(-1)^{n-1}\dbar_\sigma$ sont
dans $(\scrE_{\distr}^+)_\dbC$.

\proclaim{Lemme 6.11}
Le sous-espace $W$ de $\scrE_\dbC$ engendr\'e
par les $d_\sigma\pm (-1)^{n-1}\dbar_\sigma$
est de dimension au moins \'egale au nombre de
caract\`eres de $(\dbZ/N\dbZ)^*$ de parit\'e
$(-1)^{n-1}$.
\endproclaim

D\'eduisons 6.8 de 6.11.
La borne inf\'erieure 6.11 pour la dimension
de $W$ co{\umi}ncide tant avec la dimension de
$\omega(\gru^{\ab})_n$ (2.3.1) qu'avec la
dimension de $\scrE_{\distr}^+$ (Kubert (1979)
App.).
Puisque $W$ est contenu dans l'image
complexifi\'ee $\scrF_\dbC$, $W$ est \'egal
\`a cette image, (6.1.1) est injectif.
De m\^{e}me, puisque $S$ est contenu dans
$(\scrE_{\distr}^+)_{\dbC}$, on a \'egalit\'e,
prouvant 6.8.

\demo{Preuve de $6.11$}
On peut supposer que $k=\dbQ(\exp(2\pi
i/N))\Subset \dbC$.
On identifie l'ensemble des plongements
complexes de $k$ \`a $(\dbZ/N)^*$ par
$\zeta\mapsto \zeta^a$, et $\mu_N$ \`a
$\dbZ/N$ par $b\mapsto\exp(2\pi i\,b/N)$.
Avec ces notations, $W$ est l'image de
l'application
$$
F\colon\, \dbC^{(\dbZ/N)^*}\to
\dbC^{\dbZ/N}
$$
de matrice
$$
\sum\limits_{m}\,
\frac{\exp(-2\pi i\,abm)}{m^n}
+(-1)^{n-1}
\frac{\exp(-2\pi i\,abm)}{m^n}\,\,.
$$
Il nous faut minorer le rang de $F$.
Il sera plus commode de consid\'erer
l'application transpos\'ee ${}^tF$.
Les applications $F$ et ${}^tF$ sont
\'equivariantes, pour les actions $e_x\mapsto
e_{ax}^{-1}$ et $e_x\mapsto e_{ax}$ de
$(\dbZ/N)^*$ sur les deux membres, et il
suffit de v\'erifier que pour tout caract\`ere
$\chi$ de $(\dbZ/N)^*$ de parit\'e
$(-1)^{n-1}$, la $\chi$-partie de ${}^tF$
est non nulle.
Soient $M$ le diviseur de $N$ tel que $\chi$
provienne d'un caract\`ere primitif $\chi_0$
de $(\dbZ/M)^*$, et prolongeant $\chi_0$ par $0$
sur $\dbZ/M$.
V\'erifions que l'image par ${}^tF$ du vecteur
$(\chi_0(b))_{b\in\dbZ/N}$ est non nulle.
Prenant la composante $a=1$, on se ram\`ene
\`a v\'erifier que
$$
\sum\limits_{b,m}\,
\frac{\exp(2\pi i\,bm).\chi_0(b)}{m^n}
\not=0
$$
C'est en effet $N/M$ fois le produit de la
somme de Gau{\ss} $\sum\exp(2\pi
i\,b)\chi_0(b)$ (somme sur $(\dbZ/M)^*$)
par $L(\chi_0,n)$.
\enddemo

\remark{\bf Remarque 6.12}
L'argument utilis\'e pour $n>1$, utilisant
tous les plongements complexes et (6.10.2)
modulo $(2\pi i)^n\scrE\otimes\dbR$
plut\^{o}t qu'un seul plongement et (6.10.2)
modulo $2\pi i\scrE$, s'applique aussi pour
$n=1$ si $N$ n'a qu'un seul facteur premier.
Dans ce cas en effet, l'application
r\'egulateur de $\scrO_N^*[1/N]$ dans $\dbR^m$,
$m$ le nombre de places \`a l'infini, a un
noyau fini et  une image discr\`ete.
\endremark

\remark{\bf Remarque 6.13}
Combinant les r\'esultats de cet article avec
ceux de Goncharov (2001), on obtient
imm\'ediatement l'analogue motivique de
r\'esultats de loc. cit. sur les alg\`ebres de
Lie de Galois: remplacer syst\'ematiquement
l'alg\`ebre de Lie de Galois $g_N^\ell$
par $\im(\Lie\,U_\omega\to\Lie\,V_\omega)$
(voir 10.5.3) qui est son analogue motivique,
et travailler sur $\dbQ$ plut\^{o}t que
$\dbQ_\ell$.
\endremark

\newpage


\dspace
\subhead
Appendice. Rappels sur les groupes
alg\'ebriques unipotents
\endsubhead

Dans tout cet appendice, nous supposons
\^{e}tre en caract\'eristique z\'ero.

\subhead
A.1
\endsubhead
En caract\'eristique z\'ero, le th\'eor\`eme de
Poincar\'e-Birkhoff-Witt peut \^{e}tre
\'enonc\'e comme suit. 
Pour $\scrL$ une
alg\`ebre de Lie sur $K$, soit
$\varphi_n\colon\,\Sym^n\scrL\to\scrU\scrL$ le
produit sym\'etris\'e, de la puissance
symm\'etrique $n^{\ieme}$ de $\scrL$ dans
l'alg\`ebre enveloppante de $\scrL$.
Pour $x$ dans $\scrL$, il envoie
$x^n\in\Sym^n\scrL$ sur $x^n\in\scrU\scrL$.
Le morphisme d'espaces vectoriels
$$
\varphi_*\colon\,\Sym^*\scrL\to\scrU\,\scrL
\tag{A.1.1}
$$
est un isomorphisme de cog\`ebres.
Il transforme la filtration croissante de
$\Sym^*\scrL$ par les $\oplusop_{n\Le
p}\Sym^n\scrL$ en la filtration de
$\scrU\scrL$ par les images des
$\oplusop_{n\Le p}\otimes^n\scrL$.

Sous cette forme, le th\'eor\`eme est valable
dans toute cat\'egorie tensorielle (de
caract\'eristique $0$) pseudo-ab\'elienne
(pour donner un sens \`a $\Sym^n$) qui admet
des sommes d\'enombrables. 
Le cas qui nous int\'eresse est celui des
espaces vectoriels filtr\'es.
Soit $\scrL$ une alg\`ebre de Lie munie d'une
filtration d\'ecroissante $F$.
La filtration correspondante de $\scrU\,\scrL$
est caract\'eris\'ee par la propri\'et\'e qu'une
filtration sur un $\scrL$-module $M$ est
compatible \`a la filtration de $\scrL$ si et
seulement si elle l'est \`a celle de
$\scrU\,\scrL$.
La version cat\'egorique de
Poincar\'e-Birkhoff-Witt implique
 la compatibilit\'e de (A.1.1) aux filtrations.

Une filtration $F$ sur un
$\scrL$-module $M$ est compatible la
filtration centrale descendante $\grZ$ de
$\scrL$, si et
seulement si les $F^pM$ sont des sous-modules
et que $\scrL$ agit trivialement sur les
$\Gr_F^pM$.
La filtration de $\scrU\,\scrL$ induite par le
filtration $\grZ$ de $\scrL$ est donc la
filtration par les puissances de l'id\'eal
d'augmentation $I$.
Par compatibilit\'e de (A.1.1) aux
filtrations, on a

\proclaim{Corollaire A.2}
Si l'alg\`ebre de Lie $\scrL$, de s\'erie
centrale descendante  $\grZ$, est nilpotente
de classe $c$ $(\grZ^{c+1}=0)$, alors, si on
identifie $\Sym^*\scrL$ et $\scrU\,\scrL$ par
(A.1.1), on a
$$
I^n\Supset\oplusop_{p\Ge n}\Sym^p\scrL\Supset
I^{nc}\tag{A.2.1}
$$
\endproclaim

Notons $\scrU^{\uphat}\scrL$ le compl\'et\'e
$I$-adique de $\scrU\,\scrL$. 

\proclaim{Corollaire A.3}
Si l'alg\`ebre de Lie $\scrL$ est nilpotente,
l'isomorphisme (A.1.1) induit par compl\'etion
un isomorphisme d'espace vectoriels
$$
\prod\limits_{n\Ge 0}\Sym^n\scrL\to
\scrU^{\uphat}\scrL.
\tag{A.3.1}
$$
\endproclaim

L'isomorphisme (A.3.1) est compatible aux
topologies des deux membres.
A gauche, la topologie est la topologie
produit (d'espaces discrets).
A droite, celle de limite projective des
(espaces discrets) $\scrU\,\scrL/I^n$.

\subhead
A.4
\endsubhead
Le foncteur ``alg\`ebre de Lie'' est une
\'equivalence de la cat\'egorie des groupes
alg\'ebriques unipotents avec celle des
alg\`ebres de Lie nilpotentes de dimension
finie.
De plus, le foncteur $\rho\mapsto d\rho$ est
une \'equivalence de la cat\'egorie des
repr\'esentations lin\'eaires (de dimension
finie) de $V$, suppos\'e unipotent, avec celle
des repr\'esentations nilpotentes de son
alg\`ebre de Lie $\bfv$.

On dispose d'une application exponentielle
$$
\exp\colon\,\bfv\to V\tag{A.4.1}
$$
caract\'eris\'ee par la propri\'et\'e que pour
toute repr\'esentation lin\'eaire  $\rho$ de
$V$, on ait
$$
\rho(\exp(x))=\exp(d\rho(x)).\tag{A.4.2}
$$
C'est un isomorphisme de $\bfv$, muni de la
loi de groupe de Campbell-Hausdorff, avec $V$.

\subhead
A.5
\endsubhead
Soit $V$ un groupe alg\'ebrique unipotent
d'alg\`ebre de Lie $\bfv$.
L'alg\`ebre enveloppante $\scrU\bfv$ est
l'alg\`ebre des op\'erateurs diff\'erentiels
sur $V$ invariants par translations \`a gauche.
Comme cog\`ebre, c'est le dual topologique du
compl\'et\'e de l'anneau local de $V$ en $1$.
L'application exponentielle (A.4.1) induit un
isomorphisme de la cog\`ebre $\scrU\bfv$ avec
la cog\`ebre $\Sym^*\bfv$ des op\'erateurs
diff\'erentiels en $0$ sur l'espace affine
$\bfv$.

Par r\'eduction au cas o\`u $\bfv$ est de
dimension un, on v\'erifie que cet
isomorphisme envoie $x^n\in\Sym^n\bfv$ sur
$x^n\in\scrU\bfv$.
C'est donc l'isomorphisme (A.1.1).
Pour l'espace affine $\bfv$, l'alg\`ebre
affine de $\bfv$ est le dual topologique de la
cog\`ebre compl\'et\'ee
$\prod\limits_{n>0}\Sym^n\bfv$.
Appliquant A3, on obtient:

\proclaim{Proposition A.6}
L'alg\`ebre affine de $V$ est le dual
topologique de la cog\`ebre $\scrU\bfv$,
compl\'et\'ee pour la topologie $I$-adique.
\endproclaim

Dualement, le dual (muni de la topologie
faible) de l'alg\`ebre affine de $V$ est
$\scrU^{\uphat} \bfv$.
Les points de $V$ \`a valeurs dans une
alg\`ebre $A$ (homomorphismes
de l'alg\`ebre affine dans $A$) sont les {\it
\'el\'ements groupaux} $g$ de
$A\widehat{\otimes}\scrU^{\uphat}\bfv$: les
\'el\'ements d'augmentation $1$ tels que
$\Delta g=g\otimes g$.
Compatibilit\'es: l'application exponentielle
(A.4.1) est l'exponentielle, de $A\otimes g$
dans $V(A)\Subset
A\widehat{\otimes}\scrU^{\uphat}\bfv$.
La loi de groupe est le produit dans $\scrU^{\uphat}
\bfv$.

\subhead
A.7
\endsubhead
Un sch\'ema en groupe affine $G$ est la
limite projective de ses quotients $G_\alpha$
qui sont des groupes alg\'ebriques.
Il y a lieu de d\'efinir son alg\`ebre de Lie
$\Lie(G)$ comme \'etant la pro-alg\`ebre de
Lie limite projective des $\Lie(G_\alpha)$.
De fa\c{c}on \'equivalente, c'est la limite
projective des $\Lie(G_\alpha)$, qui sont de
dimension finie, munie de sa topologie de
limite projective d'espaces discrets.
C'est le dual de la limite inductive
$\colim(\Lie(G_\alpha)\upvee)$, muni de sa
topologie faible.

\subhead
A.8
\endsubhead
Un sch\'ema en groupe affine $G$ est {\it
pro-unipotent} s'il est limite projective de
groupes alg\'ebriques unipotents, i.e. si les
$G_\alpha$ de A.7 sont unipotents.
On d\'efinit l'alg\`ebre enveloppante
compl\'et\'ee $\scrU^{\uphat} \Lie(G)$ de
$\Lie(G)$ par
$$
\scrU^{\uphat}\Lie\,G=\lim\nolimits_\alpha
\scrU^{\uphat}\Lie(G_\alpha)=\lim\nolimits_{\alpha,n}
\scrU(\Lie(G_\alpha))/I^n,
$$
cette limite \'etant munie de sa topologie de
limite projective.
Avec ces d\'efinitions, A.6 reste valable,
comme on le v\'erifie par passage aux limites.

\subhead
A.9
\endsubhead
Cas particulier. Soit $\Lib((x_i)_{i\in I})$
l'alg\`ebre de Lie libre sur $K$ engendr\'ee
par une famille finie de g\'en\'erateurs
$x_i$.
Notons-la simplement $\Lib$.
Son alg\`ebre enveloppante est l'alg\`ebre
associative libre $K\left<(x_i)_{i\in I}\right>$.

Soit $\Lib^{\uphat}$ le compl\'et\'e de Lib pour
la filtration centrale descendante.
C'est l'alg\`ebre de Lie du sch\'ema en groupe
pro-unipotent
$$
\exp(\Lib^{\uphat})=\lim\,\exp(\Lib/\grZ^n).
$$
Son alg\`ebre enveloppante compl\'et\'ee est
l'alg\`ebre des s\'eries formelles
associatives $K\left<\!\left<(x_i)\right>\!\right>$.
Les points de $\exp(\Lib^{\uphat})$ \`a valeurs
dans une $K$-alg\`ebre $A$ sont les \'el\'ements
groupaux $g$ de 
$K\left<\!\left<(x_i)\right>\!\right>$:
terme constant $1$ et $\Delta g=g\otimes g$,
$\Delta$ \'etant le coproduit $K\left<\!\left<
(x_i)\right>\!\right>\to K\left<\!\left<
(x_i)\right>\!\right>\widehat{\otimes}K\left<\!\left<
(x_i)\right>\!\right>$ tel que $\Delta
x_i=x_i\otimes 1+1\otimes x_i$.

\subhead
A.10
\endsubhead
Soit $\bfv$ l'alg\`ebre de Lie d'un groupe
alg\'ebrique nilpotent $V$, et soit
$D$ une filtration
d\'ecroissante finie de l'alg\`ebre de Lie
$\bfv$, telle que $D^0(\bfv)=\bfv$.
On notera encore $D$ la filtration de $V$ par
les sous-groupes d'alg\`ebres de Lie les
$D^p\bfv$.
Puisque $[D^p\bfv, D^q\bfv]\Subset
D^{p+q}\bfv$, ce sont des sous-groupes
distingu\'es, et pour $p\Ge 1$, $D^pV/D^{p+1}V$
est ab\'elien.

On notera encore $D$ les filtrations induites
sur $\scrU\bfv$, $\scrU^{\uphat} \bfv$ et son
dual $\scrO(V)$.
La filtration $n\mapsto D^{-n}\scrO(V)$
de $\scrO(V)$ est croissante, et
$D^1\scrO(V)=0$.

Pour $\bfv$ muni du crochet nul,
$\exp(\bfv,\text{crochet }0)$ est le sch\'ema
pro-vectoriel $\bfv$.

Par A.1, l'application exponentielle $\exp\colon\,
\bfv\to V$ respecte les filtrations $D$ des
alg\`ebres affines.
Par r\'eduction au cas lin\'eaire, on en
d\'eduit que l'alg\`ebre affine de
$V/D^{n+1}V$ est la sous-alg\`ebre de
$\scrO(V)$ engendr\'ee par $D^{-n}\scrO(V)$.
En particulier, $D^0\scrO(V)$ est l'alg\`ebre
affine de $V/D^1V$.

Si $M$ est une repr\'esentation de $V$, les
conditions suivantes sur une filtration $F$ de $M$
sont \'equivalentes:

\medskip\noindent
\underbar{$\,\,\,$}
compatiblit\'e \`a $D$ sur $\bfv$:
$D^p\bfv\otimes F^qM\to F^{p+q}M$,

\noindent
\underbar{$\,\,\,$}
compatibilit\'e \`a $D$ sur $\scrU\bfv$:
$D^p\scrU\bfv\otimes F^qM\to F^{p+q}M$,

\noindent
\underbar{$\,\,\,$}
compatibilit\'e \`a $D$ sur $\scrO(V)$:
$F^pM\to\sum D^{-r}\scrO(V)\otimes F^{p+r}M$

\medskip
La filtration $D$ de $\scrO(V)$ est stable par
translations \`a gauche et \`a droite.
Si $P$ est un $V$-torseur, elle fournit encore
une filtration $D$ sur l'alg\`ebre affine
$\scrO(P)$ de $P$.

\subhead
A.11 Fonctorialit\'es
\endsubhead
(i)\enspace
Soient $\bfv_1$ une sous-alg\`ebre de Lie de
$\bfv$ (resp. un quotient) et $V_1$ le
sous-groupe (resp. quotient) correspondant de
$V$.
Notons $D_1$ la filtration de $\bfv_1$ induite
par (resp. quotient de) la filtration $D$ de
$\bfv$, et les filtrations qui s'en
d\'eduisent sur $\Sym^*\bfv_1$, $\scrU\bfv_1$,
$\scrU^{\uphat}\bfv_1$ et $\scrO(V_1)$.
La puissance symm\'etrique $\Sym^n(\bfv_1)$
est un sous-espace (resp. quotient) de
$\Sym^n(\bfv)$, et sa filtration $D_1$ est
induite par (resp. quotient de) la filtration
$D$ de $\Sym^n(\bfv)$.
Il r\'esulte de A.1 que la m\^{e}me assertion
vaut pour $\scrU v_1$ et $\scrU^{\uphat}v_1$
et que dualement, la filtration $D_1$ du
quotient (resp. sous-espace) $\scrO(V_1)$ de
$\scrO(V)$ est quotient de (resp. induite
par) la filtration $D$ de $\scrO(V)$.

\noindent
(ii)\enspace
Supposons que la filtration $D$ de $\bfv$
provienne d'une graduation
$\bfv=\oplusop_{n\Ge 0}\bfv_n$.
Cette graduation se propage \`a
$\Sym^*(\bfv)$, $\scrU\bfv$,
$\scrU^{\uphat}\bfv$ et $\scrO(V)$.
Pour $\scrU^{\uphat}\bfv$, ``graduation'' et
\`a prendre au sens de ``d\'ecomposition en
produit'' plut\^{o}t que ``d\'ecomposition en
somme''.
La filtration $D$ de $\Sym^n(\bfv)$ provient
de la graduation de $\Sym^n(\bfv)$.
Par A.1, la m\^{e}me assertion vaut pour
$\scrU\bfv$, $\scrU^{\uphat}\bfv$ et
$\scrO(V)$.

\noindent
(iii)\enspace
Par passage \`a la limite, A.10 et A.11
restent valable dans le cas pro-unipotent.

\proclaim{A.12 Exemples} \rm
(i)\enspace
Prenons l'alg\`ebre de Lie pro-unipotente
$\Lib^{\uphat}$ de A.9, et la filtration centrale
descendante $\grZ$.
Elle provient de la graduation pour laquelle
les $x_i$ sont de degr\'e $1$.
Chaque mon\^{o}me $m$ en les $x_i$ d\'efinit
sur $\exp(\Lib^{\uphat})\Subset k\left<\!\left<
(x_i\right>\!\right>$ une fonction $c_m$
``coefficient de $m$ dans $g$'', et
l'alg\`ebre affine a pour base les $c_m$.
Le produit est le produit de m\'elange.
Le degr\'e de $c_m$ est l'oppos\'e de celui de
$m$.

\noindent
(ii)\enspace
Parmi les g\'en\'erateurs $x_i$, distinguons en
un, not\'e $x_0$, et posons $I=\{0\}\amalg J$.
On a 
$$
\Lib:=\Lib((x_i)_{i\in I})=k.x_0\ltimes
\Lib((\ad\,x_0)^n(x_j),\,\,
(n\Ge 0,\,j\in J)).
$$
Soit $D$ la filtration de Lib pour laquelle
$D^p$ est, pour $p\Ge 1$, la s\'erie centrale
descendante du second facteur.
La filtration $D$ provient de la graduation
pour laquelle $x_0$ est de degr\'e $0$ et les
$x_j$ $(j\in J)$ de degr\'e $1$.
Compl\'etant cette filtration, on obtient une
filtration, encore not\'ee $D$, de
$\Lib^{\uphat}$.
La filtration correspondante de l'alg\`ebre
affine
$$
\scrO(\exp(\Lib^{\uphat}))
$$
est d\'eduite de la graduation pour laquelle
le degr\'e de $c_m$ est l'oppos\'e du degr\'e
du mon\^{o}me $m$ en les $x_j$ $(j\in J)$.
\endproclaim

\subhead
A.13
\endsubhead
Soit $G$ un groupe alg\'ebrique produit
semi-direct du groupe multiplicatif $\dbG_m$
par un groupe alg\'ebrique unipotent $V$,
d'alg\`ebre de Lie $\bfv$.
L'action de $\dbG_m$ sur $\bfv$ gradue cette
alg\`ebre de Lie; sur $\bfv_n$,
$\lambda\in\dbG_m$ agit par multiplication par
$\lambda^n$.

Supposons $\bfv$ \`a degr\'es $>0$.
Le foncteur $\rho\mapsto d\rho$ est alors une
\'equivalence de la cat\'egorie des
repr\'esentations lin\'eaires de $G$ avec
celle des espaces vectoriels gradu\'es $E$,
muni d'une action compatible aux graduations
de l'alg\`ebre de Lie $\bfv$.
La graduation de $E$ donne l'action de
$\dbG_m$ et l'action de $\bfv$,
automatiquement nilpotente, donne celle de $V$.

Notons $k(a)$ l'espace vectoriel $k$
plac\'e en degr\'e $a$.
Pour $a<b$, une extension de $k(a)$ par
$k(b)$ dans la cat\'egorie des
repr\'esentations de $G$ est donn\'ee par un
homomorphisme $\bfv\to
\Hom(k(a),k(b))$, i.e.
par un \'el\'ement du dual de la partie
homog\`ene de degre\'e $b-a$ de $\bfv^{\ab}:=
\bfv/[\bfv,\bfv]$:
$$
\Ext_{\Rep(G)}^1(k(a),k(b))=
(\bfv^{\ab})_{b-a}{\upvee}.
\tag{A.13.1}
$$

\subhead
A.14
\endsubhead
Soit $G$ un sch\'ema en groupes produit
semi-direct de  $\dbG_m$ par $V$
pro-unipotent.
Ecrivons-le comme limite projective de
quotients par des sous-groupes distingu\'es
contenus dans $V$:
$$
G=\dbG_m\ltimes V=\lim\,\dbG_m\ltimes
V_\alpha.
$$

Supposons les alg\`ebres de Lie gradu\'ees
$\bfv_\alpha$ des $V_\alpha$ (A.13) \`a
degr\'e $>0$.
Passant \`a la limite, on d\'eduit de (A.13.1)
que
$$
\Ext_{\Rep(G)}^1(k(a),k(b))=
\colim(v_\alpha^{\ab})_{b-a}{\upvee}.
\tag{A.14.1}
$$

Supposons en outre que pour tout,
$$
\dim\,\Ext^1(k(0),k(n))<\infty.
$$
Les limites projectives (A.14.1) sont alors
stationnaires et, puisque si on rel\`eve dans les
$(\bfv_\alpha)_n$ des bases des
$(\bfv_\alpha^{\ab})_n$, on obtient un syst\`eme
g\'en\'erateur de $\bfv_\alpha$, la dimension de
$(\bfv_\alpha)_n$ est, pour chaque $n$, born\'ee
ind\'ependamment de $n$.
On notera $\bfv_{\GR}$ l'alg\`ebre de Lie
gradu\'ee de composantes homog\`enes les
limites projectives (stationnaires) des
composantes homog\`enes des $v_\alpha$.
Passant \`a la limite dans A.13, on obtient:

\proclaim{Proposition A.15}
Sous les hypoth\`eses de A.14:
$G=\dbG_m\ltimes V$, Lie $V$ \`a degr\'es
$>0$, $\Ext_{\Rep(G)}^1(k(0),k(n))$
de dimension finie, le foncteur $\rho\mapsto
d\rho$ est une \'equivalence de la cat\'egorie
des repr\'esentations de $G$ avec celle des
espaces vectoriels gradu\'es munis d'une 
action de l'alg\`ebre de Lie gradu\'ee
$\bfv_{\GR}$.
On a
$$
\Ext_{\Rep(G)}^1(k(0),k(n))=
(\bfv_{\GR}^{\ab})_n\upvee.
$$
\endproclaim

\newpage

\Heading{Bibliographie}{0}

{\nspace{
\ref{}
H. Bass.
{\it Generators and relations for cyclotomic
units}, Nagoya Math. J. {\bf 27} (1966),\break
p.~401--407.

\medskip
\ref{}
A. Beilinson, J. Bernstein, et P. Deligne.  
{\it Faisceaux pervers},  in: 
 Analyse et topologie sur
les espaces singuliers, Ast\'erisque 100,
SMF (1982).

\medskip
\ref{}
S. Bloch.
{\it Algebraic cycles and algebraic
$K$-theory},
Adv. in Math. {\bf 61}, 3 (1986), p.~267--304.

\medskip
\ref{}
\hbox to 35pt{\hrulefill}~.
{\it The moving lemma for higher Chow groups},
J. Alg. Geom. {\bf 3}, 3 (1994), p.~537--568.

\medskip
\ref{}
K. T. Chen.
{\it Iterated integrals of differential forms
and loop space homology}, Ann. of Math. {\bf
97} (1973), p.~217--246.

\medskip
\ref{}
\hbox to 35pt{\hrulefill}~.
{\it Reduced Bar Constructions on de Rham
complexes}, in: Algebra, Topology and Category
Theory, a collection of papers in honor of
Samuel Eilenberg, p.~19--32,
Academic Press (1976).

\medskip
\ref{}
P. Deligne.
{\it Le groupe fondamental de la droite
projective moins trois points} in:
Galois group over $\dbQ$, MSRI Publ. 16,
p.~79--313, Springer Verlag (1989).

\medskip
\ref{}
\hbox to 35pt{\hrulefill}~.
{\it Cat\'egories tannakiennes} in:
Grothendieck Festschrift vol. 2, p.~111--195.
Progress in Math. 87, Birkh\"auser (1990).

\medskip
\ref{}
M. Demazure et P. Gabriel.
{\it Groupes alg\'ebriques}, Masson (1970).

\medskip
\ref{}
A.~B.~Goncharov.
{\it Polylogarithms in arithmetic and
geometry}, Proc. ICM Zurich 1994, p.~374--387,
Birkhauser.

\medskip
\ref{}
A.~B.~Goncharov.
{\it The dihedral Lie algebras and Galois
symmetries of $\pi_1(\dbP^1-\{0,\infty\}\CUP
\mu_N)$}, Duke Math. J. {\bf 110} 3 (2001), p.~397--487.

\medskip
\ref{}
R. Hain and M. Matsumoto.
{\it Weighted completion of Galois groups and
Galois actions on the fundamental group of
$\dbP^1-\{0,1,\infty\}$}, preprint (2001).

\medskip
\ref{}
R. Hain and S. Zucker.
{\it Unipotent variations of mixed Hodge
structure}, Inv. Math. {\bf 88}, (1987)
p.~83--124.

\medskip
\ref{}
M. Hanamura.
{\it Mixed motives and algebraic cycles I.},
Math. Res. Lett. {\bf 2}, 6 (1995),
p.~811--821.
See also II, preprint (1996).

\medskip
\ref{}
A. Huber.
{\it Mixed Motives and their Realization in
Derived Categories}, Lecture Notes in Math.
1604, Springer Verlag (1995).

\medskip
\ref{}
\hbox to 35pt{\hrulefill}~.
{\it Realization of Voevodsky's motives}, J.
Alg. Geom. {\bf 9} (2000), p.~755--799.

\medskip
\ref{}
U. Jannsen.
{\it Mixed motives and algebraic $K$-theory},
Lecture Notes in Math. 1400, Springer Verlag
(1990).

\medskip
\ref{}
D. Kubert.
{\it The universal ordinary distribution},
Bull. SMF {\bf 107} (1979), p.~179--202.

\medskip
\ref{}
M. Levine.
{\it Tate motives and the vanishing
conjectures for algebraic $K$-theory}, in:
Algebraic $K$-theory and algebraic topology
(Lake Louise, 1991), p.~167--188.
NATO Adv. Sci. but. Ser. C Math. Phys. {\bf
407}, Kluwer (1993).

\medskip
\ref{}
\hbox to 35pt{\hrulefill}~.
{\it Bloch's higher Chow groups revisited}, in:
$K$-theory (Strasbourg 1992), p.~235--320,
Ast\'erisque 226, SMF (1994).

\medskip
\ref{}
\hbox to 35pt{\hrulefill}~.
{\it Mixed Motives}, Math. Surveys and
Monographs {\bf 57}, AMS (1998).

\medskip
\ref{}
Chr. Reutenauer.
{\it Free Lie algebras}, LMS monographs new ser.
{\bf 7},
Oxford University Press (1993).

\medskip
\ref{}
T. Terasoma.
{\it Multiple zeta values and mixed Tate
motives}.

\medskip
\ref{}
V. Voevodsky.
{\it Triangulated categories of motives over a
field}, in: Cycles, transfer and motivic
homology theories, p.~188--238,
Ann. of Math. Studies 143, Princeton University
Press (2000).

\medskip
\ref{}
L. C. Washington.
{\it Introduction to cyclotomic fields --
graduate texts in Math.} {\bf 81}, Springer
Verlage (1997).

\medskip
\ref{}
Z. Wojtkowiak.
{\it Cosimplicial objects in algebraic
geometry}, in Algebraic $K$-theory and
algebraic topology (Lake Louise, 1991), p.~287--327.
NATO Adv. Sci. but. Ser. C Math. Phys. {\bf
407}, Kluwer (1993).

}}

\enddocument